\documentclass[pdflatex]{article}

\usepackage{graphicx}%
\usepackage[backref=page]{hyperref}
\usepackage{multirow}%
\usepackage{amsmath,amssymb,amsfonts}%
\usepackage{amsthm}%
\usepackage{mathrsfs}%
\usepackage[title]{appendix}%
\usepackage{xcolor}%
\usepackage{textcomp}%
\usepackage{manyfoot}%
\usepackage{booktabs}%
\usepackage{algorithm}%
\usepackage{algorithmicx}%
\usepackage{algpseudocode}%
\usepackage{listings}%

\usepackage[utf8]{inputenc}
\usepackage[T1]{fontenc}
\usepackage{lmodern,xcolor}
\usepackage{microtype}
\numberwithin{equation}{section}

\hypersetup{breaklinks,hidelinks,bookmarksnumbered}
\usepackage[capitalize,noabbrev,nameinlink]{cleveref}

\usepackage{thmtools}
\declaretheorem[within=section]{theorem}
\declaretheorem[sibling=theorem]{proposition}

\declaretheorem[sibling=theorem]{corollary}
\declaretheorem[sibling=theorem]{lemma}
\declaretheorem[sibling=theorem,style=definition,qed=\qedsymbol]{definition}
\declaretheorem[sibling=theorem,style=definition,qed=\qedsymbol]{construction}

\declaretheorem[sibling=theorem,style=remark,qed=\qedsymbol]{example}
\declaretheorem[sibling=theorem,style=remark,qed=\qedsymbol]{explication}
\declaretheorem[sibling=theorem,style=remark,qed=\qedsymbol]{remark}

\usepackage{nicematrix}
\NiceMatrixOptions{cell-space-limits=3pt}

\RequirePackage{etoolbox,xparse}
\RequirePackage{amsmath,amssymb,amsthm,mathtools,dsfont,amsfonts,mathrsfs}
\RequirePackage{quiver}

\newcommand{\define}[1]{\textbf{#1}}

\newcommand{\To}{\Rightarrow}
\newcommand{\xto}[1]{\xrightarrow{#1}}
\newcommand{\xrightsquigarrow}[1]{\overset{#1}{\rightsquigarrow}}
\newcommand{\cat}[1]{\mathsf{#1}}
\newcommand{\op}{\mathrm{op}}
\newcommand{\Set}{\cat{S}\mathsf{et}}

\newcommand{\WalkingArrow}{\cat{Arr}}
\newcommand{\ob}{\mathrm{ob}}
\newcommand{\mor}{\mathrm{mor}}
\DeclareMathOperator{\id}{id}
\DeclareMathOperator{\dom}{dom}

\DeclareMathOperator{\Ob}{Ob}

\newcommand{\twocat}[1]{\mathbf{#1}}
\newcommand{\Cat}{\twocat{Cat}}
\newcommand{\CAT}{\twocat{CAT}}

\def\p#1{\mathrel{\ooalign{\hfil$\mapstochar\mkern 5mu$\hfil\cr$#1$}}}
\newcommand{\proto}{\p\rightarrow}
\newcommand{\xproto}[1]{\overset{#1}{\p\rightarrow}}
\newcommand{\proTo}{\p\Rightarrow}

\newcommand{\dbl}[1]{\ifstrequal{#1}{1}{\mathds{1}}{\mathbb{#1}}}
\newcommand{\dbllong}[2]{\ifstrequal{#1}{1}{\mathds{1}_{#2}}{\mathbb{#1}\mathsf{#2}}}
\newcommand{\Span}{\dbllong{S}{pan}}
\newcommand{\Prof}{\dbllong{P}{rof}}
\newcommand{\Prom}{\dbllong{P}{rom}}
\newcommand{\Lax}{\twocat{Lax}}
\newcommand{\LaxTrans}{\twocat{LaxTrans}}

\newcommand{\Cart}{\twocat{Cart}}
\newcommand{\WalkingOb}{\dbllong{O}{b}}
\newcommand{\WalkingLoose}{\dbllong{L}{oose}}
\newcommand{\LeftOb}{\vdash}
\newcommand{\RightOb}{\dashv}
\newcommand{\TopOb}{\top}
\newcommand{\BotOb}{\bot}
\newcommand{\WalkingTight}{\dbllong{T}{ight}}
\newcommand{\WalkingSquare}{\dbllong{S}{quare}}

\DeclareMathOperator{\Loose}{Loose}

\DeclareMathOperator{\Inst}{Inst}
\DeclareMathOperator{\Dopf}{Dopf}
\DeclareMathOperator{\CartInst}{Inst}
\DeclareMathOperator{\CartDopf}{Dopf}

\newcommand{\xinlinecell}[9]{%
  \begin{tikzcd}[ampersand replacement=\&, column sep=small, row sep=small, arrows={font=\scriptsize}]
    #1 \arrow[r, "#7","\shortmid"{marking},] \arrow[d, "#5"'] \& #2 \arrow[d, "#6"] \\
    #3 \arrow[r, "#8"',"\shortmid"{marking},] \& #4
    \arrow[from=1-1, to=2-2, phantom, "#9"]
  \end{tikzcd}%
}

\begin{document}

\title{Presheaves on lax double functors; \\ Or, instances of models of double theories}

\author{Kevin Carlson\thanks{Topos Institute (\href{mailto:kevin@topos.institute}{kevin@topos.institute})} \and Evan Patterson\thanks{Topos Institute (\href{mailto:evan@topos.institute}{evan@topos.institute})}}

\date{}

\maketitle

\begin{abstract}
We introduce a notion of (co)presheaf on a lax double functor $X$, which we generally call
an instance. In the terminology of double-categorical logic,
a lax double functor valued in sets, possibly preserving finite products,
is called a \emph{model} of a double (Lawvere) theory.
By varying the double theory, we uniformly define a well-behaved notion of instances of categories,
profunctors, monads, monoidal categories, multicategories, and more, and we recover for instance 
the multifunctors into the category of sets in the last example.
We show that instances of $X$ can be described either in terms of modules from the terminal 
model $I$ to $X,$ satisfying an additional condition on triviality of the left action, or as 
loose natural transformations from $I$ to $X.$ We propose a notion of discrete opfibration between models
of a double theory, establish a comprehensive factorization system, and prove an elements correspondence giving 
an equivalence between the category of instances of and the category of discrete opfibrations over 
a model $X.$ We describe properties of the resulting categories of instances, relying on a 
``collage'' construction which we characterize as a lax colimit of a model of a double theory. 
An appendix gives a detailed treatment of certain morphisms of lax functors relevant also for bicategory theory: 
(loose) transformations versus modules and modifications versus modulations.

\end{abstract}

\clearpage
\tableofcontents

\section{Introduction}

In his Yoneda theory for double categories \cite{yoneda-2011}, Robert Paré established that the correct notion of a (tight) \emph{copresheaf} on a small double category $\dbl{D}$ is a lax double functor $X:\dbl{D}\to \Span,$ where $\Span$ is the double category of sets, functions, and spans.
Lambert and Patterson \cite{cartesian-double-theories-2024} further developed the theory
 of span-valued lax functors on $\dbl{D},$ while generalizing to other ``double doctrines''.
 In the \emph{simple} doctrine, there is no additional structure on $\dbl{D}$, while in the \emph{cartesian} doctrine,
 $\dbl{D}$ is a cartesian double category (Definition \ref{def:cartesianDoubleCategory})
 and the lax double functors are required to preserve finite products.
 We follow Lambert and Patterson in our focus on these two doctrines, though many others may be conceived.

This program is an approach to double-categorical logic. To highlight the analogy to 1-dimensional categorical logic, 
we often refer to a small, strict double category in some fixed doctrine as a \emph{double (Lawvere) theory}, and to a 
lax double functor $X:\dbl{D}\to \Span$, respecting any structure native to the doctrine, as a \emph{model} of that theory. In the simplest possible case, of the simple double theory $\dbl{D}=\dbl{1}$ given by the terminal double category, 
a model is just a category. Moving towards our object of interest in this work, a loose morphism, or \emph{module}, $H:X\proto Y$ in the virtual double category of models of $\dbl{1}$
corresponds to a \emph{profunctor} $H:X^\op\times Y\to \Set.$ 
In particular, we can recover presheaves and copresheaves on categories as loose morphisms out of or into the terminal model $1.$ 

While a \emph{model} $X:\dbl{D}\to \Span$ is itself a kind of double copresheaf, 
we are primarily interested here in considering (co)presheaves \emph{on} a model, moving
a dimension up. As we just saw, the classical notion of presheaf on a 1-category does arise as a 
``presheaf on a presheaf on $\dbl{1}$'' in this way, in terms of modules between lax double functors. 
This question of ``presheaves on models'' is naturally interesting, and it might appear that we already have an answer:
simply define a presheaf on a model $X$ to be a module $X\proto 1$ into the terminal model. 
A key observation motivating this paper, though, is that this is not quite right.
When $X$ is a model of a double theory containing nontrivial loose morphisms, the model $1$,
though terminal in terms of tight morphisms of models, is nonetheless rich enough to \emph{act nontrivially} on the
left of a module. For example, as we shall show, when $X$ is a multicategory, viewed as a model of a certain 
cartesian double theory, a module $1\proto X$ gives a multifunctor
into the category of augmented simplicial sets, rather than into bare sets.

We must therefore investigate the correct notion of  (co)presheaf on a model of a double theory.
We will in fact focus on the copresheaves, calling them 
\emph{instances} of the model. We shall define a notion of instance which does trivialize the left action of the terminal 
model as desired, recovering for example the plain multifunctors into $\Set$ as instances of multicategories. 
We find that this trivialization of the left action in fact precisely reduces a module out of $1$ into a loose 
transformation out of $1$, thus producing a connection to the general theory of the higher-categorical structure
on the collection of lax functors between two fixed double categories.

One might have tried to define instances via the internal notion of discrete opfibration in the 2-category of 
models of a double theory, but this sometimes fails to recover the desired notion---for instance, internal discrete 
opfibrations of multicategories only admit lifts for unary arrows. We therefore also work to show that the instances of 
models admit a well-behaved Grothendieck construction. The result is an extension of the theory of discrete opfibrations and 
comprehensive factorization systems to categories of models of double theories, which recovers the appropriate notions 
when specialized to cases like multicategories that have already been studied.

\subsection*{Historical background}

The modern theory of weak double categories begins with Robert Paré and Marco
Grandis' investigation into limits in double categories \cite{limits-in-double-1999}.
This paper is seemingly the first to define a double category to be weak in \emph{only one direction},
though Verity had already proposed the less tractable doubly-weak
double categories in his 1993 thesis \cite{verity-thesis-1992}. Already in this 
first paper, the ``leitmotif'' of weak double categories is
``summarised as follows: arrows which are too relaxed (like
profunctors, spans, relations) or too strict (like adjunctions) to have limits, can be
studied in a (pseudo) double category, correlating them with more ordinary
(horizontal) arrows'' \cite[p.\ 165]{limits-in-double-1999}.
This prescient choice of
\emph{leitmotif} has been well borne out in the quarter century since, as consequences
of the interacting strict and weak structures in a double category have grown steadily and
now seem to be accelerating rapidly.

Along with his series in collaboration with Grandis, Paré has written numerous solo papers
on double categories, of which the seminal paper on Yoneda theory for double categories \cite{yoneda-2011} 
is a particular influence on this work.
Here, Paré first draws
attention to lax double functors $\dbl{D}^\op\to\Span$ into the double
category of sets as the correct notion of a \emph{double presheaf} over a double category $\dbl{D}$.
In order to give a Yoneda embedding $\dbl{D}\to \Lax(\dbl{D}^\op,\Span),$ he is then led to define a notion of
loose morphism for the putative double category of lax functors in the codomain, resulting in the notion of \emph{module} between
double presheaves. (The codomain of the Yoneda embedding, it should be mentioned, is merely a
\emph{virtual} double category in general, as modules do not always compose; Paré later
studied this problem in \cite{modules-2013}.)

In a different vein, Spivak and Kent introduced a new approach to categorical database theory \cite{ologs-2012}, viewing a
small category $C$ as an \emph{ontology log} or \emph{database schema}, and a concrete
\emph{database instance} as a copresheaf on $C.$ It is from this background that we draw the word ``instance'' for ``copresheaf.'' This idea has since been extended in a few directions within applied category theory, including to the \emph{algebraic databases} of Schultz et al
\cite{algebraicdatabases-2016} and the \emph{attributed $C$-sets} or
``acsets'' of Patterson et al \cite{acsets-2022}.\footnote{Both these generalizations are designed to give formalisms
that respect the difference between a database \emph{row}, which is really just an
element of a set---databases' semantics are not varied if their rows are permuted---and the
values in a database \emph{column}, which are real data elements that should not be permuted
by database morphisms, and also generally admit interesting and important algebraic structure.}

\subsection*{Motivation from software}

In 2024, the authors and other collaborators at Topos Institute began work on
\href{https://www.catcolab.org/help}{CatColab}, a web-based application for formal, compositional
scientific modeling based directly on Paré's span-valued lax double functors. CatColab implements the simplest kind of morphism between models of a double
theory. These are the strict, tight natural transformations, comprising the tight morphisms
in Paré's virtual double category $\Lax(\dbl{D}^\op,\Span).$ Such morphisms allow the user to study, for instance, feedback loops in causal loop diagrams, 
an important question for systems dynamicists, via an efficient and mathematically rigorous implementation 
\cite{reg-nets-2023,polarities-2025}. In terms of CatColab, the study of instances is aimed at incorporating 
categorical database theory as a third level of abstraction in the system, above models and theories themselves. 
The system already allows for implicit use of instances presented as model morphisms, 
which via the comprehensive factorization studied here can be used to generate an instance; thus this mathematical 
work has been needed to justify that implementation and understand its relationship to an ongoing implementation 
of instances directly, rather than through presentations.

\subsection*{Related work}

In general, the comprehensive factorization systems constructed here differ
from that given to an arbitrary
locally presentable 2-category in \cite{diagrammatic-equations-2024}. 
The latter factorization system was mentioned at least as early as 2012 
on Shulman's \href{https://ncatlab.org/michaelshulman/show/comprehensive%2Bfactorization}{nLab page};
it uses the representable discrete opfibrations, which are not the best choice for 2-categories such as 
that of multicategories, where the requirement of lifting only against parameterized unary morphisms 
means that the representable discrete opfibrations do not correspond to multifunctors into $\Set.$

For another related direction, Riehl and Verity \cite{Riehl_Verity_2022} develop a theory of
initial functors in a rather wide class of virtual equipments induced from so-called ``$\infty$-cosmoi.''
Initial functors of categories, famously, can be characterized not only in terms of orthogonality 
to discrete opfibrations but also in terms of commuting with taking limits. Riehl and Verity's
definition focuses on the latter point of view. As we have not yet considered the internal notion
of limits in the virtual equipment of models of a double theory, the relationship between the two
definitions and the extent of any overlap is not yet known.

Additional projects having some overlap with this one are Moeller and Vasilakopoulou's \cite{monoidal-grothendieck-2020},
which generalizes the standard Grothendieck construction only for monoidal categories, but covers general fibrations, rather than just discrete ones as here; and Cigoli, Mantovani, and Metere's \cite{cigoli-mantovani-metere-2021}, which studies
discrete opfibrations in a 2-category of fibrations over a fixed base.

In future work, we will explore \emph{modal} (virtual) double theories as a more
convenient foundation for encoding non-simple double Lawvere theories and 
their models into software than the finite product theories discussed in \cite{products-revisited-2024}
and even the cartesian theories discussed above. This work will aim to capture motivating 
examples such as the algebraic profunctors
of \cite{algebraicdatabases-2016}.

We now summarize the contents of the paper in more detail.

\subsection*{Contents of the paper}

We begin in \cref{sec:background} by providing mathematical background. This
section can be skimmed by the reader familiar with double category theory
through the notion of a module between lax double functors, though we offer new formulations
of some old definitions.

In \cref{sec:instances}, we define an \emph{instance} of a model of a double theory.
We show that the category of instances of any fixed model $X$ is, if the 
theory is simple, of presheaf type 
(\cref{prop:instances_presheaf_type}) and, if the theory is cartesian, algebraic 
(\cref{prop:cartesian_instances_algebraic}) by giving an explicit construction of a category
$\kappa(X)$, called the collage of $X,$ whose copresheaves are instances of $X.$ This construction also
enables us to construct an adjoint triple between categories of instances induced by
any morphism of models. We note that the category of instances of $X$ may be seen as the category
of lax loose transformations and modulations from the terminal model $I$ to $X,$ though we emphasize
the module point of view. We prove that $\kappa(X)$ is the lax colimit of $X,$ that is, the 
initial category admitting a lax natural transformation $X\to \Delta \kappa(X)$, 
in \cref{prop:collage_universal_property}.

Our main theorem (\cref{thm:instances-discrete-opfibrations}) arrives in \cref{sec:discrete_opfibrations},
where we define \emph{discrete opfibrations} between models of double theories. These are morphisms
enjoying a lifting property analogous to that of ordinary discrete opfibrations of categories, relativized
to every loose morphism in the theory. \cref{thm:instances-discrete-opfibrations} itself gives a model-of-elements 
correspondence between instances of a model and discrete opfibrations over that model 
generalizing the classical equivalence between copresheaves on a category $C$ and discrete opfibrations over $C.$ 
Using local presentability of categories of models, proven in Appendix \ref{app:local-presentability}, 
we also show that the category of models admits a comprehensive factorization system, introduced for 
plain categories by Street and Walters in \cite{street-walters-1973}. 
This result allows one to present an instance of a model $X$ in terms of an arbitrary morphism $D:J\to X$, by 
using the comprehensive factorization to generate from $D$ a discrete opfibration, equivalently an instance.
In the left class of the factorization system, we also find a notion of initial morphism of models. 
All these results apply in both the simple and the cartesian doctrine.

Appendix \ref{app:loose-transformations} may be of independent interest, comprising a self-contained treatment of the relationship 
between lax loose transformations and modules, as well as between the modifications and modulations 
between these objects. We generalize the key results of Cockett, Koslowski, Seely, and Wood in \cite{modules-2003}
to the double-categorical setting: most notably, loose transformations and modulations embed fully faithfully 
in modules and modulations (\cref{thm:loose-transformations-into-modules}). Appendix 
\ref{app:proof-of-discrete-opfibrations} contains the rather long proof of \cref{thm:instances-discrete-opfibrations}.

\paragraph*{Acknowledgments}

The first author thanks Jason Brown, David Jaz Myers, Nathanael Arkor, and David
Spivak for helpful conversations related to this work. In particular, Myers
suggested \cref{prop:instances-as-loose-transformations} during a visit to Topos
Oxford. Arkor deserves further thanks for creating the Quiver software tool,
which finally allows an author to typeset diagrams which were long not worth the 
labor of transfer from the chalkboard. The authors are also grateful to an anonymous 
referee for valuable suggestions that have greatly improved the paper's exposition, and 
to Bryce Clarke for pointing out \cref{rem:connection-to-double-limits}.

The authors acknowledge support from the Air Force Office of Scientific Research
(AFOSR) Young Investigator Program (YIP) through Award FA9550-23-1-0133.

\section{Background: standard definitions from double category theory}\label{sec:background}

Our double categories are pseudo by default. We refer to the weaker direction as ``loose''; 
the loose arrows are variously known in the literature as ``proarrows'', ``horizontal arrows'', or, 
unfortunately, ``vertical arrows''. We draw loose arrows horizontally. 
We will, however, strictify double categories (using \cite[Section
7.5]{limits-in-double-1999}) at will, frequently assuming away the unitors and associators in 
definitions which would be complicated by their presence.

We generally denote a double category by a blackboard bold letter such as $\dbl{D}.$
Its underlying graph of categories is $\dbl{D}_1\rightrightarrows \dbl{D}_0$, so that $\dbl{D}_0$ 
is the category of objects and tight arrows and $\dbl{D}_1$ is the category of loose arrows and cells. 
The loose identities on an object or tight arrow are notated as $\id_x$ or $\id_f,$ while the tight identities 
on an object or loose arrow are denoted as $1_x$ and $1_m$ respectively. The loose
composition is notated as $m\odot n$ for composable loose arrows $x\xproto{m} y \xproto{n} z$. 
We notate loose composition in diagrammatic order exclusively. Tight composition of arrows $x\xto{f} y\xto{g} z$
will generally be notated as $f\cdot g$ (note the diagrammatic order.)
The use of $1$ for tight identities will lead to our use of $I$ for terminal objects.

Our double functors are often lax. We recall the definition of a lax functor in full here,
as we shall need to refer to its components in detail. 

\begin{definition}[Lax functor] \label{def:lax-functor}
  A \define{lax double functor} $F: \dbl{D} \to \dbl{E}$ between double
  categories $\dbl{D}$ and $\dbl{E}$ consists of
  \begin{itemize}
    \item a morphism in $\cat{Graph}(\cat{Cat})$ as below:
    \[\begin{tikzcd}
      {\dbl{D}_1} & {\dbl{D}_0} \\
      {\dbl{E}_1} & {\dbl{E}_0}
      \arrow[shift left, from=1-1, to=1-2]
      \arrow[shift right, from=1-1, to=1-2]
      \arrow["{F_1}"', from=1-1, to=2-1]
      \arrow["{F_0}", from=1-2, to=2-2]
      \arrow[shift left, from=2-1, to=2-2]
      \arrow[shift right, from=2-1, to=2-2]
    \end{tikzcd}\]
    \item two 2-morphisms in $\twocat{Cat}/\dbl{E}_0\times \dbl{E}_0$ (thus, with globular components) as below:
    \[\begin{tikzcd}
      {\dbl{D}_1\times_{\dbl{D}_0}\dbl{D}_1} & {\dbl{D}_1} && {\dbl{D}_1} & {\dbl{D}_0} \\
      && {\dbl{E}_0\times \dbl{E}_0} \\
      {\dbl{E}_1\times_{\dbl{E}_0}\dbl{E}_1} & {\dbl{E}_1} && {\dbl{E}_1} & {\dbl{E}_0}
      \arrow["\odot"', from=1-1, to=1-2]
      \arrow[""{name=0, anchor=center, inner sep=0}, "{F_1\times_{F_0}F_1}"', from=1-1, to=3-1]
      \arrow[from=1-2, to=2-3]
      \arrow[""{name=1, anchor=center, inner sep=0}, "{F_1}"{description}, from=1-2, to=3-2]
      \arrow[from=1-4, to=2-3]
      \arrow[""{name=2, anchor=center, inner sep=0}, "{F_1}"{description}, from=1-4, to=3-4]
      \arrow["\id", from=1-5, to=1-4]
      \arrow[""{name=3, anchor=center, inner sep=0}, "{F_0}", from=1-5, to=3-5]
      \arrow["\odot"', from=3-1, to=3-2]
      \arrow[from=3-2, to=2-3]
      \arrow[from=3-4, to=2-3]
      \arrow["\id", from=3-5, to=3-4]
      \arrow["{F_{-,-}}"', between={0.2}{0.8}, Rightarrow, from=0, to=1]
      \arrow["{F_{-}}", between={0.2}{0.8}, Rightarrow, from=3, to=2]
    \end{tikzcd}\]
    The component $F_{m,n}:F(m)\odot F(n)\to F(m\odot n)$ is called
    the \define{laxator} at $m$ and $n$ 
    and the component $F_x:\id_{Fx}\to F\id_x$ is called the \define{unitor} at $x$.
  \end{itemize}
  The following axioms must be satisfied.
  \begin{itemize}
    \item Associativity: we have the equality of 2-cells 
    \[\begin{tikzcd}[column sep=tiny]
      {\dbl{D}_1\times_{\dbl{D}_0}\dbl{D}_1\times_{\dbl{D}_0}\dbl{D}_1} & {\dbl{E}_1\times_{\dbl{E}_0} \dbl{E}_1\times_{\dbl{E}_0} \dbl{E}_1} && {\dbl{D}_1\times_{\dbl{D}_0}\dbl{D}_1\times_{\dbl{D}_0}\dbl{D}_1} & {\dbl{E}_1\times_{\dbl{E}_0} \dbl{E}_1\times_{\dbl{E}_0} \dbl{E}_1} \\
      {\dbl{D}_1\times_{\dbl{D}_0}\dbl{D}_1} & {\dbl{E}_1\times_{\dbl{E}_0} \dbl{E}_1} & {=} & {\dbl{D}_1\times_{\dbl{D}_0}\dbl{D}_1} & {\dbl{E}_1\times_{\dbl{E}_0} \dbl{E}_1} \\
      {\dbl{D}_1} & {\dbl{E}_1} && {\dbl{D}_1} & {\dbl{E}_1}
      \arrow[""{name=0, anchor=center, inner sep=0}, from=1-1, to=1-2]
      \arrow["{{\odot \times_{\dbl{D}_0} \dbl{D}_1}}"', from=1-1, to=2-1]
      \arrow["{{\odot\times_{\dbl{E}_0} \dbl{E}_1}}", from=1-2, to=2-2]
      \arrow[""{name=1, anchor=center, inner sep=0}, from=1-4, to=1-5]
      \arrow["{{\dbl{D}_1\times_{\dbl{D}_0}\odot}}"', from=1-4, to=2-4]
      \arrow["{{\dbl{E}_1\times_{\dbl{E}_0}\odot}}", from=1-5, to=2-5]
      \arrow[""{name=2, anchor=center, inner sep=0}, from=2-1, to=2-2]
      \arrow["\odot"{description}, from=2-1, to=3-1]
      \arrow["\odot"{description}, from=2-2, to=3-2]
      \arrow[""{name=3, anchor=center, inner sep=0}, from=2-4, to=2-5]
      \arrow["\odot"{description}, from=2-4, to=3-4]
      \arrow["\odot"{description}, from=2-5, to=3-5]
      \arrow[""{name=4, anchor=center, inner sep=0}, "{{F_1}}"{description}, from=3-1, to=3-2]
      \arrow[""{name=5, anchor=center, inner sep=0}, "{{F_1}}"{description}, from=3-4, to=3-5]
      \arrow["{{F_{-,-}\times_{F_0}F_1}}"{description}, draw=none, from=0, to=2]
      \arrow["{{F_1\times_{F_0}F_{-,-}}}"{description}, draw=none, from=1, to=3]
      \arrow["{{F_{-,-}}}"{description}, draw=none, from=2, to=4]
      \arrow["{{F_{-,-}}}"{description}, draw=none, from=3, to=5]
\end{tikzcd}\]
    so that there is a canonical ternary 
    laxator with components of the form $F_{l,m,n}:F(l)\odot F(m)\odot F(n)\to F(l\odot m\odot n).$
    (Note that here, as in many places below, we tacitly strictify $\dbl{D}$ and $\dbl{E}$ to avoid mentioning
    associators.) 
    \item Unitality: 
    We have the equalities of 2-cells: 
    \[\begin{tikzcd}[column sep=small]
      {\dbl{D}_0\times_{\dbl{D}_0}\dbl{D}_1} & {\dbl{E}_0\times_{\dbl{E}_0} \dbl{E}_1} && {\dbl{D}_1} && {\dbl{D}_1\times_{\dbl{D}_0}\dbl{D}_0} & {\dbl{E}_0\times_{\dbl{E}_0} \dbl{E}_1} \\
      {\dbl{D}_1\times_{\dbl{D}_0} \dbl{D}_1} & {\dbl{E}_1\times_{\dbl{E}_0} \dbl{E}_1} & {=} && {=} & {\dbl{D}_1\times_{\dbl{D}_0} \dbl{D}_1} & {\dbl{E}_1\times_{\dbl{E}_0} \dbl{E}_1} \\
      {\dbl{D}_1} & {\dbl{E}_1} && {\dbl{E}_1} && {\dbl{D}_1} & {\dbl{E}_1}
      \arrow[""{name=0, anchor=center, inner sep=0}, from=1-1, to=1-2]
      \arrow["{{\id_{-}\times_{\dbl{D}_0} \dbl{D}_1}}"', from=1-1, to=2-1]
      \arrow["{{\id_{-}\times_{\dbl{E}_0} \dbl{E}_1}}", from=1-2, to=2-2]
      \arrow[""{name=1, anchor=center, inner sep=0}, "{{F_1}}"{description}, shift right=4, from=1-4, to=3-4]
      \arrow[""{name=2, anchor=center, inner sep=0}, "{{F_1}}"{description}, shift left=4, from=1-4, to=3-4]
      \arrow[""{name=3, anchor=center, inner sep=0}, from=1-6, to=1-7]
      \arrow["{{ \dbl{D}_1\times_{\dbl{D}_0}\id_{-}}}"', from=1-6, to=2-6]
      \arrow["{{\dbl{E}_1\times_{\dbl{E}_0} \id_{-}}}", from=1-7, to=2-7]
      \arrow[""{name=4, anchor=center, inner sep=0}, from=2-1, to=2-2]
      \arrow["\odot"{description}, from=2-1, to=3-1]
      \arrow["\odot"{description}, from=2-2, to=3-2]
      \arrow[""{name=5, anchor=center, inner sep=0}, from=2-6, to=2-7]
      \arrow["\odot"{description}, from=2-6, to=3-6]
      \arrow["\odot"{description}, from=2-7, to=3-7]
      \arrow[""{name=6, anchor=center, inner sep=0}, "{{F_1}}"{description}, from=3-1, to=3-2]
      \arrow[""{name=7, anchor=center, inner sep=0}, "{{F_1}}"{description}, from=3-6, to=3-7]
      \arrow["{{F_{-}\times_{F_0} F_1}}"{description}, draw=none, from=0, to=4]
      \arrow[between={0.3}{0.7}, equals, from=1, to=2]
      \arrow["{{F_1\times_{F_0} F_{-}}}"{description}, draw=none, from=3, to=5]
      \arrow["{{F_{-,-}}}"{description}, draw=none, from=4, to=6]
      \arrow["{{F_{-,-}}}"{description}, draw=none, from=5, to=7]
    \end{tikzcd}\]
    in $\dbl{D}$ and $\dbl{E},$ stating that laxators act trivially at unitors on either side.
  \end{itemize}
  If the laxators and unitors are isomorphisms in $\dbl{E}_1$, the double
  functor is called \define{pseudo}; if they are identities, the double functor
  is \define{strict}. If just the unitors are invertible, then $F$ is said to be
  \define{normal}.
\end{definition}

We now turn to discussion of the morphisms and higher morphisms between lax double functors. 
For this it is expedient to introduce notation for a few of the smallest and most important double categories.

\begin{definition}[Walking cells]
We define the following double categories:
\begin{itemize}
\item Let $\WalkingOb$ denote the terminal double category, also known as the walking object: thus $\WalkingOb$ has a single object,
equipped only with its identity morphisms and their identity cell.
\item Let $\WalkingTight = \{\TopOb\xto{t}\BotOb\}$ denote the walking tight arrow, i.e., the double category freely generated by two objects and a single tight morphism between them.
\item Similarly, let $\WalkingLoose = \{\LeftOb\xproto{\ell}\RightOb\}$ denote the walking loose arrow.
\item Finally, let $\WalkingSquare$ be the walking square double category, as shown below:
\[\begin{tikzcd}
	\ulcorner & \urcorner \\
	\llcorner & \lrcorner
	\arrow[""{name=0, anchor=center, inner sep=0}, "{\ell_\TopOb}"{inner sep=.8ex}, "\shortmid"{marking}, from=1-1, to=1-2]
	\arrow["{t_\LeftOb}"', from=1-1, to=2-1]
	\arrow["{t_\RightOb}", from=1-2, to=2-2]
	\arrow[""{name=1, anchor=center, inner sep=0}, "{\ell_\BotOb}"'{inner sep=.8ex}, "\shortmid"{marking}, from=2-1, to=2-2]
	\arrow["\sigma"{description}, draw=none, from=0, to=1]
\end{tikzcd} \qedhere\]
\end{itemize}
\end{definition}

In this paper, we privilege the \emph{tight, strict} transformations as our
default notion of natural transformation between lax double functors, though almost all 
kinds of such morphisms will make cameo appearances: we shall make use of 
tight, \emph{lax} transformations above \cref{prop:collage_universal_property}, and also of 
\emph{loose} transformations in Appendix \ref{app:loose-transformations}. We define this 
first most important kind of 2-morphism now.

\begin{definition}[2-category of double categories and lax functors]\label{def:morphism-lax-functors}
Given lax double functors $F,G:\dbl{D}\to \dbl{E}$, a \define{morphism} $\alpha:F\to G$ is a 
lax double functor $\alpha:\dbl{D}\times \WalkingTight\to \dbl{E}$ that restricts to
$F$ and $G$ on $\dbl{D}\times \{\TopOb\}$ and $\dbl{D}\times \{\BotOb\}$, respectively.

We denote the 2-category of double categories, lax double functors, and 
morphisms between them by $\Lax.$
\end{definition}

Unpacking this definition, such a morphism between lax double functors
consists of a 2-morphism in $\cat{Graph}(\cat{Cat})$,
\[\begin{tikzcd}[column sep=large]
	{\dbl{D}_1} & {\dbl{D}_0} \\
	{\dbl{E}_1} & {\dbl{E}_0}
	\arrow[shift left, from=1-1, to=1-2]
	\arrow[shift right, from=1-1, to=1-2]
	\arrow[""{name=0, anchor=center, inner sep=0}, "{F_1}"', curve={height=6pt}, from=1-1, to=2-1]
	\arrow[""{name=1, anchor=center, inner sep=0}, "{G_1}", curve={height=-6pt}, from=1-1, to=2-1]
	\arrow[""{name=2, anchor=center, inner sep=0}, "{F_0}"', curve={height=6pt}, from=1-2, to=2-2]
	\arrow[""{name=3, anchor=center, inner sep=0}, "{G_0}", curve={height=-6pt}, from=1-2, to=2-2]
	\arrow[shift right, from=2-1, to=2-2]
	\arrow[shift left, from=2-1, to=2-2]
	\arrow["{\alpha_1}", between={0.2}{0.8}, Rightarrow, from=0, to=1]
	\arrow["{\alpha_0}", between={0.2}{0.8}, Rightarrow, from=2, to=3]
\end{tikzcd},\]
such that the following two squares commute in the categories indicated below each square:
\[\begin{tikzcd}
	{\mathrm{id}_{F_0(-)}} & {F_1(\mathrm{id}_{-})} & {F_1\odot F_1} & {F_1(-\odot -)} \\
	{\mathrm{id}_{G_0(-)}} & {G_1(\mathrm{id}_{-})} & {G_1\odot G_1} & {G_1(-\odot -)} \\
	{\dbl{D}_0} & {\dbl{E}_1} & {\dbl{D}_1\times_{\dbl{D}_0}\dbl{D}_1} & {\dbl{E}_1}
	\arrow["{F_{-}}", from=1-1, to=1-2]
	\arrow["{\mathrm{id}_{\alpha_{-}}}"', from=1-1, to=2-1]
	\arrow["{\alpha_{\mathrm{id}_{-}}}", from=1-2, to=2-2]
	\arrow["{F_{-,-}}", from=1-3, to=1-4]
	\arrow["{\alpha\odot \alpha}"', from=1-3, to=2-3]
	\arrow["{\alpha_{-\odot -}}", from=1-4, to=2-4]
	\arrow["{G_{-}}"', from=2-1, to=2-2]
	\arrow["{G_{-,-}}"', from=2-3, to=2-4]
	\arrow[from=3-1, to=3-2]
	\arrow[from=3-3, to=3-4]
\end{tikzcd}\]

Recall that modules, as defined in \cite{yoneda-2011}
and \cite[Definition 9.1]{cartesian-double-theories-2024}, are the loose morphisms
between lax double functors. Specifically:
\begin{definition}[Module between double functors]\label{def:module}
Given lax double functors $F,G:\dbl{D}\to \dbl{E},$ a \define{module} 
$F\proto G$ is a lax functor $M:\dbl{D}\times\WalkingLoose\to \dbl{E}$ that
restricts to $F$ and $G$ on the endpoints.
\end{definition}

The \emph{components} of the module, as listed in \cite[Definition 9.1]{cartesian-double-theories-2024},
are obtained by restricting the functor $M_1:\dbl{D}_1\times \WalkingLoose_1\to \dbl{E}_1$ along the inclusion
$\{\ell\}\to \WalkingLoose_1.$ Thus $M$ extends $F,G$ with a new functor, also denoted $M,$
of signature $\dbl{D}_1\to \dbl{E}_1,$ sending each loose arrow $m:x\proto y$ in $\dbl{D}$ to a loose arrow
$M(m):F(x)\proto G(y)$ in $\dbl{E}.$ The laxators and unitors of $M$ provide
actions of $F$ and $G$ on this $M.$ Thus:

\begin{explication}[Module between double functors, unpacked]\label{explication:module_unpacked}
Given lax double functors $F,G:\dbl{D}\to \dbl{E},$ a module $M:F\proto G$ consists of the following data:
\begin{itemize}
\item A functor $M:\dbl{D}_1\to \dbl{E}_1,$ sending each loose arrow $m:x\proto y$ in $\dbl{D}$ to a loose arrow
$M(m):F(x)\proto G(y)$ in $\dbl{E},$ called the \define{component} of $M$ at $m,$ and similarly 
for cells. 
\item A globular natural transformation 
$M^\ell:F_1\odot M\To M(-\odot -):\dbl{D}_1\times_{\dbl{D}_0}\dbl{D}_1\to \dbl{E}_1,$ 
and similarly $M^r:M\odot G_1\To M(-\odot -),$ called the \define{left} and \define{right actions} of $M.$
\item The actions are associative and unital with respect to loose composition, in that the following diagrams commute: 
\[\begin{tikzcd}
	{F_1\odot F_1\odot M} & {F_1\odot M(-\odot -)} &&& \\
	{F(-\odot-)\odot M} & {M(-\odot -\odot -)} && {\mathrm{id}\odot M} & {F_1\odot M} \\
	{M\odot G_1\odot G_1} & {M(-\odot -)\odot G_1} &&& M \\
	{M\odot G(-\odot-)} & {M(-\odot -\odot -)} && {M\odot \mathrm{id}} & {M\odot G_1} \\
	{F_1\odot M\odot G_1} & {F_1\odot M} &&& M \\
	{M\odot G_1} & M && {} & {} \\
	{} & {}
	\arrow["{{1_F\odot M^\ell}}", from=1-1, to=1-2]
	\arrow["{{F_{-,-}\odot 1_M}}"', from=1-1, to=2-1]
	\arrow["{{M^\ell_{-,-\odot -}}}", from=1-2, to=2-2]
	\arrow["{{M^\ell_{-\odot -,-}}}"', from=2-1, to=2-2]
	\arrow["{{F_{-}\odot 1_M}}", from=2-4, to=2-5]
	\arrow[equals, from=2-4, to=3-5]
	\arrow["{{M^\ell}}", from=2-5, to=3-5]
	\arrow["{{M^r\odot 1_G}}", from=3-1, to=3-2]
	\arrow["{{1_M\odot G_{-,-}}}"', from=3-1, to=4-1]
	\arrow["{{M^r_{-\odot -,-}}}", from=3-2, to=4-2]
	\arrow["{{M^r_{-,-\odot -}}}"', from=4-1, to=4-2]
	\arrow["{{1_M\odot G_{-}}}", from=4-4, to=4-5]
	\arrow[equals, from=4-4, to=5-5]
	\arrow["{{M^r}}", from=4-5, to=5-5]
	\arrow["{F_1\odot M^r}"{inner sep=.8ex}, "\shortmid"{marking}, from=5-1, to=5-2]
	\arrow["{M^\ell\odot G_1}"'{inner sep=.8ex}, "\shortmid"{marking}, from=5-1, to=6-1]
	\arrow["{M^\ell}"{inner sep=.8ex}, "\shortmid"{marking}, from=5-2, to=6-2]
	\arrow["{M^r}"'{inner sep=.8ex}, "\shortmid"{marking}, from=6-1, to=6-2]
	\arrow["{{\dbl{D}_0\times_{\dbl{D}_0} \dbl{D}_1\cong \dbl{D}_1\to \dbl{E}_1}}"{description}, draw=none, from=6-4, to=6-5]
	\arrow["{{\dbl{D}_1\times_{\dbl{D}_0}\dbl{D}_1\times_{\dbl{D}_0} \dbl{D}_1\to \dbl{E}_1}}"{description}, draw=none, from=7-1, to=7-2]
\end{tikzcd}\]
\end{itemize}
\end{explication}

We recall, finally, the notion of cell among morphisms of lax double functors\footnote{We do not define 
multimodulations as they will not arise in this paper.}:

\begin{definition}[Modulation]\label{def:modulation}
Given double categories $\dbl{D}$ and $\dbl{E},$ lax functors
$F_\ulcorner,F_\urcorner,F_\llcorner,F_\lrcorner:\dbl{D}\to \dbl{E},$
tight transformations $\alpha_\LeftOb:F_\ulcorner\to F_\llcorner$ and 
$\alpha_\RightOb:F_\urcorner\to F_\lrcorner,$ and modules 
$M_\TopOb:F_\ulcorner\proto F_\urcorner$ and $M_\BotOb:F_\llcorner\proto F_\lrcorner,$
a \define{modulation} $\Theta$ with boundary $(\alpha_\LeftOb,\alpha_\RightOb,M_\TopOb,M_\BotOb)$
is a lax double functor $\Theta:\dbl{D}\times \WalkingSquare\to \dbl{E},$ restricting to the
given data on the boundary of $\WalkingSquare.$
\end{definition}

By restricting $\Theta$ along the inclusion of $\{\sigma\}$, we find that the only new
data in a modulation relative to its boundary is the choice 
$\Theta(\gamma,\sigma)$ for every cell $\xinlinecell{x}{y}{z}{w}{f}{g}{m}{n}{\gamma}$ in $\dbl{D}.$ 
Now, $(\gamma,\sigma)=(1_m,\sigma)\cdot (\gamma,1_{\ell_\BotOb})=(\gamma,1_{\ell_{\TopOb}})\cdot (1_n,\sigma),$ 
so the functoriality of $\Theta_1$ implies that these values are determined (subject to the conditions) by the cells 
$\Theta m:=\Theta(1_m,\sigma).$ This leads to the unpacked definition of modulation:

\begin{explication}[Modulation, unpacked]\label{explication:modulation_unpacked}
A modulation with boundary $(\alpha_\LeftOb,\alpha_\RightOb,M_\TopOb,M_\BotOb)$ as defined above
consists of a natural transformation $\Theta:M_\TopOb \To M_\BotOb: \dbl{D}_0\to \dbl{E}_1$
over $\alpha_\LeftOb$ and $\alpha_\RightOb,$
such that the equivariance diagrams below commute in the functor category $\Cat(\dbl{D}_1\times_{\dbl{D}_0}\dbl{D}_1,\dbl{E}_1)$:
\[\begin{tikzcd}
	{M_\top \odot F_{\urcorner}} & {M_\bot\odot F_{\lrcorner}} & {F_{\ulcorner}\odot M_\top} & {F_{\llcorner}\odot M_\bot} \\
	{M_\top(-\odot-)} & {M_\bot(-\odot -)} & {M_\top(-\odot-)} & {M_\bot(-\odot -)}
	\arrow["{\Theta\odot \alpha_{\dashv}}", from=1-1, to=1-2]
	\arrow["{M_\top^r}"', from=1-1, to=2-1]
	\arrow["{M_\bot^r}", from=1-2, to=2-2]
	\arrow["{\alpha_{\vdash}\odot \Theta}", from=1-3, to=1-4]
	\arrow["{M_\top^\ell}"', from=1-3, to=2-3]
	\arrow["{M_\bot^\ell}", from=1-4, to=2-4]
	\arrow["{\Theta_{-\odot -}}"', from=2-1, to=2-2]
	\arrow["{\Theta_{-\odot -}}"', from=2-3, to=2-4]
\end{tikzcd}\]

\end{explication}

We shall need the definition of an equipment on one occasion:

\begin{definition}[Equipment]\label{def:equipment}
  An \define{equipment} (also known as a \define{proarrow equipment} or a \define{fibrant double category}) is a double category $\dbl{D}$ such that
  the pairing $\dbl{D}_1\rightarrow \dbl{D}_0^2$ of the source and target functors is a fibration.
\end{definition}

For explication on the subject of equipments, see \cite{relative-monadicity-2025}
and \cite{generalized-multicategories-2010}.

Tabulators will also make a brief appearance:

\begin{definition}[Tabulator]\label{def:tabulator}
If $m:x\proto y$ is a loose morphism in a double category $\dbl{E},$ a \define{tabulator} of $m$ is an 
object $\top m$ equipped with a universal cell of the form below:
\[\begin{tikzcd}
	{\top m} && {\top m} \\
	x && y
	\arrow[""{name=0, anchor=center, inner sep=0}, "\shortmid"{marking}, equals, from=1-1, to=1-3]
	\arrow["s"', from=1-1, to=2-1]
	\arrow["t", from=1-3, to=2-3]
	\arrow[""{name=1, anchor=center, inner sep=0}, "m"'{inner sep=.8ex}, "\shortmid"{marking}, from=2-1, to=2-3]
	\arrow["{\tau_m}"{description}, draw=none, from=0, to=1]
\end{tikzcd}\]
The universal property is that, for any cell
\[\begin{tikzcd}
	z && z \\
	x && y
	\arrow[""{name=0, anchor=center, inner sep=0}, "\shortmid"{marking}, equals, from=1-1, to=1-3]
	\arrow["f"', from=1-1, to=2-1]
	\arrow["g", from=1-3, to=2-3]
	\arrow[""{name=1, anchor=center, inner sep=0}, "m"'{inner sep=.8ex}, "\shortmid"{marking}, from=2-1, to=2-3]
	\arrow["\alpha"{description}, draw=none, from=0, to=1]
\end{tikzcd}\]
of the same form, there is a unique tight morphism $h:z\to \top m$ such that
$h\cdot s=f$, $h\cdot t=g$, and $\alpha = \id_h\cdot \tau_m$.
\end{definition}

For example, the tabulator of a span $X\leftarrow A\to Y$ is just $A$ itself.

We shall also need the notion of cartesian double category \cite{aleiferi2018}. In an unavoidable difficulty
of terminology, we note that this is not the same notion as a double category ``with finite products,''
for which see \cite{products-revisited-2024}.

\begin{definition}[Cartesian double category]\label{def:cartesianDoubleCategory}
A double category $\dbl{D}$ is \define{cartesian} if the canonical double functors
$!:\dbl{D}\to \dbl{1}$ and $\Delta:\dbl{D}\to \dbl{D}\times \dbl{D}$ 
admit pseudo right adjoints. (It is \define{precartesian} if these adjoints exist but are only lax.)

In more elementary terms, $\dbl{D}$ is cartesian if the underlying categories $\dbl{D}_0$ and $\dbl{D}_1$ have finite products,
and the source, target, identity-assignment, and composition functors preserve them. 

A lax functor $F:\dbl{D}\to \dbl{E}$ between cartesian double categories is
\define{cartesian} if the underlying functors $F_0: \dbl{D}_0 \to \dbl{E}_1$ and
$F_1: \dbl{D}_1 \to \dbl{E}_1$ both preserve finite products.
The tight transformations
between cartesian lax functors are no different than in the simple case, while 
a module $M:F\proTo G:\dbl{D}\to \dbl{E}$ between cartesian lax functors is cartesian if the functor 
$M:\dbl{D}_1\to \dbl{E}_1$ preserves finite products. Modulations between cartesian modules 
have no further constraints over general modulations.
\end{definition}

We will exhibit our constructions below especially on two extreme classes of double categories,
those discrete in one direction or the other.

\begin{definition}[Discrete double categories]\label{def:discreteDoubleCategory}
There are two functors $\dbl{T},\dbl{L}:\cat{Cat}\to \cat{Dbl}$ embedding categories as double categories:
\begin{itemize}
\item The \define{tight embedding} $\dbl{T}C$ of a category $C$ is the double category with 
objects and tight morphisms coming from $C,$ with only identity loose morphisms and cells.
\item The \define{loose embedding} $\dbl{L}C$ of a category $C$ is the double category with
objects and loose morphisms coming from $C,$ with only identity tight morphisms and cells.
\end{itemize}

We call a double category \define{tightly discrete} if it is isomorphic to $\dbl{T}C$ for some category $C,$
and \define{loosely discrete} if it is isomorphic to $\dbl{L}C$ for some category $C.$
\end{definition}

Note that $\dbl{L}C$ is always a strict double category; it is arguably more natural to consider 
only the tight embedding of categories and the loose embedding of \emph{bi}categories, but we will 
not find use for the latter in our examples. 

Let us write $\Lax_n$ for the 2-category of double categories, \emph{normal} lax functors, and 
morphisms. A normal lax functor is a lax functor whose unitors are invertible. The discrete embeddings
send categories into double categories with normal lax functors between them.

\begin{proposition}[Functors out of discrete double categories]\label{prop:functors-out-of-discrete}
The 2-functor $\dbl{T}:\Cat\to \Lax_n$ is left biadjoint to the functor 
sending a double category $\dbl{D}$ to its category of objects and tight morphisms $\dbl{D}_0.$ 
\end{proposition}
\begin{proof}
Every normal lax functor is canonically isomorphic, using the unitors themselves, to one 
which strictly preserves identities, and it is easy to see that such strictly unitary 
double functors $\dbl{T}C\to \dbl{D}$ are uniquely determined by their tight part. 
\end{proof}
\section{Instances of models}\label{sec:instances}

We turn now to instances, the main subject of this paper. First, we recall a
terminological convention relied on especially in
\cite{cartesian-double-theories-2024}. The function of the terminology below is
to evoke the interpretation from categorical logic that lax double functors are
models of a theory, and more pragmatically, to allow us to use the same term
``model'' to encompass both lax functors out of ordinary double categories and
also cartesian lax functors out of cartesian double categories.

\begin{definition}[Double theories and their models]\label{def:theoryAndModel}
A \define{simple double theory} is a small, strict double category. A \define{cartesian double theory} 
is a small, strict double category equipped with cartesian structure. 
For the purposes of this paper, a \define{double theory} is either a simple or a cartesian double theory.

A \define{model} of a double theory $\dbl{D}$ is a lax double functor $X:\dbl{D}\to \Span,$
which preserves finite products in the cartesian case. We have no special logical terminology 
for modules between models of double theories, but when we use such modules between 
models of \emph{cartesian} double theories,
we shall always mean cartesian modules in the sense of \cref{def:cartesianDoubleCategory}.
\end{definition}

\begin{remark}[Strictness of double theories]
The strictness assumption on double theories is no real constraint due to the strictification theorem for
double categories \cite[Section 7.5]{limits-in-double-1999}. We make the restriction mainly to 
emphasize that we think of double theories as syntactic constructs. Regarding our 
running examples of tightly and loosely discrete double categories, we can remark that 
both cases produce simple double theories, as both are strict.
\end{remark}

We now aim to define a notion of
\emph{instance} of a model (\cref{def:theoryAndModel}) that is to arbitrary modules (\cref{def:module}) as, in the case of categories (models of
the terminal double theory), copresheaves are to arbitrary
profunctors. Indeed, over the terminal (simple) double theory $\WalkingOb$, we can recover 
functors $C\to \Set$ for a category $C:\WalkingOb\to\Span$ as 
modules $I\proto C,$ where $I$ is the terminal category.
Dually, presheaves on $C$ are modules $C \proto I$. 

In view of this motivating example, we might attempt to
define an instance of a model $X:\dbl{D}\to\Span$ as a module $I\proto X.$ 
However, an extra constraint is needed to obtain the correct definition for
models of an arbitrary double theory, whenever
the double theory contains nontrivial loose morphisms. (We shall illustrate the issue 
below when we arrive at examples.)
To solve this problem, we make the following definition, which will be explicated
and unfolded below.

\begin{definition}[Instance of a model, compressed definition]\label{def:compressed-instance}
  Let $\dbl{D}$ be a double theory and let $\dbl{E}$ be a double category with a
  terminal object $I$. Then $\dbl{D}$ has a terminal model $I$ in $\dbl{E}$
  sending every object, tight morphism, loose morphism, and cell in $\dbl{D}$ to
  $I$, $1_I$, $\id_I$, and $1_{\id_I}$ respectively.
  
  An \define{instance} of a model $X$ of $\dbl{D}$ in $\dbl{E}$ is a
  module $ H:I\proto X $ out of the terminal model $ I $ of $\dbl{D}$ in
  $\dbl{E}$ such that ``$ I $ acts trivially on the left'' in the sense that all
  laxators of the form below are identities:
  \[\begin{tikzcd}
      I & I & {X(z)} \\
      I && {X(z)}
      \arrow["{I(m)}"{inner sep=.8ex}, "\shortmid"{marking}, from=1-1, to=1-2]
      \arrow[equals, from=1-1, to=2-1]
      \arrow["{H(n)}"{inner sep=.8ex}, "\shortmid"{marking}, from=1-2, to=1-3]
      \arrow[equals, from=1-3, to=2-3]
      \arrow[""{name=0, anchor=center, inner sep=0}, "{H(mn)}"'{inner sep=.8ex}, "\shortmid"{marking}, from=2-1, to=2-3]
      \arrow["{H_{m,n}}"{description, pos=0.4}, draw=none, from=1-2, to=0]
    \end{tikzcd}\]
  A \define{co-instance}, or a presheaf on $X,$ is a module $ X\proto I$
  satisfying the dual conditions.
\end{definition}

\begin{explication}[Instance of a model, unpacked]
        Let us unfold this definition. Then from this lax functor 
        $H:\dbl{D}\times \WalkingLoose\to\dbl{E}$ we have:
        \begin{itemize}
          \item For each loose morphism $ m:x\proto y $ in $ \dbl{D} $, a loose morphism
            $ H(m):=H(m,\ell):I\proto X(y) $ in $ \dbl{E} $.
          \item For every cell as on the left below in $ \dbl{D} $, a corresponding
            cell as on the right in $\dbl{E},$ where we write $H(\alpha)$ for $H(\alpha,1_\ell)$:
            \begin{equation*}
              \begin{tikzcd}
                x & y \\
                w & z
                \arrow[""{name=0, anchor=center, inner sep=0}, "m"{inner sep=.8ex}, "\shortmid"{marking}, from=1-1, to=1-2]
                \arrow["f"', from=1-1, to=2-1]
                \arrow["g", from=1-2, to=2-2]
                \arrow[""{name=1, anchor=center, inner sep=0}, "n"'{inner sep=.8ex}, "\shortmid"{marking}, from=2-1, to=2-2]
                \arrow["\alpha"{description}, draw=none, from=0, to=1]
              \end{tikzcd}
              \qquad\leadsto\qquad
              \begin{tikzcd}
                I & Xy \\
                I & Xz
                \arrow[""{name=0, anchor=center, inner sep=0}, "{H(m)}"{inner sep=.8ex}, "\shortmid"{marking}, from=1-1, to=1-2]
                \arrow[equals, from=1-1, to=2-1]
                \arrow["Xg", from=1-2, to=2-2]
                \arrow[""{name=1, anchor=center, inner sep=0}, "{H(n)}"'{inner sep=.8ex}, "\shortmid"{marking}, from=2-1, to=2-2]
                \arrow["{H(\alpha)}"{description}, draw=none, from=0, to=1]
              \end{tikzcd}
            \end{equation*}
            \item For every composable pair $ x \xproto{m} y \xproto{n} z $ in
            $ \dbl{D} $, action cells in $ \dbl{E} $ as below, arising respectively from the laxators 
            $H_{(m,\id_\vdash),(n,\ell)}$ and $H_{(m,\ell),(n,\id_\dashv)}$ of the module $H:\dbl{D}\times \WalkingLoose\to \dbl{E}$:
      \[\begin{tikzcd}
        I & I & Xz & I & Xy & Xz \\
        I && Xz & I && Xz
        \arrow["I"{inner sep=.8ex}, "\shortmid"{marking}, from=1-1, to=1-2]
        \arrow[equals, from=1-1, to=2-1]
        \arrow["Hn"{inner sep=.8ex}, "\shortmid"{marking}, from=1-2, to=1-3]
        \arrow[equals, from=1-3, to=2-3]
        \arrow["Hm"{inner sep=.8ex}, "\shortmid"{marking}, from=1-4, to=1-5]
        \arrow[equals, from=1-4, to=2-4]
        \arrow["Xn"{inner sep=.8ex}, "\shortmid"{marking}, from=1-5, to=1-6]
        \arrow[equals, from=1-6, to=2-6]
        \arrow[""{name=0, anchor=center, inner sep=0}, "{H(mn)}"'{inner sep=.8ex}, "\shortmid"{marking}, from=2-1, to=2-3]
        \arrow[""{name=1, anchor=center, inner sep=0}, "{H(mn)}"'{inner sep=.8ex}, "\shortmid"{marking}, from=2-4, to=2-6]
        \arrow["{H^\ell_{m,n}}"{description, pos=0.3}, draw=none, from=1-2, to=0]
        \arrow["{H^r_{m,n}}"{description, pos=0.3}, draw=none, from=1-5, to=1]
      \end{tikzcd}\]
\end{itemize}
To the usual properties of a module as a lax functor, we have added the
assumption that the left actions $ H^\ell_{m,n} $ be identities.

This assumption leads to several reductions of structure, as follows:
\begin{enumerate}
  \item First note that setting $ H^\ell $ to the identity makes no sense unless we also assume that $ H(mn)=H(n) $
  for all composable pairs $ x\xproto{m}y\xproto{n} z $ of loose morphisms in $ \dbl{D} $. In particular, $ H(m)=H(m\id_y)=H(\id_y) $, using strictness of $ \dbl{D} $.
  This means we can forget about 
  $ H(n) $ except in case $ n=\id_x $ for some object $ x $ of $ \dbl{D} $. We therefore shift notation and write $ H(x) $ for $ H(\id_x) $.

  \item Similarly, consider the naturality of the laxator $H^\ell$, which says
    that given loosely composable cells in $ \dbl{D} $
    $
      \begin{tikzcd}
        x & y & z \\
        {x'} & {y'} & {z'}
        \arrow[""{name=0, anchor=center, inner sep=0}, "m"{inner sep=.8ex}, "\shortmid"{marking}, from=1-1, to=1-2]
        \arrow["f"', from=1-1, to=2-1]
        \arrow[""{name=1, anchor=center, inner sep=0}, "n"{inner sep=.8ex}, "\shortmid"{marking}, from=1-2, to=1-3]
        \arrow["g"{description}, from=1-2, to=2-2]
        \arrow["h", from=1-3, to=2-3]
        \arrow[""{name=2, anchor=center, inner sep=0}, "{m'}"'{inner sep=.8ex}, "\shortmid"{marking}, from=2-1, to=2-2]
        \arrow[""{name=3, anchor=center, inner sep=0}, "{n'}"'{inner sep=.8ex}, "\shortmid"{marking}, from=2-2, to=2-3]
        \arrow["\alpha"{description}, draw=none, from=0, to=2]
        \arrow["\beta"{description}, draw=none, from=1, to=3]
      \end{tikzcd},
    $
    we have the equality in $\dbl{E}$:

  \[\begin{tikzcd}
    I & I & Xz && I & I & Xz \\
    I & I & {Xz'} & {=} & I && Xz \\
    I && {Xz'} && I && {Xz'}
    \arrow[""{name=0, anchor=center, inner sep=0}, "{\id_I}"{inner sep=.8ex}, "\shortmid"{marking}, from=1-1, to=1-2]
    \arrow[equals, from=1-1, to=2-1]
    \arrow[""{name=1, anchor=center, inner sep=0}, "Hn"{inner sep=.8ex}, "\shortmid"{marking}, from=1-2, to=1-3]
    \arrow[equals, from=1-2, to=2-2]
    \arrow["Xh", from=1-3, to=2-3]
    \arrow["{\id_I}"{inner sep=.8ex}, "\shortmid"{marking}, from=1-5, to=1-6]
    \arrow[equals, from=1-5, to=2-5]
    \arrow["Hn"{inner sep=.8ex}, "\shortmid"{marking}, from=1-6, to=1-7]
    \arrow[equals, from=1-7, to=2-7]
    \arrow[""{name=2, anchor=center, inner sep=0}, "{\id_I}"'{inner sep=.8ex}, "\shortmid"{marking}, from=2-1, to=2-2]
    \arrow[equals, from=2-1, to=3-1]
    \arrow[""{name=3, anchor=center, inner sep=0}, "{{Hn'}}"'{inner sep=.8ex}, "\shortmid"{marking}, from=2-2, to=2-3]
    \arrow[equals, from=2-3, to=3-3]
    \arrow[""{name=4, anchor=center, inner sep=0}, "{{H(mn)}}"'{inner sep=.8ex}, "\shortmid"{marking}, from=2-5, to=2-7]
    \arrow[equals, from=2-5, to=3-5]
    \arrow["Xh"{description}, from=2-7, to=3-7]
    \arrow[""{name=5, anchor=center, inner sep=0}, "{{H(m'n')}}"'{inner sep=.8ex}, "\shortmid"{marking}, from=3-1, to=3-3]
    \arrow[""{name=6, anchor=center, inner sep=0}, "{{H(m'n')}}"'{inner sep=.8ex}, "\shortmid"{marking}, from=3-5, to=3-7]
    \arrow["{I\alpha}"{description}, draw=none, from=0, to=2]
    \arrow["{{H\beta}}"{description}, draw=none, from=1, to=3]
    \arrow["{{H^\ell_{m,n}}}"{description}, draw=none, from=1-6, to=4]
    \arrow["{{H^\ell_{m',n'}}}"{description}, draw=none, from=2-2, to=5]
    \arrow["{{H(\alpha\beta)}}"{description}, draw=none, from=4, to=6]
  \end{tikzcd}\]

  Since $ I\alpha $ and the $ H^\ell $s are all identities, this amounts to the condition that
  $ H\beta=H(\alpha\beta) $ and, in particular, that $ H\alpha=H(\id_g) $.
  We can again thus forget about $ H\alpha $ except in the case $ \alpha=\id_f $ for a tight morphism $ f $. 
  We therefore again shift notation and write $ H(f) $ for $ H(\id_f) $.

  \item Finally, given the pair of loose arrows
    $ x\xproto{m} y\xproto{n} z $, consider the associativity of the laxators of $H$ at the triple composite
    $(m,\id_0)\odot (\id_y,\ell)\odot (n,\id_1)$, which says
    that the following equation holds:
    \[\begin{tikzcd}
      I & I & Xy & Xz && I & I & Xy & Xz \\
      I && Xy & Xz & {=} & I & I && Xz \\
      I &&& Xz && I &&& Xz
      \arrow["\shortmid"{marking}, equals, from=1-1, to=1-2]
      \arrow[equals, from=1-1, to=2-1]
      \arrow["Hy"{inner sep=.8ex}, "\shortmid"{marking}, from=1-2, to=1-3]
      \arrow[""{name=0, anchor=center, inner sep=0}, "Xn"{inner sep=.8ex}, "\shortmid"{marking}, from=1-3, to=1-4]
      \arrow[equals, from=1-3, to=2-3]
      \arrow[equals, from=1-4, to=2-4]
      \arrow["\shortmid"{marking}, equals, from=1-6, to=1-7]
      \arrow[equals, from=1-6, to=2-6]
      \arrow["Hy"{inner sep=.8ex}, "\shortmid"{marking}, from=1-7, to=1-8]
      \arrow[equals, from=1-7, to=2-7]
      \arrow["Xn"{inner sep=.8ex}, "\shortmid"{marking}, from=1-8, to=1-9]
      \arrow[equals, from=1-9, to=2-9]
      \arrow[""{name=1, anchor=center, inner sep=0}, "Hy"'{inner sep=.8ex}, "\shortmid"{marking}, from=2-1, to=2-3]
      \arrow[equals, from=2-1, to=3-1]
      \arrow[""{name=2, anchor=center, inner sep=0}, "Xn"'{inner sep=.8ex}, "\shortmid"{marking}, from=2-3, to=2-4]
      \arrow[equals, from=2-4, to=3-4]
      \arrow["I"{inner sep=.8ex}, "\shortmid"{marking}, equals, from=2-6, to=2-7]
      \arrow[equals, from=2-6, to=3-6]
      \arrow[""{name=3, anchor=center, inner sep=0}, "{{Hz}}"'{inner sep=.8ex}, "\shortmid"{marking}, from=2-7, to=2-9]
      \arrow[equals, from=2-9, to=3-9]
      \arrow[""{name=4, anchor=center, inner sep=0}, "Hz"'{inner sep=.8ex}, "\shortmid"{marking}, from=3-1, to=3-4]
      \arrow[""{name=5, anchor=center, inner sep=0}, "{{H(\id_yn)=Hz}}"'{inner sep=.8ex}, "\shortmid"{marking}, from=3-6, to=3-9]
      \arrow["{{H^\ell_{m,y}}}"{description, pos=0.4}, draw=none, from=1-2, to=1]
      \arrow["{1_{Xn}}"{description}, draw=none, from=0, to=2]
      \arrow["{{H^r_{y,n}}}"{description, pos=0.4}, draw=none, from=1-8, to=3]
      \arrow["{{H^r_{m\id_y,n}=H^r_{m,n}}}"{description}, draw=none, from=2-3, to=4]
      \arrow["{{H^\ell_{m,\id_yn}}}"{description}, draw=none, from=2-7, to=5]
    \end{tikzcd}\]
    Since the left actions $ H^\ell $ are identities by assumption, this amounts to
    the condition on the right actions that $ H^r_{m,n}=H^r_{\id_y,n}.$ Thus we
    may simply write $ Hn=H^r_{\id_y,n} $ and forget about the other action
    cells. \qedhere
  \end{enumerate}
\end{explication}

This reduces the material of an instance as follows.

\begin{definition}[Instance of a model, explicit definition] \label{def:instance}
    Let $\dbl{D}$ be a double theory and let $\dbl{E}$ be a double category with a terminal object $I.$
    An \define{instance} $H$ of a model $X:\dbl{D}\to\dbl{E} $ of the theory $\dbl{D}$ is given by
    the following data.
    \begin{itemize}
        \item A functor $H_0:\dbl{D}_0\to \dbl{E}_1$ such that $H_0\cdot s:\dbl{D}_0\to\dbl{E}_0$ is 
        constant at $I$ and $H_0 \cdot t = X_0$.
        We denote $H_0(d)$ by $H(d):I\proto X(d)$ for each object $d$ of $\dbl{D}$,
        and similarly for tight morphisms.
      \item A natural transformation
        $H_1: s\cdot H_0\odot X_1\To t\cdot H_0: \dbl{D}_1\to \dbl{E}_1$ such that 
        $H_1* s$ is constant at $I$ and $H_1*t$ is constant at $t\cdot X_0.$
        Thus $H_1$ consists of components
        $H(m) \coloneqq H_1(m)$ for each $m: d \to d'$ in $\dbl{D}$, with boundaries 
        as displayed below:
        \[\begin{tikzcd}
          I & Xd & {Xd'} \\
          I && {Xd'}
          \arrow["Hd"{inner sep=.8ex}, "\shortmid"{marking}, from=1-1, to=1-2]
          \arrow[equals, from=1-1, to=2-1]
          \arrow["Xm"{inner sep=.8ex}, "\shortmid"{marking}, from=1-2, to=1-3]
          \arrow[equals, from=1-3, to=2-3]
          \arrow[""{name=0, anchor=center, inner sep=0}, "{Hd'}"'{inner sep=.8ex}, "\shortmid"{marking}, from=2-1, to=2-3]
          \arrow["Hm"{description}, draw=none, from=1-2, to=0]
        \end{tikzcd}\]
    \end{itemize}

    These data are associative and unital in the sense that the following two diagrams commute, in the functor categories
    indicated below each diagram, and where $s,c,t:\dbl{D}_1\times_{\dbl{D}_0}\dbl{D}_1\to\dbl{D}_0$ are the 
    functors picking out the source, center, and target object of a path of loose arrows. 
\[\begin{tikzcd}
	{(s\cdot H_0\odot X_1)\odot X_1} & {c\cdot H_0\odot X_1} & {H_0\odot \mathrm{id}_{X -}} & {H_0\odot X(\mathrm{id}_{-})} \\
	{s\cdot H_0\odot X_1} & {t\cdot H_0} && {H_0} \\
	{\dbl{D}_1\times_{\dbl{D}_0}\dbl{D}_1} & {\dbl{E}_1} & {\dbl{D}_0} & {\dbl{E}_1}
	\arrow["{H_1\odot \mathrm{id}_{X_1}}", from=1-1, to=1-2]
	\arrow["{\mathrm{id}_{H_0}\odot X_{-,-}}"', from=1-1, to=2-1]
	\arrow["{H_1}", from=1-2, to=2-2]
	\arrow["{\mathrm{id}_{H_0}\odot X_{-}}", from=1-3, to=1-4]
	\arrow[equals, from=1-3, to=2-4]
	\arrow["{H_1(\mathrm{id}_{-})}", from=1-4, to=2-4]
	\arrow["{H_1}"', from=2-1, to=2-2]
	\arrow[from=3-1, to=3-2]
	\arrow[from=3-3, to=3-4]
\end{tikzcd}\]

  When $\dbl{D}$ and $X$ are cartesian, we require that $H$ is also cartesian in that the canonical 
  comparison cells below are isomorphisms:
  \begin{equation*}
    \begin{tikzcd}
      I && {X(d_1 \times d_2)} \\
      I && {X(d_1) \times X(d_2)}
      \arrow[""{name=0, anchor=center, inner sep=0}, "{H(d_1 \times d_2)}"{inner sep=.8ex}, "\shortmid"{marking}, from=1-1, to=1-3]
      \arrow[equals, from=1-1, to=2-1]
      \arrow["{\langle X\pi_{d_1}, X\pi_{d_2} \rangle}", from=1-3, to=2-3]
      \arrow[""{name=1, anchor=center, inner sep=0}, "{H(d_1) \times H(d_2)}"'{inner sep=.8ex}, "\shortmid"{marking}, from=2-1, to=2-3]
      \arrow["{\langle H\pi_{d_1}, H\pi_{d_2} \rangle}"{description}, draw=none, from=0, to=1]
    \end{tikzcd}
    \qquad\qquad
    \begin{tikzcd}
      I & {X(I)} \\
      I & I
      \arrow[""{name=0, anchor=center, inner sep=0}, "{H(I)}"{inner sep=.8ex}, "\shortmid"{marking}, from=1-1, to=1-2]
      \arrow[equals, from=1-1, to=2-1]
      \arrow["{!_{X(I)}}", from=1-2, to=2-2]
      \arrow[""{name=1, anchor=center, inner sep=0},  "\shortmid"{marking}, equals, from=2-1, to=2-2]
      \arrow["{!_{H(I)}}"{description}, draw=none, from=0, to=1]
    \end{tikzcd}
    \qedhere
  \end{equation*}
\end{definition}

A \define{morphism} of instances of the model $X$ is a globular modulation (\cref{def:modulation}) between
the corresponding modules.
As in the  definition of instances themselves, we can substantially simplify the data involved.
In fact, the naturality of laxators of a modulation $\mu:H\proto K$ 
between instances
of the same model $X$ of theory $\dbl{D}$ (and over $\id_X$) implies that $\mu_n=\mu_
{m\odot n}$ for any composable loose arrows $m$ and $n$ in $\dbl{D},$ thus that
$\mu_m=\mu_{\id_x\odot m}=\mu_{\id_x}$ if $m:x\proto y.$
Thus it suffices to give the components of $\mu$ at objects,
and we recover the following definition.

\begin{definition}[Category of instances]
If $H$ and $K$ are instances of a fixed model $X:\dbl{D}\to \dbl{E},$ then a \define{morphism} $\mu:H\to K$ is
given by a globular natural transformation $\mu: H_0\To K_0:\dbl{D}_0\to \dbl{E}_1$ which is equivariant in the 
sense that the following square commutes in the category of functors 
$\dbl{D}_1\to \dbl{E}_1$: 
\[\begin{tikzcd}
	{s\cdot H_0\odot X_1} & {s\cdot K_0\odot X_1} \\
	{t\cdot H_0} & {t\cdot K_0}
	\arrow["{s\cdot \mu\odot X_1}", from=1-1, to=1-2]
	\arrow["{H_1}"', from=1-1, to=2-1]
	\arrow["{K_1}", from=1-2, to=2-2]
	\arrow["{t\cdot \mu}"', from=2-1, to=2-2]
\end{tikzcd}\]
We denote the \define{category of instances} of $X$ by $\Inst(X)$.
\end{definition}

The category of instances, as constructed above, is actually recoverable in terms of \emph{loose} transformations,
a notion of 2-morphism of lax functor which, like ordinary transformations, has components indexed by objects, 
but, like modules, has as these components loose rather than tight arrows. 
Since we won't use loose transformations otherwise in this paper, we defer supporting results to 
Appendix \ref{app:loose-transformations}. 
\begin{proposition}[Instances as loose transformations]\label{prop:instances-as-loose-transformations}
  Let $\dbl{D}$ be a double theory, let $\dbl{E}$ be a double category with terminal object $I,$
  and let $X:\dbl{D}\to \dbl{E}$ be a model of $\dbl{D}$ in $\dbl{E}.$
  Then the category $\Inst(X)$ of $X$ in the sense of \cref{def:compressed-instance} 
  is naturally equivalent to the category of loose transformations and modulations $\cat{Loose}(I,X).$
\end{proposition}
\begin{proof}
Since the terminal model $I$ is strict, by \cref{thm:loose-transformations-into-modules}, 
the category of loose transformations and modulations $I\to X$ is equivalent to the category 
of modules $I\proto X$ with trivial left action and modulations between them, that is, to 
$\Inst(X).$
\end{proof}

We next state the most obvious sense in which $\Inst(X)$ is functorial in $X$.

\begin{proposition}[Functoriality of instances]\label{prop:instances_functoriality}
  Consider models $X,Y:\dbl{D}\to \dbl{E}$ such that $\dbl{E}$ is an equipment 
  (\cref{def:equipment}) with a terminal object. For any morphism of models $\alpha:X\to Y$,
  there is an induced ``substitution'' functor $\alpha^*:\Inst(Y)\to \Inst(X).$
\end{proposition}
\begin{proof}
First, we define $\alpha^*$ on objects. Given a $Y$-instance $H,$ we construct. an $X$-instance $\alpha^*H$ as follows:
\begin{itemize}
\item We define $\alpha^*H(d)$ as the domain of the restriction cell 
  \[\begin{tikzcd}
    I & Xd \\
    I & Yd
    \arrow[""{name=0, anchor=center, inner sep=0}, "{\alpha^*Hd}"{inner sep=.8ex}, "\shortmid"{marking}, from=1-1, to=1-2]
    \arrow[equals, from=1-1, to=2-1]
    \arrow["{\alpha_d}", from=1-2, to=2-2]
    \arrow[""{name=1, anchor=center, inner sep=0}, "Hd"'{inner sep=.8ex}, "\shortmid"{marking}, from=2-1, to=2-2]
    \arrow["{\mathrm{res}}"{description}, draw=none, from=0, to=1]
  \end{tikzcd}\] 
\item We define $\alpha^*H(f)$ as the unique solution of the equation below, using the naturality
of $\alpha$ and the universal property of restriction: 
  \[\begin{tikzcd}
    I & Xd && I & Xd \\
    I & Yd & {=} & I & {Xd'} \\
    I & {Yd'} && I & {Yd'}
    \arrow[""{name=0, anchor=center, inner sep=0}, "{\alpha^*Hd}"{inner sep=.8ex}, "\shortmid"{marking}, from=1-1, to=1-2]
    \arrow[equals, from=1-1, to=2-1]
    \arrow["{\alpha_d}", from=1-2, to=2-2]
    \arrow[""{name=1, anchor=center, inner sep=0}, "{\alpha^*Hd}"{inner sep=.8ex}, "\shortmid"{marking}, from=1-4, to=1-5]
    \arrow[equals, from=1-4, to=2-4]
    \arrow["Xf", from=1-5, to=2-5]
    \arrow[""{name=2, anchor=center, inner sep=0}, "Hd"'{inner sep=.8ex}, "\shortmid"{marking}, from=2-1, to=2-2]
    \arrow[equals, from=2-1, to=3-1]
    \arrow["Yf"{description}, from=2-2, to=3-2]
    \arrow[""{name=3, anchor=center, inner sep=0}, "{\alpha^*Hd'}"{inner sep=.8ex}, "\shortmid"{marking}, from=2-4, to=2-5]
    \arrow[equals, from=2-4, to=3-4]
    \arrow["{\alpha_{d'}}", from=2-5, to=3-5]
    \arrow[""{name=4, anchor=center, inner sep=0}, "{Hd'}"'{inner sep=.8ex}, "\shortmid"{marking}, from=3-1, to=3-2]
    \arrow[""{name=5, anchor=center, inner sep=0}, "{Hd'}"'{inner sep=.8ex}, "\shortmid"{marking}, from=3-4, to=3-5]
    \arrow["{\mathrm{res}}"{description}, draw=none, from=0, to=2]
    \arrow["{\alpha^*Hf}"{description}, draw=none, from=1, to=3]
    \arrow["Hf"{description}, draw=none, from=2, to=4]
    \arrow["{\mathrm{res}}"{description}, draw=none, from=3, to=5]
  \end{tikzcd}\]
\item We define $\alpha^*H(m)$ for a loose morphism $m:d\proto d'$ in $\dbl{D}$ 
as the unique solution of the equation below:
  \[\begin{tikzcd}
    I & Xd & {Xd'} && I & Xd & {Xd'} \\
    I & Yd & {Yd'} & {=} & I && {Xd'} \\
    I && {Yd'} && I && {Yd'}
    \arrow[""{name=0, anchor=center, inner sep=0}, "{\alpha^*Hd}"{inner sep=.8ex}, "\shortmid"{marking}, from=1-1, to=1-2]
    \arrow[equals, from=1-1, to=2-1]
    \arrow[""{name=1, anchor=center, inner sep=0}, "Xm"{inner sep=.8ex}, "\shortmid"{marking}, from=1-2, to=1-3]
    \arrow["{\alpha_d}"{description}, from=1-2, to=2-2]
    \arrow["{\alpha_{d'}}", from=1-3, to=2-3]
    \arrow["{\alpha^*Hd}"{inner sep=.8ex}, "\shortmid"{marking}, from=1-5, to=1-6]
    \arrow[equals, from=1-5, to=2-5]
    \arrow["Xm"{inner sep=.8ex}, "\shortmid"{marking}, from=1-6, to=1-7]
    \arrow[equals, from=1-7, to=2-7]
    \arrow[""{name=2, anchor=center, inner sep=0}, "Hd"'{inner sep=.8ex}, "\shortmid"{marking}, from=2-1, to=2-2]
    \arrow[equals, from=2-1, to=3-1]
    \arrow[""{name=3, anchor=center, inner sep=0}, "Ym"', from=2-2, to=2-3]
    \arrow[equals, from=2-3, to=3-3]
    \arrow[""{name=4, anchor=center, inner sep=0}, "{\alpha^*Hd'}"'{inner sep=.8ex}, "\shortmid"{marking}, from=2-5, to=2-7]
    \arrow[equals, from=2-5, to=3-5]
    \arrow["{\alpha_{d'}}"', from=2-7, to=3-7]
    \arrow[""{name=5, anchor=center, inner sep=0}, "{Hd'}"'{inner sep=.8ex}, "\shortmid"{marking}, from=3-1, to=3-3]
    \arrow[""{name=6, anchor=center, inner sep=0}, "{Hd'}"'{inner sep=.8ex}, "\shortmid"{marking}, from=3-5, to=3-7]
    \arrow["{\mathrm{res}}"{description}, draw=none, from=0, to=2]
    \arrow["{\alpha_m}"{description}, draw=none, from=1, to=3]
    \arrow["{\alpha^*Hm}"{description}, draw=none, from=1-6, to=4]
    \arrow["Hm"{description}, draw=none, from=2-2, to=5]
    \arrow["{\mathrm{res}}"{description}, draw=none, from=4, to=6]
  \end{tikzcd}\]
\end{itemize} 
Furthermore, if $\mu:H\to K$ is a morphism of $Y$-instances, we define $\alpha^*\mu: \alpha^* H \to \alpha^* K$ to have as
components the unique solutions of the equations 
\[\begin{tikzcd}
	I & Xd && I & Xd \\
	I & Yd & {=} & I & Xd \\
	I & Yd && I & Yd
	\arrow[""{name=0, anchor=center, inner sep=0}, "{\alpha^*Hd}"{inner sep=.8ex}, "\shortmid"{marking}, from=1-1, to=1-2]
	\arrow[equals, from=1-1, to=2-1]
	\arrow["{\alpha_d}", from=1-2, to=2-2]
	\arrow[""{name=1, anchor=center, inner sep=0}, "{\alpha^*Hd}"{inner sep=.8ex}, "\shortmid"{marking}, from=1-4, to=1-5]
	\arrow[equals, from=1-4, to=2-4]
	\arrow[equals, from=1-5, to=2-5]
	\arrow[""{name=2, anchor=center, inner sep=0}, "Hd"'{inner sep=.8ex}, "\shortmid"{marking}, from=2-1, to=2-2]
	\arrow[equals, from=2-1, to=3-1]
	\arrow[equals, from=2-2, to=3-2]
	\arrow[""{name=3, anchor=center, inner sep=0}, "{\alpha^*Kd}"'{inner sep=.8ex}, "\shortmid"{marking}, from=2-4, to=2-5]
	\arrow[equals, from=2-4, to=3-4]
	\arrow["{\alpha_{d}}", from=2-5, to=3-5]
	\arrow[""{name=4, anchor=center, inner sep=0}, "Kd"'{inner sep=.8ex}, "\shortmid"{marking}, from=3-1, to=3-2]
	\arrow[""{name=5, anchor=center, inner sep=0}, "Kd"'{inner sep=.8ex}, "\shortmid"{marking}, from=3-4, to=3-5]
	\arrow["{\mathrm{res}}"{description}, draw=none, from=0, to=2]
	\arrow["{\alpha^*\mu d}"{description}, draw=none, from=1, to=3]
	\arrow["{\mu_d}"{description}, draw=none, from=2, to=4]
	\arrow["{\mathrm{res}}"{description}, draw=none, from=3, to=5]
\end{tikzcd}\]
It is straightforward to check that $\alpha^*$ sends instances to instances (including in the case 
that $\dbl{D}$ is cartesian) and morphisms to 
morphisms; we leave the details to the reader, noting that in the main case of instances of $\Span$-valued 
models, we will reconstruct this result more efficiently below.
\end{proof}

\subsection{Examples of instances}

All of our examples of instances will be of models valued in the double category
of spans. With this focus it will be beneficial to describe the categories of
models of discrete double theories, using the notation for the discrete
embeddings from \cref{def:discreteDoubleCategory}:

\begin{proposition}[Span-valued models of discrete double theories]\label{prop:span_models_of_discrete_theories}
For any small category $C:$ 
\begin{enumerate}
\item The category $\Lax(\dbl{T}C,\Span)$ is canonically equivalent to the category 
$\CAT(C,\Cat)$ of functors from $C$ to $\Cat$ and (strict) natural transformations.
\item The category $\Lax(\dbl{L}C,\Span)$ is canonically equivalent to the slice category $\cat{Cat}/C.$

\end{enumerate}
\end{proposition}
\begin{proof}

\begin{enumerate}
\item This follows from \cref{prop:functors-out-of-discrete}, together with the fact (\cite[Proposition 5.14]{generalized-multicategories-2010})
that lax functors into $\Span$ are equivalent to normal lax functors into $\dbl{P}\mathsf{rof}.$
\item Using again the fact that $\Lax(\dbl{L}C,\Span)\simeq \Lax_{n}(\dbl{L}C,\dbl{P}\mathsf{rof}),$ 
we note that in this loosely discrete case, we may forget the tight direction entirely and reduce to 
normal lax functors of bicategories $C\to \twocat{Prof}.$ It is well-known that such normal lax functors 
$C\to \twocat{Prof}$ are equivalent to categories over $C$ (see \cite[Section 7]{distAtWork}).

However, to get the right correspondence on morphisms, we must
work directly with $\Span$-valued lax functors. Given such a functor
$F:\dbl{L}C\to \Span,$ we define a category $\int F\to C$ as follows:
\begin{itemize}
\item The objects of $\int F$ over $c$ are the elements of the set $F(c)$, i.e., the objects of the category $F(c).$
\item A morphism from $x\in F(c)$ to $y\in F(c')$ over $f:c\to c'$ is an element of the set $F(f)(x,y).$
\item Composition arises from the laxators of $F,$ and identity 
morphisms from the unitors, while associativity and unitality 
correspond to those properties of $F$'s laxators. 
\end{itemize}

Given a morphism $\alpha:F\To G:\dbl{L}C\to \Span,$ we get a functor 
$\int \alpha:\int F\to \int G$ over $C$ as follows: 
\begin{itemize}
\item On objects, $\int \alpha$ acts according to the object components of $\alpha.$ 
\item On morphisms, $\int \alpha$ acts according to the loose morphism components of $\alpha.$
\item The intertwining of $\alpha$ with laxators implies that $\int \alpha$ respects composition,
and similarly for unitors and identities.
\end{itemize}

Conversely, given a functor $A: \int F\to \int G$ over $C,$ we get 
a morphism $\alpha^A:F\to G$ with object components given by the action of $A$ on objects,
loose morphism components given by the action of $A$ on morphisms, and 
the intertwining axioms of $\alpha^A$ following from the functoriality of $A.$ 
By construction $\alpha^{\int\alpha}=\alpha$ and $\int \alpha^A=A,$ so we have the equivalence stated.
\end{enumerate}
\end{proof}

In particular, we see that $\Lax(\dbl{T} C,\Span)$ is a non-full subcategory of $\Lax(\dbl{L} C,\Span)$, via the Grothendieck construction interchanging a functor $C\to \Cat$ with a \emph{split} fibration 
over $C.$ 

\begin{example}[Instances of categories and functors]
As has been suggested, an instance of a model $X:\WalkingOb\to\Span$ of the terminal double theory $\WalkingOb$ is 
a co-presheaf on $X.$ Indeed, if $\WalkingOb=\{\ob\xproto{\mor}\ob\}$ with 
$\mor=\id_\ob$, then an instance $H$ of $X$ consists of a fibered set 
$H\ob\to X\ob$ and an action cell $H\ob\times_{X\ob} X\mor \to H\ob$, which is precisely 
the data of a functor $X\to \Set$; associativity and unitality constraints for either concept interconvert. 
Similarly, a co-instance is a presheaf. Note that for this theory, instances and coinstances are the same as
arbitrary modules from and to $ I $, as will turn out to be true whenever the 
double theory being modeled has only trivial loose morphisms (see \cref{thm:loose-transformations-into-modules}).

Now consider the theory $\WalkingTight=\{\TopOb\xto{t}\BotOb\}$ of a functor.
A module $ M:F\proto G $ is a lax functor $\WalkingTight\times\WalkingLoose\to\Span.$
Since $\WalkingTight\times\WalkingLoose\cong\WalkingSquare,$ this is simply
a modulation with $F,G$ as its tight boundary components. As in the case of the terminal theory, 
there is nothing for the more restricted
definition of instance to do here: an instance $M:I\proto (G:C\to D)$ is 
given by a pair of functors $X:C\to\Set$ and $Y:D\to\Set$ and a natural transformation $X\to G^*Y.$
\end{example}

Next we give the minimal example in which instances differ from modules from $I.$ 
\begin{example}[Instances of profunctors]
Consider the theory $\WalkingLoose$ of a loose morphism and a model $X:\WalkingLoose\to \Span,$ 
which may be identified with a profunctor $X_\ell:X_0\proto X_1$ of ordinary categories. 
An instance $H$ of $X$ consists of: 
\begin{itemize}
\item Functors $H_0:X_0\to \Set$ and $H_1:X_1\to \Set.$
\item A natural transformation $H_0\odot X_\ell\to H_1,$ where we compose 
$H_0:I\proto X_0$ with $X_\ell$ as profunctors; thus 
$H_0\odot X_\ell (x_1)$ is the coend $\int^{x_0\in X_0} H_0(x_0)\times X_\ell(x_0,x_1).$
\end{itemize}
This recovers the notion of ``attributed $C$-set'' (acset) from 
categorical database theory \cite{acsets-2022}, with the caveat that we
must slice the category of instances over an instance $X_1\to\Set$ of datatypes to recover the 
intended morphisms of acsets.
\end{example}

We now highlight how instances differ from arbitrary modules from $1$ in this case.
\begin{remark}[Instances versus modules from $1$, for profunctors]
Consider an arbitrary module
$ M:W\proto X:\WalkingLoose\to\Span. $ Thus $M$ is itself a lax functor $\WalkingLoose^2\to \Span.$ 
Now $\WalkingLoose^2$ looks like this, containing in particular five non-identity loose arrows: 
\[\begin{tikzcd}
	{(\LeftOb,\LeftOb)} & {(\LeftOb,\RightOb)} \\
	{(\RightOb,\LeftOb)} & {(\RightOb,\RightOb)}
	\arrow["\shortmid"{marking}, from=1-1, to=1-2]
	\arrow["\shortmid"{marking}, from=1-1, to=2-1]
	\arrow["\shortmid"{marking}, from=1-1, to=2-2]
	\arrow["\shortmid"{marking}, from=1-2, to=2-2]
	\arrow["\shortmid"{marking}, from=2-1, to=2-2]
\end{tikzcd}.\]
The two triangles shown commute on the nose in $\WalkingLoose^2,$ but they receive laxators under a lax
functor. Therefore $M$ comprises
a diagram of the following shape in the bicategory of profunctors:
\[\begin{tikzcd}
	{W_0} & {W_1} \\
	{X_0} & {X_1}
	\arrow["{W_\ell}"{inner sep=.8ex}, "\shortmid"{marking}, from=1-1, to=1-2]
	\arrow["\shortmid"{marking}, from=1-1, to=2-1]
	\arrow[""{name=0, anchor=center, inner sep=0}, "\shortmid"{marking}, from=1-1, to=2-2]
	\arrow["\shortmid"{marking}, from=1-2, to=2-2]
	\arrow["{X_\ell}"{inner sep=.8ex}, "\shortmid"{marking}, from=2-1, to=2-2]
	\arrow[between={0}{0.8}, Rightarrow, from=1-2, to=0]
	\arrow[between={0}{0.8}, Rightarrow, from=2-1, to=0]
\end{tikzcd}\]

Now suppose that $ W=I $ is terminal. Then such a module consists of 
functors $ M_0: X_0\to \Set,M_1,M_1': X_1\to \Set $ and morphisms 
$ M_0\odot X_\ell\To M_1$ and $M_1'\To M_1$ in $\Cat(X_1,\Set).$
An instance of $ X $ throws away precisely $M_1'$ and its morphism into 
$M_1,$ rather mysterious extra data which has no obviously useful interpretation.
\end{remark}

We now give a few examples in the cartesian setting.
Recall from \cref{def:instance} that an instance of a model of a cartesian double theory 
is one whose object-component functor preserves finite products. 
We don't strictly need the following result for the general theory, but it's important
in practice for checking when one has finished defining a would-be cartesian instance.

\begin{lemma}[Action cells at products are determined]\label{lem:action-cells-products}
Given a finite family $m_i:d_i\proto d'_i$ of loose arrows in a cartesian double theory $\dbl{D}$
and an instance $P:I\proto X$ of a model $X$ of $\dbl{D},$
the action cell $P_{\Pi m_i}$ is canonically determined by the
action cells $P_{m_i}$ for each $i.$ 
\end{lemma}
\begin{proof}
By the naturality of action cells for $P$ taken at the $i$th projection
cell $\pi_i: \prod m_i \to m_i,$ we have the equation
\[\begin{tikzcd}
	I & {X(\prod d_i)} & {X(\prod d'_i)} && I & {X(\prod d_i)} & {X(\prod d'_i)} \\
	I && {X(\prod d'_i)} & {=} & I & {X(d_i)} & {X(d'_i)} \\
	I && {X(d'_i)} && I && {X(d'_i)}
	\arrow["{P(\prod d_i)}", "\shortmid"{marking}, from=1-1, to=1-2]
	\arrow[equals, from=1-1, to=2-1]
	\arrow["{X(\prod m_i)}", "\shortmid"{marking}, from=1-2, to=1-3]
	\arrow[equals, from=1-3, to=2-3]
	\arrow[""{name=0, anchor=center, inner sep=0}, "{P(\prod d_i)}", "\shortmid"{marking}, from=1-5, to=1-6]
	\arrow[from=1-5, to=2-5]
	\arrow[""{name=1, anchor=center, inner sep=0}, "{X(\prod m_i)}", "\shortmid"{marking}, from=1-6, to=1-7]
	\arrow[from=1-6, to=2-6]
	\arrow[from=1-7, to=2-7]
	\arrow[""{name=2, anchor=center, inner sep=0}, "{P(\prod d'_i)}"', "\shortmid"{marking}, from=2-1, to=2-3]
	\arrow[equals, from=2-1, to=3-1]
	\arrow["{X(\pi_i)}", from=2-3, to=3-3]
	\arrow[""{name=3, anchor=center, inner sep=0}, "{P(d_i)}"', "\shortmid"{marking}, from=2-5, to=2-6]
	\arrow[from=2-5, to=3-5]
	\arrow[""{name=4, anchor=center, inner sep=0}, "{X(m_i)}"', "\shortmid"{marking}, from=2-6, to=2-7]
	\arrow[from=2-7, to=3-7]
	\arrow[""{name=5, anchor=center, inner sep=0}, "{P(d'_i)}"', "\shortmid"{marking}, from=3-1, to=3-3]
	\arrow[""{name=6, anchor=center, inner sep=0}, "{P(d'_i)}"', "\shortmid"{marking}, from=3-5, to=3-7]
	\arrow["{P(\prod m_i)}"{description, pos=0.4}, draw=none, from=1-2, to=2]
	\arrow["{P(\pi_i)}"{description}, draw=none, from=0, to=3]
	\arrow["{X(\pi_i)}"{description}, draw=none, from=1, to=4]
	\arrow["{P(\pi_i)}"{description}, draw=none, from=2, to=5]
	\arrow["{P(m_i)}"{description}, draw=none, from=2-6, to=6]
\end{tikzcd}\]
which gives the result, by varying $i$ and using the universal property of 
$\prod m_i.$
\end{proof}

Similarly, the action of $P$ on a tight arrow with product codomain
is uniquely determined by its projections and internal functoriality.
Thus to specify an instance $P:I\proto X$ it is enough to 
give $P$ on objects, tight arrows, and loose arrows of 
$\dbl{D}$ generating $\dbl{D}$ under products and pairings. Of course
in general, not every assignment on such a datum will extend to a 
correct instance.

A particularly important example in the cartesian setting is multicategories.

\begin{example}[Instances of multicategories]
  Consider the theory $\Prom$ of promonoids \cite[Theory 6.9]{cartesian-double-theories-2024},
  the cartesian double theory containing a single object $x$ and a single loose morphism
  $p_n:x^n\proto x$ for each $n \geq 0$. As discussed at \cite[Theory 6.9]{cartesian-double-theories-2024},
  a model $M$ of $\Prom$ is precisely a multicategory.

  Now consider an instance $X:1\proto M$ of such a multicategory $M.$ 
  By \cref{lem:action-cells-products}, we can satisfy ourselves with giving the single 
  loose morphism $Xx:1\proto Mx,$ which provides a family of sets over the 
  object set of $M$; for each number $n\geq 0$, we must also provide action cells
  of the form
\[\begin{tikzcd}
	1 & {Mx^n} & Mx \\
	1 && Mx
	\arrow["{Xx^n}"{inner sep=.8ex}, "\shortmid"{marking}, from=1-1, to=1-2]
	\arrow[equals, from=1-1, to=2-1]
	\arrow["{Mp_n}"{inner sep=.8ex}, "\shortmid"{marking}, from=1-2, to=1-3]
	\arrow[equals, from=1-3, to=2-3]
	\arrow[""{name=0, anchor=center, inner sep=0}, "Xx"'{inner sep=.8ex}, "\shortmid"{marking}, from=2-1, to=2-3]
	\arrow[between={0.3}{0.7}, Rightarrow, from=1-2, to=0]
\end{tikzcd}.\]
Thus, for each natural number $n$, objects $m_1,\ldots,m_n\in M$, multimorphism
$f:m_1,\ldots,m_n\to m,$ and elements $x_1,\ldots,x_n\in Xx$ over the objects
$m_i,$ we can produce an object $\bar m$ over $m.$ It is easy to see
that associativity and unitality of the action cells amounts to assembling
them into a multifunctor $M\to \Set$, recovering the most familiar notion
of copresheaf on a multicategory.

\end{example}

\begin{remark}[Instances versus modules from $1$ for multicategories]

It may be enlightening to observe just how far an instance is from 
a module out of $1$ in this case, with a larger theory.

Consider a multicategory $\cat{M}$ as a model of the theory $\Prom$ above.
The terminal multicategory $I$ has a single object $\bullet$ and a single multimorphism
$(\bullet)^k\to \bullet$ for each arity $k$. A (cartesian) module $H: I \proto \cat{M}$ contains,
by definition, a span $H_p: I^m\proto M_\ob^n$ for each loose arrow $p:x^m\proto x^n$ in $\Prom$.
Since $H$ preserves products, these are all uniquely determined by the choices of
$H_{p_m}:I\proto M_\ob$, and most of the action cells will be similarly
redundant.

Taking this reduction into account, the only other data is the action
cells for each $p:x^m\proto x^n.$ We draw in the loose bicategory of $\Span$ 
as there are no nonstructural tight arrows in $\Prom$:
\[\begin{tikzcd}[sep=large]
	{I^m} & {I^n} \\
	{\cat{M}_\ob^n} & {\cat{M}_\ob}
	\arrow["Ip"{inner sep=.8ex}, "\shortmid"{marking}, from=1-1, to=1-2]
	\arrow["Hp"'{inner sep=.8ex}, "\shortmid"{marking}, from=1-1, to=2-1]
	\arrow[""{name=0, anchor=center, inner sep=0}, "{Hp_m}"{description, pos=0.3}, "\shortmid"{marking}, from=1-1, to=2-2]
	\arrow["{Hp_n}"{inner sep=.8ex}, "\shortmid"{marking}, from=1-2, to=2-2]
	\arrow["{Mp_n}"'{inner sep=.8ex}, "\shortmid"{marking}, from=2-1, to=2-2]
	\arrow[between={0}{0.7}, Rightarrow, from=1-2, to=0]
	\arrow[between={0}{0.7}, Rightarrow, from=2-1, to=0]
\end{tikzcd}\]

We can thus think of $H$ as follows. First, equip each object $a\in\cat{M}$ with
a set of hetero-multimorphisms $(\bullet)^k\to a$ for each $k\ge 0$, given by the
value of $Hp_k$ at $a.$ The lower-left action cell shown above then says how to compose 
$n$ such heteromorphisms whose
codomains form a list $(a_i)_{i=1}^n$ with a morphism $(a_i)\to b$ in $\cat M$.
Thus far, together with the compatibility of the operation above with composition in 
$\cat M,$ we have a multifunctor from $\cat{M}$ into a multicategory
of $\mathbb{N}$-graded sets corresponding to the monoidal product on
such sets in which $(A\otimes B)_c=\sum_{a+b=c} A_a\times B_b$.\footnote{This
is the Day convolution product for presheaves on the discrete monoidal
category $(\mathbb{N}, +, 0)$.}

Now, the upper-right action cell in the diagram above 
(which is just what we trivialize in an $\cat{M}$-instance, forcing $m$-ary heteromorphisms
to coincide with unary ones)
shows how every map $p:x^\ell\proto x^m$ in $\Prom$ produces a corresponding map from $m$-ary
to $\ell$-ary heteromorphisms into each object, functorial in $p.$ Observe that $\Prom$ has as 
loose category (being a strict double category, there is a loose 1-category here) a skeletal category
of finite totally ordered sets and order-preserving maps. Since the empty order is allowed in this 
case, this horizontal category is that usually called the augmented simplex category $\Delta_+.$
This upgrades the
$\mathbb{N}$-graded sets $(X_m)$ discussed above into presheaves on $\Delta_+,$ that is,
augmented simplicial sets.

In summary, a module out of $1$ over the theory of multicategories 
is a multifunctor from its codomain into augmented simplicial sets.
This is interesting enough, and illustrates something interesting 
about how richly structured are (two-sided) modules between multicategories, 
but this is not the usually desirable notion of instance for a multicategory.\footnote{In this case, the divergence between instances and modules from $I$ seems parallel
to the divergence between the terminal multicategory and the walking object, 
but we do not attempt to firm up this analogy here.} 
While the fervent partisan of double category theory might argue that this in fact shows that
one \emph{should} privilege multifunctors into augmented simplicial sets as well as, or in 
preference to,  those into sets, we can counter at least partly that such multifunctors will 
not admit the Grothendieck construction we describe below. 

\end{remark}

We finally give two less crucial examples of instances of models for theories with 
no nontrivial loose morphisms, to flesh out the range of applicability of the theory.

\begin{example}[Instances of cartesian categories]
Let $\dbl{C}$ be the cartesian double theory of cartesian categories as in 
\cite[Theory 6.14]{cartesian-double-theories-2024},
generated under products by an object $x$ equipped with right adjoints
$\times:x^2\to x$ and $I:1\to x$ to its diagonal and 
terminal arrows. 

Then, as expected, a model of $\dbl{C}$ is a cartesian category $X,$ and its
instances are the cartesian (finite product preserving) functors $X\to\Set.$
Again, with no nontrivial loose morphisms present, the instances
are no different than the modules out of $1.$
\end{example}

\begin{example}[Instances of monads]
Consider the simple double theory $\dbl{M}$ of monads 
(\cite[Theory 3.8]{cartesian-double-theories-2024}), generated
by an object $x$ equipped with a tight monad $t:x\to x.$ A
model of $\dbl{M}$ in $\Span$ is a monad $T:X\to X$ 
on some category $X.$ An instance $H$ of $T$ consists of: 
\begin{itemize}
\item A fibered set $Hx\to Xx$ underlying a functor $H:X\to \Set.$
\item A fibered set map 
\[\begin{tikzcd}
	Hx & Xx \\
	Hx & Xx
	\arrow[from=1-1, to=1-2]
	\arrow["Ht"', from=1-1, to=2-1]
	\arrow["Xt", from=1-2, to=2-2]
	\arrow[from=2-1, to=2-2]
\end{tikzcd}\]
giving the data of a natural transformation $Ht:H\to H\circ T$. 
Tight functoriality
then implies we must take $Ht^2= Ht\cdot Ht\circ T.$
\end{itemize}

Associativity and unitality of actions are encompassed by the
functoriality of $H$, so that only naturality of actions, connecting 
to the monad cells
\begin{equation*}
  \begin{tikzcd}
	x & x \\
	x & x
	\arrow[""{name=0, anchor=center, inner sep=0}, "\shortmid"{marking}, equals, from=1-1, to=1-2]
	\arrow["{t^2}"', from=1-1, to=2-1]
	\arrow["t", from=1-2, to=2-2]
	\arrow[""{name=1, anchor=center, inner sep=0}, "\shortmid"{marking}, equals, from=2-1, to=2-2]
	\arrow["\mu"{description}, draw=none, from=0, to=1]
  \end{tikzcd}
  \qquad\text{and}\qquad
  \begin{tikzcd}
	x & x \\
	x & x
	\arrow[""{name=0, anchor=center, inner sep=0}, "\shortmid"{marking}, equals, from=1-1, to=1-2]
	\arrow[equals, from=1-1, to=2-1]
	\arrow["t", from=1-2, to=2-2]
	\arrow[""{name=1, anchor=center, inner sep=0}, "\shortmid"{marking}, equals, from=2-1, to=2-2]
	\arrow["\eta"{description}, draw=none, from=0, to=1]
  \end{tikzcd},
\end{equation*}
adds any further constraint.
Naturality of actions for the unit $\eta$ says the following, when drawn in $\dbl{P}\mathsf{rof}$:
\[\begin{tikzcd}
	1 & X & X && 1 & X & X \\
	1 & X & X & {=} & 1 && X \\
	1 && X && 1 && X
	\arrow[""{name=0, anchor=center, inner sep=0}, "H"{inner sep=.8ex}, "\shortmid"{marking}, from=1-1, to=1-2]
	\arrow[equals, from=1-1, to=2-1]
	\arrow[""{name=1, anchor=center, inner sep=0}, "\shortmid"{marking}, equals, from=1-2, to=1-3]
	\arrow[equals, from=1-2, to=2-2]
	\arrow["T", from=1-3, to=2-3]
	\arrow["H"{inner sep=.8ex}, "\shortmid"{marking}, from=1-5, to=1-6]
	\arrow[equals, from=1-5, to=2-5]
	\arrow["\shortmid"{marking}, equals, from=1-6, to=1-7]
	\arrow[equals, from=1-7, to=2-7]
	\arrow[""{name=2, anchor=center, inner sep=0}, "H"'{inner sep=.8ex}, "\shortmid"{marking}, from=2-1, to=2-2]
	\arrow[equals, from=2-1, to=3-1]
	\arrow[""{name=3, anchor=center, inner sep=0}, "\shortmid"{marking}, equals, from=2-2, to=2-3]
	\arrow[equals, from=2-3, to=3-3]
	\arrow[""{name=4, anchor=center, inner sep=0}, "H"'{inner sep=.8ex}, "\shortmid"{marking}, from=2-5, to=2-7]
	\arrow[equals, from=2-5, to=3-5]
	\arrow["T", from=2-7, to=3-7]
	\arrow[""{name=5, anchor=center, inner sep=0}, "H"'{inner sep=.8ex}, "\shortmid"{marking}, from=3-1, to=3-3]
	\arrow[""{name=6, anchor=center, inner sep=0}, "H"'{inner sep=.8ex}, "\shortmid"{marking}, from=3-5, to=3-7]
	\arrow["{1_H}"{description}, draw=none, from=0, to=2]
	\arrow["\eta"{description}, draw=none, from=1, to=3]
	\arrow["H"{description}, draw=none, from=1-6, to=4]
	\arrow["H"{description}, draw=none, from=2-2, to=5]
	\arrow["Ht"{description}, draw=none, from=4, to=6]
\end{tikzcd}\]
Thus if $h\in H_a$ and $f:a\to b$ in $X,$ we have the equal elements
$h\cdot (f\cdot \eta_b)$ and $(Ht)_b(h\cdot f)$ of $H_{Tb}.$ Setting $f$ 
to an identity, we see that in fact $Ht_a$ is just action
by $\eta_a.$

Now for naturality of the action of $\mu$, we find:
\[\begin{tikzcd}
	1 & X & X && 1 & X & X \\
	1 & X & X & {=} & 1 && X \\
	1 && X && 1 && X
	\arrow[""{name=0, anchor=center, inner sep=0}, "H"{inner sep=.8ex}, "\shortmid"{marking}, from=1-1, to=1-2]
	\arrow[equals, from=1-1, to=2-1]
	\arrow[""{name=1, anchor=center, inner sep=0}, "\shortmid"{marking}, equals, from=1-2, to=1-3]
	\arrow["{T^2}"{description}, from=1-2, to=2-2]
	\arrow["T", from=1-3, to=2-3]
	\arrow["H"{inner sep=.8ex}, "\shortmid"{marking}, from=1-5, to=1-6]
	\arrow[equals, from=1-5, to=2-5]
	\arrow["\shortmid"{marking}, equals, from=1-6, to=1-7]
	\arrow[equals, from=1-7, to=2-7]
	\arrow[""{name=2, anchor=center, inner sep=0}, "H"'{inner sep=.8ex}, "\shortmid"{marking}, from=2-1, to=2-2]
	\arrow[equals, from=2-1, to=3-1]
	\arrow[""{name=3, anchor=center, inner sep=0}, "\shortmid"{marking}, equals, from=2-2, to=2-3]
	\arrow[equals, from=2-3, to=3-3]
	\arrow[""{name=4, anchor=center, inner sep=0}, "H"'{inner sep=.8ex}, "\shortmid"{marking}, from=2-5, to=2-7]
	\arrow[equals, from=2-5, to=3-5]
	\arrow["T", from=2-7, to=3-7]
	\arrow[""{name=5, anchor=center, inner sep=0}, "H"'{inner sep=.8ex}, "\shortmid"{marking}, from=3-1, to=3-3]
	\arrow[""{name=6, anchor=center, inner sep=0}, "H"'{inner sep=.8ex}, "\shortmid"{marking}, from=3-5, to=3-7]
	\arrow["{Ht^2}"{description}, draw=none, from=0, to=2]
	\arrow["\mu"{description}, draw=none, from=1, to=3]
	\arrow["H"{description}, draw=none, from=1-6, to=4]
	\arrow["H"{description}, draw=none, from=2-2, to=5]
	\arrow["Ht"{description}, draw=none, from=4, to=6]
\end{tikzcd}\]
The left-hand side says, given $h\in H_a$ and $f:a\to b$ in $X,$
to compute the action 
\[h'=H(t^2)(h)\cdot \mu_b=h\cdot \eta_a\cdot \eta_{Ta}\cdot T^2(f)\cdot \mu_b,\]
which is in $H_{Tb},$ as we see here: 
\[\begin{tikzcd}
	h &&&& {h'} \\
	a & Ta & {T^2a} & {T^2b} & Tb
	\arrow[between={0.3}{0.7}, maps to, from=1-1, to=1-5]
	\arrow["{\eta_a}", from=2-1, to=2-2]
	\arrow["{\eta_{Ta}}", from=2-2, to=2-3]
	\arrow["{T^2f}", from=2-3, to=2-4]
	\arrow["{\mu_b}", from=2-4, to=2-5]
\end{tikzcd}\] 
The right-hand side says to compute
$Ht(h\cdot f)=h\cdot f\cdot \eta_b,$ and naturality thus says $h'=h\cdot f\cdot \eta_b.$
Now, using the naturality of $\eta$ and of $\mu$ and the unitality of $T,$ we calculate 
\begin{align*}
h'&= h\cdot \eta_a\cdot T(f)\cdot \eta_{Tb}\cdot \mu_b\\
&= h\cdot \eta_a\cdot T(f)\\
&= h\cdot f\cdot \eta_b,
\end{align*}
so that this axiom is redundant. 

Therefore, we have derived a natural notion of instance of a category 
$X$ equipped with a monad $T$: it is just an functor $H:X\to\Set,$ 
considered along with the natural transformation $H\eta:H\to H\circ T.$
\end{example}

\section{The category of instances for a span-valued model}

We study the structure of the category of instances of a model in the important
case that the model is valued in $\Span$. We will show that, for a fixed
span-valued model of a simple theory, the category of instances is equivalent to
a presheaf category, and then extend the result to models of cartesian theories.
We first construct the category whose actions will turn out to coincide with
instances of a given model.

\begin{construction}[Collage of a span-valued model]\label{cons:collage}
  Given a model $X:\dbl{D}\to \Span$ of a double theory, we define its
  \define{collage} $\kappa X$ to be the following category:
  \begin{itemize}
    \item For each $d\in \dbl{D}$ and $x\in X(d),$ there is an object $[x]$ in $\kappa X.$
    \item Morphisms of $\kappa X$ are generated by two cases: 
      \begin{enumerate}
      \item For each tight morphism $f:d\to d'$ in $\dbl{D}$ and $x\in X(d),$ there is a morphism
      $[f]_x:[x]\to [Xf(x)]$ in $\kappa X.$
      \item For each loose morphism $m:d\proto d'$ in $\dbl{D}$, objects
      $x \in X(d)$ and $x'\in X(d'),$ and $h\in X(m)(x,x'),$
      there is a morphism $[h]:[x]\to [x']$ in $\kappa X.$
      \end{enumerate}
    \item Relations of $\kappa X$ are generated by four cases:
      \begin{enumerate}
        \item For composable tight morphisms $d\xto{f} d'\xto{g} d''$ in $\dbl{D}$ and
        any $x\in X(d),$ we have $[f]_x\cdot [g]_{Xf(x)}=[f\cdot g]_x$ in $\kappa X,$ while also 
        $[\id_d]_x=1_{[x]}.$ (This provides a functor $\int X_0\to \kappa X.$)
        \item For composable loose morphisms $d_0\xproto{m_1} d_1\xproto{m_2} d_2$ in $\dbl{D},$ $x_i\in X(d_i)$ and $h_i\in X(m_i)(x_{i-1},x_i),$
        the triangle below commutes in $\kappa X:$
        \[\begin{tikzcd}
          {[x_0]} & {[x_1]} \\
          & {[x_2]}
          \arrow["{[h_1]}", from=1-1, to=1-2]
          \arrow["{[X_{m_1,m_2}(h_1,h_2)]}"', from=1-1, to=2-2]
          \arrow["{[h_2]}", from=1-2, to=2-2]
        \end{tikzcd}\]
        (These relations, in case $m_1=\id_{d_1}$ or $m_2=\id_{d_1},$ provide a functor from the collage 
        of the profunctor $X(m)$ to $\kappa X.$)
        \item For each cell $\xinlinecell{d_0}{d_1}{d_2}{d_3}{f}{g}{m}{n}{\alpha}$ in $\dbl{D},$ objects $x_0 \in X(d_0)$ and $x_1 \in X(d_1),$
        and $h\in X(m)(x_0,x_1),$ the square below commutes in $\kappa X$:
        \[\begin{tikzcd}
            {[x_0]} & {[x_1]} \\
            {[Xf(x_0)]} & {[Xg(x_1)]}
            \arrow["{[h]}", from=1-1, to=1-2]
            \arrow["{[f_{x_0}]}"', from=1-1, to=2-1]
            \arrow["{[g_{x_1}]}", from=1-2, to=2-2]
            \arrow["{[X\alpha(h)]}"', from=2-1, to=2-2]
          \end{tikzcd}\]
        \item For each $d\in \dbl{D}$ and $x\in X(d),$ the morphism $[X_d(x)]:[x]\to [x]$ arising from
        application $X_d(x)\in X(\id_d)$ of the unitor
        coincides with the identity $1_{[x]}$ in $\kappa X.$
      \end{enumerate}
  \end{itemize}
\end{construction}

\begin{example}[Collage of terminal model]\label{ex:collage_terminal_model}
It is clarifying to consider the collage $\kappa \dbl{D} := \kappa I_{\dbl{D}}$ of the terminal
model mapping all objects and loose arrows to the terminal set. To wit, 
$\kappa\dbl D$ has as objects the objects of $\dbl{D}$ and as morphisms both the tight
\emph{and} the loose morphisms of $\dbl{D}$. As for relations, $\kappa\dbl D$ respects
both tight and loose composition in $\dbl{D},$ flattens cells in $\dbl{D}$ into commutative squares, 
and glues the loose and tight identities of $d\in \dbl{D}$ together into the one identity 
of $d$ in $\kappa\dbl{D}.$ For instance, an instance of the terminal model of the walking cell is simply a commutative square of functions. 

The collage thus gives a natural truncation of $\dbl{D}$ into
a 1-category, and highlights that an instance of $I_{\dbl{D}}$ is given by simultaneous
actions of the tight and the loose categories of $\dbl{D}$ on sets (truncating the loose category 
to be strictly associative if necessary), cohering with each other according to the cells of $\dbl{D}.$
\end{example}

We will prove shortly that the collage construction extends to a functor
$\kappa: \cat{Lax}(\dbl{D}, \Span) \to \cat{Cat}$, but for the moment we focus
on the object assignment.

To illustrate the construction, we compute the collage in general for models of
discrete double theories (\cref{def:discreteDoubleCategory}).

\begin{proposition}[Collages for models of discrete theories]
  \begin{enumerate} 
    \item Let $\dbl{D}=\dbl{L}C$ be the loosely discrete double theory associated to a category $C.$
    Then the collage functor $\kappa:\cat{Lax}(\dbl{D},\Span)\to \cat{Cat}$ factors as 
    the composite $\cat{Lax}(\dbl{D},\Span)\simeq \cat{Cat}/C\xto{\Sigma} \cat{Cat},$
    where $\Sigma$ is the forgetful functor from the slice.

    \item Let $\dbl{D}=\dbl{T}C$ be the tightly discrete double theory associated to a category $C.$
    Then the collage functor $\kappa:\cat{Lax}(\dbl{D},\Span)\to \cat{Cat}$ factors as 
    the composite $\cat{Lax}(\dbl{D},\Span)\simeq \CAT(C,\cat{Cat})\xto{\int} \cat{Cat}/C\xto{\Sigma} \cat{Cat}.$
  \end{enumerate}
\end{proposition}

\begin{proof}
Recall the equivalences that were given in \cref{prop:span_models_of_discrete_theories}.

\begin{enumerate}
\item In this case, the collage construction can be reduced:
case 1 of the morphism generators does not occur, and nor 
do cases 1 and 3 of the relations. Thus the collage $\kappa X$
of a model of $\dbl{L}C$ is a quotient of the category 
freely generated by all the values of the $Xd$ and the $Xm$ by 
the relations determined by the unitors and laxators of $X,$ 
which is the same construction that gives the equivalence between 
the normal lax functors on $\dbl{L}C$ and the slice $\cat{Cat}/C.$
\item In this case, type-2 morphisms in the collage occur only 
at loose identities of $\dbl{D},$ and together with the only occurrences
of type-2 and type-4 relations, we get a map $\coprod_{c\in C}Xc\to \kappa X$
for a model $X$ viewed as a normal lax functor $X:\dbl{T}C\to \Prof.$
The type-3 relations do not occur, so all that remains is to account for the 
type-1 morphisms and type-1 relations, which are exactly those that go into the Grothendieck 
construction $ \CAT(C,\cat{Cat})\xto{\int} \cat{Cat}/C.$
\end{enumerate}
\end{proof}

\begin{proposition}[Instances are of presheaf type]\label{prop:instances_presheaf_type}
  Given a model $X:\dbl{D}\to \Span$ of a simple double theory, the category of instances $\Inst(X)$ is 
  equivalent to the category of functors $\kappa X\to \Set,$ where the collage  $\kappa X$ is 
  as defined above.

  We will denote the equivalence $\Inst(X)\to \Cat(\kappa X,\Set)$ by $\bar\kappa$.
\end{proposition}
\begin{proof}
  Given an instance $H:1\proto X$ of $X,$ define $\bar\kappa(H):\kappa X\to \Set$ 
  as follows: 
  \begin{itemize}
    \item For $x\in X(d),$ send $[x]$ to the fiber $Hd_x.$
    \item For $f:d\to d'$ in $\dbl{D}$ and $x\in X(d),$ send $[f]_x:[x]\to [Xf(x)]$ to the function
    $Hd_x\to Hd'_{Xf(x)}$ given by the fiber of the span map $Hf.$ 
    \item For $m:d\proto d'$ in $\dbl{D},$ $x\in X(d),x'\in X(d'),$ and $h\in X(m)(x,x'),$ 
    send $[h]:[x]\to [x']$ to the function $Hd_x\to Hd'_{x'}$ given by the appropriate fiber of 
    the action cell $Hm$ at $(x,h,x').$
  \end{itemize}
  This operation is invertible up to isomorphism, by summing values of a 
  functor $\kappa X\to \Set$
  over various fibers. We next consider the axioms on $H$ versus $\bar\kappa H.$ 

  \begin{enumerate}
  \item $\bar\kappa H$ respects composition of arrows associated to tight
  arrows if and only if $H$ is functorial on arrows. 
  \item $\bar\kappa H$ respects composition of arrows associated to 
  loose arrows if and only if $H$ has associative actions. 
  \item $\bar\kappa H$ respects the relations associated to cells 
  if and only if $H$ has natural actions. 
  \item $\bar\kappa H$ respects the relation $[X_d(x)]=1_{[x]}$ if and only if $H$ has unital actions.
  \end{enumerate}

  Next, we consider the action of $\bar\kappa$ on morphisms. 
  Given a morphism $\mu: H\to K$ of instances of $X:\dbl{D}\to\Span$,
  define $\bar\kappa\mu:\bar\kappa H\to \bar\kappa K$ such that, 
  for $x\in X(d),$ the component 
  $\bar\kappa_{[x]}:\bar\kappa H([x])\to \bar\kappa K([x])$ is the fiber
  of $\mu_d:Hd\to Kd$ at $x.$ As on objects, this operation is 
  obviously invertible by summing over fibers. And 
  we see that $\bar\kappa\mu$ is natural at morphisms of the form 
  $[f]_x$ if and only if $\mu$ is natural, while $\bar\kappa\mu$ is 
  natural on morphisms of the form $[h]$ if and only if $\mu$ 
  is equivariant. This completes the proof.
\end{proof}

We now see what we can prove for categories of instances of models of cartesian theories. 
Observe that the naturality axiom for a modulation implies that,
if $\mu:M\to N$ is any modulation of cartesian instances, then
the component $\mu_{\Pi d_i}$ has as its $i$th component 
$M\pi_i\cdot \mu_{d_i}.$ Thus, the category of
cartesian instances $\CartInst(X)$ is a \emph{full} subcategory of the category 
$\Inst(X)$ of instances and, furthermore, the components of a cartesian
instance morphism are all determined from the components at a set of objects
spanning the theory under products. 

In fact, since the cartesian
instance axiom involves only finite products, we can deduce the following,
building on the result that the category of all instances is of presheaf type
(\cref{prop:instances_presheaf_type}):

\begin{proposition}[Cartesian instances are algebraic]\label{prop:cartesian_instances_algebraic}
The category of cartesian instances is closed under sifted colimits
and all limits. The inclusion $\CartInst(X)\to \Inst(X)$
is reflective and the category of cartesian instances is algebraic---in
particular, finitely locally presentable.
\end{proposition}
\begin{proof}
We have just shown that an instance $P$ of a model $X$ of the theory $\dbl{D}$ 
is equivalently a functor $P':\kappa(X)\to \cat{Set},$ with objects 
of $\kappa(X)$ the sum $\sum_d X(d)$ and $P'(x_d)=P(d)^{-1}(x).$

Now consider the case that $\dbl{D}$ and $X$ are cartesian, where for any 
$x_d\in X(d),y_{d'}\in X(d'),$ we have the induced element 
$x_d\leftarrow (x_d,y_{d'}) \to y_{d'}$ in $\kappa(X)$ 
arising from the projections $\pi_d:d\times d'\to d$ and 
$\pi_{d'}:d\times d'\to d'.$ Similarly, we have $X(I)\cong \{\bullet\}$ and the cone 
$(\bullet)$ over the empty diagram in $\kappa(X).$\footnote{Note that $(x_d,y_{d'})$ is 
\emph{not} actually a product in $\kappa(X)$. There might 
be loose morphisms $m:e\proto d$ and $m':e\proto d'$ producing morphisms $z_e\to x_d,y_{d'}$ with no
factorization through $(x_d,y_{d'})$, since $d\times d'$ is not a loose product.
Similarly, when there are loose arrows into $I$ in $\dbl{D}$, $\bullet$ need not be terminal
in $\kappa(X).$}

The condition that $P$ be cartesian then reduces to precisely the condition that 
$P'$ send the cones $x_d\leftarrow (x_d,y_{d'})\to y_{d'}$ and $(\bullet)$ in $\kappa(X)$ to 
product cones in $\cat{Set}.$ This class of cones above equips $\kappa(X)$ with the 
structure of a finite product sketch, whose category of models 
is the category of cartesian instances of $X.$ Then 
the result follows from the facts that categories of models of finite
product sketches are closed under sifted colimits and all limits in the 
copresheaf category, are cocomplete, and from the universal property of 
the copresheaf category as a free cocompletion. All these results may be 
found in Ad\'amek, Rosický, and Vitale's \cite{algebraic-theories-2010}, 
especially Theorem 4.5 for the cocompleteness.
\end{proof}

\begin{remark}[Other semantic double categories]
For models valued in semantic double categories $\dbl{S}$ besides $\Span$,
the goal of generalizing $\kappa$ would be to produce a universal 
category object in $\dbl{S}$ associated to a model $X:\dbl{D}\to \dbl{S}.$
Then we might aim to show that the category of instances of $X$ 
is equivalent to a category of internal copresheaves on $\kappa X.$ 
We leave this generalization to future work, in particular as we are not
aware of any definition of internal copresheaf on a category object 
in an arbitrary double category yet extant. A more straightforward generalization
would be to work in $\dbl{S}=\Span(\cat{E})$, thus generalizing to classical
internal category theory.
\end{remark}

We next establish the functoriality of $\kappa.$

\begin{proposition}[Collage of a model is a functor]
The collage $\kappa X$ of a model $X:\dbl{D}\to \Span$ of a simple double theory $\dbl{D}$
extends to a functor $\kappa:\cat{Lax}(\dbl{D},\Span)\to \cat{Cat}$.
\end{proposition}
\begin{proof}
Consider a morphism $\alpha:X\to Y$ of span valued-models of $\dbl{D}.$ We define
$\kappa \alpha: \kappa X \to \kappa Y$ on objects as $\kappa\alpha([x]) \coloneqq [\alpha_d(x)]$ for $x\in X(d).$

Now for the two classes of generating morphisms in $\kappa X.$ 
\begin{enumerate}
\item Given 
a morphism of the form $[f]_x:[x]\to [Xf(x)]$ in $\kappa X,$ we must define 
$\kappa\alpha([f]_x):[\alpha_d(x)]\to [\alpha_{d'}(Xf(x))].$ Recalling that 
$\alpha_{d'}(Xf(x))=Yf(\alpha_d(x)),$ via naturality, it type-checks to 
set $\kappa\alpha([f]_x) \coloneqq [f]_{\alpha_d(x)},$ and we proceed to do so.
(Note that the use of strict 
naturality here highlights that we do not expect functoriality of $\kappa$ in
non-strict transformations of models.)

\item Given a morphism of the form $[h]:[x]\to [x']$ in $\kappa X,$ we will define
$\kappa\alpha([h]) \coloneqq [\alpha_m(h)]:[\alpha_d(x)]\to [\alpha_{d'}(x')].$
\end{enumerate}

We now check that the graph map $\kappa\alpha$ just defined respects the relations in $\kappa X.$ 
\begin{enumerate}
\item Given composable tight morphisms $d\xto{f} d'\xto{g} d''$ in $\dbl{D}$ and $x\in X(d),$ we have
\begin{align*}
\kappa\alpha([f]_x)\cdot \kappa\alpha([g]_{Xf(x)})&=[f]_{\alpha_d(x)}\cdot [g]_{\alpha_{d'}(Xf(x))}\\
&=[f]_{\alpha_d(x)}\cdot [g]_{Yf(\alpha_d(x))}\\
&= [f\cdot g]_{\alpha_d(x)}\\
&=\kappa\alpha([f\cdot g]_x),
\end{align*}
applying naturality of $\alpha$ and the tight-arrows relation in $\kappa Y.$ 
\item Given a cell $\xinlinecell{d_0}{d_1}{d_2}{d_3}{f}{g}{m}{n}{\gamma}$ in $\dbl{D}$,
  $x_0 \in X(d_0)$, $x_1\in X(d_1)$, and $h\in X(m)(x_0,x_1),$ we have
  \begin{align*}
  \kappa\alpha([h])\cdot \kappa\alpha([g]_{x_1})&= [\alpha_m(h)]\cdot [g]_{\alpha_{d_1}(x_1)}\\
  &= [f]_{\alpha_{d_0}(x_0)}\cdot [Y_\gamma(\alpha_m(h))]\\
  &= [f]_{\alpha_{d_0}(x_0)}\cdot [\alpha_n(X_\gamma(h))]\\
  &= \kappa\alpha([f]_{x_0})\cdot \kappa\alpha([X_\gamma(h)]),
  \end{align*}
  where we use the relation in $\kappa Y$ arising from $\alpha$ 
  as well as the naturality of $\alpha$ at $\gamma.$
\item Given composable loose morphisms $d_0\xproto{m_1} d_1\xproto{m_2} d_2$ in $\dbl{D},$ $x_i\in X(d_i)$ and $h_i\in X(m_i)(x_{i-1},x_i),$ we have
\begin{align*}
\kappa\alpha([h_1])\cdot \kappa\alpha([h_2])&=[\alpha_{m_1}(h_1)]\cdot [\alpha_{m_2}(h_2)]\\
&=[Y_{m_1,m_2}(\alpha_{m_1}(h_1),\alpha_{m_2}(h_2))]\\
&=[\alpha_{m_1\odot m_2}(X_{m_1,m_2}(h_1,h_2))]\\
&=\kappa\alpha([X_{m_1,m_2}(h_1,h_2)]),
\end{align*}
where we use relation arising from $m_1,m_2$ in $\kappa Y$ and respect of $\alpha$ for laxators.
\item Given $d\in \dbl{D}$ and $x\in X(d),$ we have 
\[\kappa\alpha([X_d(x)])=[\alpha_{\id_d}(X_d(x))]=[Y_d(\alpha_d(x))]=1_{[\alpha_d(x)]}.\]
\end{enumerate}
This establishes that $\kappa\alpha$ provides a functor $\kappa X\to \kappa Y,$ 
and it is immediate that $\kappa$ itself respects composition, which proves the result. 
\end{proof}

We can now fully establish the functoriality of $\Inst(X)$ in $X$, and more:

\begin{corollary}[Functoriality of span-valued instances]\label{cor:span_instances_functoriality}
  For any morphism $\alpha:X\to Y$ of span-valued models of a simple double
  theory, there is an adjoint triple
  \[\alpha_!\dashv \alpha^* \dashv \alpha_*:\Inst(Y)\to \Inst(X),\] with 
  $\alpha^*$ as defined in \cref{prop:instances_functoriality}. If $X$ and $Y$
  are models of a \emph{cartesian} theory, then
  $\alpha^*:\CartInst(Y)\to \CartInst(X)$ has a left adjoint
  $\alpha_!.$
\end{corollary}
\begin{proof}
Under the equivalence $\Inst(X)\cong \Cat(\kappa X,\Set),$ this triple is simply the 
correlate of the Kan extension triple $\kappa\alpha_!\dashv \kappa\alpha^*\dashv \kappa\alpha_*,$
and one checks directly that $\kappa\alpha^*$ behaves as in the prior definition of $\alpha^*.$
For the cartesian case, this follows from presentability of 
$\CartInst(X)$ and the fact that it is closed under filtered colimits
and limits in $\Inst(X)$, so that the restricted $\alpha^*$ is 
accesible and has a left adjoint. Note that in the cartesian case, $\alpha^*$
generally lacks a right adjoint. 
\end{proof}


\begin{remark}[$\kappa$ is neither full, faithful, conservative, continuous, nor cocontinuous]\label{rem:kappa-bad}
As a warning, the collage functor $\kappa$ is not particularly nice.
For instance, a lax functor
$\dbl{D}\to \Span$, where $\dbl{D}$ is freely generated by a spherical cell 
$\xinlinecell{x}{x}{x}{x}{\id}{\id}{\id}{\id}{\alpha},$ may be identified with a category $\cat{C}$ equipped
with a natural endomorphism $q$ of its identity functor. Then $\kappa$ will produce the category that 
identifies each morphism of $\cat{C}$ with its image under $q.$ 

Clearly $\kappa$ then identifies 
$(\cat{C},q)$ with $(\kappa\cat{C},\id).$ For instance, if $\cat{C}=\mathbf{B}C$ is an abelian group 
viewed as a one-object category, then $q$ can be any element of $C,$ and 
$\kappa(\cat{C},q)=\mathbf{B}(C/\langle q\rangle).$ 
Now we see there is no map $(\mathbf{B}\mathbb{Z}/n,0)\to (\mathbf{B}\mathbb{Z},n)$, while the 
unique map in the other direction is sent to an isomorphism, showing $\kappa$ is neither full nor 
conservative. There are, furthermore, many maps $(\mathbf{B}(\mathbb Z^2),(1,0))\to (\mathbf B \mathbb Z,1)$, 
but only one map between the images of these under $\kappa,$ so $\kappa$ is not faithful either.

As for continuity, \cref{ex:collage_terminal_model} shows that $\kappa$ does not send the terminal
model to the terminal category, unless $\dbl{D}$ is itself terminal. 
Cocontinuity is a bit more exciting. In fact, $\kappa X$ ought to be thought of as the 
lax colimit of the model $X.$
Thus we cannot expect a right adjoint to $\kappa$ landing in strict 
transformations,\footnote{One might imagine applying a lax transformation coclassifier functor, 
but having a model coclassify transformations that are lax in the tight 
direction is inconsistent with the lack of effect of this laxity on the loose direction. 
For instance, for such a theory as the $\dbl{D}$ freely generated by a tight and by a loose morphism,
one can check directly that a lax transformation coclassifier of a model $C\leftarrow A\proto B$ 
must both leave $A$ be and replace it with its comma over $C$.} but
one can think of a universal property of $\kappa$ in terms of \emph{lax} transformations. 
We will execute this idea just below, in \cref{prop:collage_universal_property}.

Thus on the one hand $\kappa$ is of little use in studying properties of the category 
$\Lax(\dbl{D},\Span)$, while on the other we see that the category of instances of a 
model does not depend on nearly all the data of the model in general. 
Nonetheless we shall show below (\cref{prop:models_locally_presentable})
that the category of models $\Lax(\dbl{D},\Span)$ is indeed 
very well-behaved, namely, locally finitely presentable.
\end{remark}

\subsection{Universal property of the collage}

We conclude this section by giving a universal property of $\kappa,$ in terms of lax 
transformations (with components given by \emph{tight} arrows), 
a structure not otherwise considered in this paper. We reproduce the 
definition here almost exactly as in \cite[Definition 7.1]{cartesian-double-theories-2024};
there is very limited opportunity to compress the definition as for stricter 
structures, since none of the components of a lax natural transformation form 
either a functor or a natural transformation on the nose. Note well that the naturality 
comparison has loose codomain $G(\id)$ rather than $\id_G,$ a softening analogous 
to the distinction between modifications and modulations of \emph{loose} 
transformations, see Appendix \ref{app:loose-transformations}. 

\begin{definition}\label{def:lax_tight_natural_transformation}
Let $F,G:\dbl{D}\to \dbl{E}$ be lax double functors. 
A \define{lax transformation} $\alpha:F\to G$ consists of the following data:
\begin{itemize}
\item For each object $d\in \dbl{D}$, a tight morphism $\alpha_d:Fd\to Gd$ in $\dbl{E}$.
\item For each tight morphism $f:d\to d'$ in $\dbl{D}$, a cell
\[\begin{tikzcd}
	Fd & Fd \\
	Gd & {Fd'} \\
	{Gd'} & {Gd'}
	\arrow[""{name=0, anchor=center, inner sep=0}, "\shortmid"{marking}, equals, from=1-1, to=1-2]
	\arrow["{\alpha_d}"', from=1-1, to=2-1]
	\arrow["Ff", from=1-2, to=2-2]
	\arrow["Gf"', from=2-1, to=3-1]
	\arrow["{\alpha_{d'}}", from=2-2, to=3-2]
	\arrow[""{name=1, anchor=center, inner sep=0}, "{G(\mathrm{id}_{d'})}"'{inner sep=.8ex}, "\shortmid"{marking}, from=3-1, to=3-2]
	\arrow["{\alpha_f}"{description}, draw=none, from=0, to=1]
\end{tikzcd}\]
\item For every loose morphism $m:d\proto d'$ in $\dbl{D}$, a cell
\[\begin{tikzcd}
	Fd & {Fd'} \\
	Gd & {Gd'}
	\arrow[""{name=0, anchor=center, inner sep=0}, "Fm"{inner sep=.8ex}, "\shortmid"{marking}, from=1-1, to=1-2]
	\arrow["{\alpha_d}"', from=1-1, to=2-1]
	\arrow["{\alpha_{d'}}", from=1-2, to=2-2]
	\arrow[""{name=1, anchor=center, inner sep=0}, "Gm"'{inner sep=.8ex}, "\shortmid"{marking}, from=2-1, to=2-2]
	\arrow["{\alpha_m}"{description}, draw=none, from=0, to=1]
\end{tikzcd}\]
\end{itemize}

such that: 
\begin{itemize}
\item $\alpha_m$ respects loose composition in that the following squares commute, 
in the categories of functors and \emph{unnatural} transformations marked below each 
square: 
\[\begin{tikzcd}
	{F_1\odot F_1} & {F_1(-\odot -)} & {\mathrm{id}_{F-}} & {F(\mathrm{id}_{-})} \\
	{G_1\odot G_1} & {G_1(-\odot -)} & {\mathrm{id}_{G-}} & {G(\mathrm{id}_{-})} \\
	{\dbl{D}_1\times_{\dbl{D}_0}\dbl{D}_1} & {\dbl{E}_1} & {\dbl{D}_0} & {\dbl{E}_1}
	\arrow["{F_{-,-}}", from=1-1, to=1-2]
	\arrow["{\alpha\odot \alpha}"', from=1-1, to=2-1]
	\arrow["{\alpha_{-\odot -}}", from=1-2, to=2-2]
	\arrow["{F_{-}}", from=1-3, to=1-4]
	\arrow["{\mathrm{id}_{\alpha_{-}}}"', from=1-3, to=2-3]
	\arrow["{\alpha_{\mathrm{id}_{-}}}", from=1-4, to=2-4]
	\arrow["{G_{-,-}}"', from=2-1, to=2-2]
	\arrow["{G_{-}}"', from=2-3, to=2-4]
	\arrow[from=3-1, to=3-2]
	\arrow[from=3-3, to=3-4]
\end{tikzcd}\]
\item $\alpha_m$ is approximately natural in $m:$ given a cell 
$\xinlinecell{d_0}{d_1}{d_2}{d_3}{f}{g}{m}{n}{\gamma}$ in $\dbl{D}$, we have 
\[\begin{tikzcd}
	{Fd_0} & {Fd_0} & {Fd_1} && {Fd_0} & {Fd_1} & {Fd_1} \\
	{Gd_0} & {Fd_2} & {Fd_3} && {Gd_0} & {Gd_1} & {Fd_3} \\
	{Gd_2} & {Gd_2} & {Gd_3} && {Gd_2} & {Gd_3} & {Gd_3} \\
	{Gd_2} && {Gd_3} && {Gd_2} && {Gd_3}
	\arrow[""{name=0, anchor=center, inner sep=0}, "\shortmid"{marking}, equals, from=1-1, to=1-2]
	\arrow["{\alpha_{d_0}}"', from=1-1, to=2-1]
	\arrow[""{name=1, anchor=center, inner sep=0}, "Fm"{inner sep=.8ex}, "\shortmid"{marking}, from=1-2, to=1-3]
	\arrow["Ff"{description}, from=1-2, to=2-2]
	\arrow["Fg", from=1-3, to=2-3]
	\arrow[""{name=2, anchor=center, inner sep=0}, "Fm"{inner sep=.8ex}, "\shortmid"{marking}, from=1-5, to=1-6]
	\arrow["{\alpha_{d_0}}"', from=1-5, to=2-5]
	\arrow[""{name=3, anchor=center, inner sep=0}, "\shortmid"{marking}, equals, from=1-6, to=1-7]
	\arrow["{\alpha_{d_1}}"{description}, from=1-6, to=2-6]
	\arrow["Fg", from=1-7, to=2-7]
	\arrow["Gf"', from=2-1, to=3-1]
	\arrow[""{name=4, anchor=center, inner sep=0}, "Fn"'{inner sep=.8ex}, "\shortmid"{marking}, from=2-2, to=2-3]
	\arrow["{\alpha_{d_2}}"{description}, from=2-2, to=3-2]
	\arrow[""{name=5, anchor=center, inner sep=0}, "{\alpha_{d_3}}", from=2-3, to=3-3]
	\arrow[""{name=6, anchor=center, inner sep=0}, "Gm"'{inner sep=.8ex}, "\shortmid"{marking}, from=2-5, to=2-6]
	\arrow[""{name=7, anchor=center, inner sep=0}, "Gf"', from=2-5, to=3-5]
	\arrow["Gg"{description}, from=2-6, to=3-6]
	\arrow["{\alpha_{d_3}}", from=2-7, to=3-7]
	\arrow[""{name=8, anchor=center, inner sep=0}, "{G(\mathrm{id}_{d_2})}"'{inner sep=.8ex}, "\shortmid"{marking}, from=3-1, to=3-2]
	\arrow[equals, from=3-1, to=4-1]
	\arrow[""{name=9, anchor=center, inner sep=0}, "Gn"'{inner sep=.8ex}, "\shortmid"{marking}, from=3-2, to=3-3]
	\arrow[equals, from=3-3, to=4-3]
	\arrow[""{name=10, anchor=center, inner sep=0}, "Gn"'{inner sep=.8ex}, "\shortmid"{marking}, from=3-5, to=3-6]
	\arrow[equals, from=3-5, to=4-5]
	\arrow[""{name=11, anchor=center, inner sep=0}, "{G(\mathrm{id}_{d_3})}"'{inner sep=.8ex}, "\shortmid"{marking}, from=3-6, to=3-7]
	\arrow[equals, from=3-7, to=4-7]
	\arrow[""{name=12, anchor=center, inner sep=0}, "Gn"'{inner sep=.8ex}, "\shortmid"{marking}, from=4-1, to=4-3]
	\arrow[""{name=13, anchor=center, inner sep=0}, "Gn"'{inner sep=.8ex}, "\shortmid"{marking}, from=4-5, to=4-7]
	\arrow["{\alpha_f}"{description}, draw=none, from=0, to=8]
	\arrow["{F\gamma}"{description}, draw=none, from=1, to=4]
	\arrow["{\alpha_m}"{description}, draw=none, from=2, to=6]
	\arrow["{\alpha_g}"{description}, draw=none, from=3, to=11]
	\arrow["{\alpha_n}"{description}, draw=none, from=4, to=9]
	\arrow[between={0.4}{0.6}, equals, from=5, to=7]
	\arrow["{G\gamma}"{description, pos=0.6}, draw=none, from=6, to=10]
	\arrow["{G_{d_2,n}}"{description}, draw=none, from=3-2, to=12]
	\arrow["{G_{n,d_3}}"{description}, draw=none, from=3-6, to=13]
\end{tikzcd}\]
\item $\alpha_f$ is approximately functorial in $f$ in that the following pastings 
coincide whenever they are of valid type: 
\[\begin{tikzcd}
	{Fd_0} & {Fd_0} & {Fd_0} && {Fd_0} & {Fd_0} & \\
	{Gd_0} & {Fd_1} & {Fd_1} && {Gd_0} & {Fd_2} \\
	{Gd_1} & {Gd_1} & {Fd_2} && {Gd_2} & {Gd_2} \\
	{Gd_2} & {Gd_2} & {Gd_2} & Fd & Fd & Fd & Fd \\
	{Gd_2} && {Gd_2} & Gd & Gd & Gd & Fd \\
	&&& Gd & Gd & Gd & Gd
	\arrow[""{name=0, anchor=center, inner sep=0}, equals, from=1-1, to=1-2]
	\arrow["{{\alpha_{d_0}}}"', from=1-1, to=2-1]
	\arrow[equals, from=1-2, to=1-3]
	\arrow["Ff"{description}, from=1-2, to=2-2]
	\arrow["Ff", from=1-3, to=2-3]
	\arrow[""{name=1, anchor=center, inner sep=0}, equals, from=1-5, to=1-6]
	\arrow["{{\alpha_{d_0}}}"', from=1-5, to=2-5]
	\arrow["{{F(f\cdot g)}}", from=1-6, to=2-6]
	\arrow["Gf"', from=2-1, to=3-1]
	\arrow[""{name=2, anchor=center, inner sep=0}, equals, from=2-2, to=2-3]
	\arrow["{{\alpha_{d_1}}}"{description}, from=2-2, to=3-2]
	\arrow[between={0.3}{0.7}, equals, from=2-3, to=2-5]
	\arrow["Fg", from=2-3, to=3-3]
	\arrow["{{G(f\cdot g)}}"', from=2-5, to=3-5]
	\arrow["{{\alpha_{d_2}}}", from=2-6, to=3-6]
	\arrow[""{name=3, anchor=center, inner sep=0}, "{{G(\mathrm{id}_{d_1})}}"'{inner sep=.8ex}, "\shortmid"{marking}, from=3-1, to=3-2]
	\arrow["Gg"', from=3-1, to=4-1]
	\arrow["Gg"{description}, from=3-2, to=4-2]
	\arrow["{{\alpha_{d_2}}}", from=3-3, to=4-3]
	\arrow[""{name=4, anchor=center, inner sep=0}, "{{G(\mathrm{id}_{d_2})}}"'{inner sep=.8ex}, "\shortmid"{marking}, from=3-5, to=3-6]
	\arrow[""{name=5, anchor=center, inner sep=0}, "{{G(\mathrm{id}_{d_2})}}"'{inner sep=.8ex}, "\shortmid"{marking}, from=4-1, to=4-2]
	\arrow[equals, from=4-1, to=5-1]
	\arrow[""{name=6, anchor=center, inner sep=0}, "{{G(\mathrm{id}_{d_2})}}"'{inner sep=.8ex}, "\shortmid"{marking}, from=4-2, to=4-3]
	\arrow[equals, from=4-3, to=5-3]
	\arrow[""{name=7, anchor=center, inner sep=0}, "\shortmid"{marking}, equals, from=4-4, to=4-5]
	\arrow["{{\alpha_d}}"', from=4-4, to=5-4]
	\arrow["{{\alpha_d}}", from=4-5, to=5-5]
	\arrow[""{name=8, anchor=center, inner sep=0}, "\shortmid"{marking}, equals, from=4-6, to=4-7]
	\arrow["{{\alpha_d}}"', from=4-6, to=5-6]
	\arrow[equals, from=4-7, to=5-7]
	\arrow[""{name=9, anchor=center, inner sep=0}, "{{G(\mathrm{id}_{d_2})}}"'{inner sep=.8ex}, "\shortmid"{marking}, from=5-1, to=5-3]
	\arrow[""{name=10, anchor=center, inner sep=0}, "\shortmid"{marking}, equals, from=5-4, to=5-5]
	\arrow[equals, from=5-4, to=6-4]
	\arrow[between={0.1}{0.9}, equals, from=5-5, to=5-6]
	\arrow[equals, from=5-5, to=6-5]
	\arrow[equals, from=5-6, to=6-6]
	\arrow["{{\alpha_d}}", from=5-7, to=6-7]
	\arrow[""{name=11, anchor=center, inner sep=0}, "{{G(\mathrm{id}_d)}}"'{inner sep=.8ex}, "\shortmid"{marking}, from=6-4, to=6-5]
	\arrow[""{name=12, anchor=center, inner sep=0}, "{{G(\mathrm{id}_d)}}"'{inner sep=.8ex}, "\shortmid"{marking}, from=6-6, to=6-7]
	\arrow["{{\alpha_f}}"{description}, draw=none, from=0, to=3]
	\arrow["{{\alpha_{f\cdot g}}}"{description}, draw=none, from=1, to=4]
	\arrow["{{\alpha_g}}"{description}, draw=none, from=2, to=6]
	\arrow["{{G(\mathrm{id}_g)}}"{description, pos=0.6}, draw=none, from=3, to=5]
	\arrow["{{G_{d_2,d_2}}}"{description}, draw=none, from=4-2, to=9]
	\arrow["{{\mathrm{id}_{\alpha_d}}}"{description}, draw=none, from=7, to=10]
	\arrow["{{\alpha_{1_d}}}"{description}, draw=none, from=8, to=12]
	\arrow["{G_d}"{description}, draw=none, from=10, to=11]
\end{tikzcd}\]
\end{itemize}
\end{definition}

Modulations between lax transformations are a bit more complicated than those between tight transformations, 
but since we will only need them with identity modules on the loose boundary, there is a compensating 
simplification, in that we use $\id$ on top in the definition below: 
\begin{definition}[Modulations of lax transformations]\label{def:modulation_lax}
A \define{modulation} $\mu:\alpha \Rrightarrow\beta:F\To G:\dbl{D}\to \dbl{E}$ between 
lax transformations $\alpha$ and $\beta $ consists of, for every object $d\in \dbl{D},$ 
a cell as below: 
\[\begin{tikzcd}
	Fd & Fd \\
	Gd & Gd
	\arrow[""{name=0, anchor=center, inner sep=0}, "\shortmid"{marking}, equals, from=1-1, to=1-2]
	\arrow["{\alpha_d}"', from=1-1, to=2-1]
	\arrow["{\beta_d}", from=1-2, to=2-2]
	\arrow[""{name=1, anchor=center, inner sep=0}, "{G(\mathrm{id}_d)}"'{inner sep=.8ex}, "\shortmid"{marking}, from=2-1, to=2-2]
	\arrow["{\mu_d}"{description}, draw=none, from=0, to=1]
\end{tikzcd}\]
The cells $\mu_d$ are subject to the conditions: 
\begin{itemize}
\item $\mu_d$ respects loose composition in that the two pastings below coincide whenever they are of valid type:
\[\begin{tikzcd}
	Fd & Fd & {Fd'} & Fd & {Fd'} & {Fd'} \\
	Gd & Gd & {Gd'} & Gd & {Gd'} & {Gd'} \\
	Gd && {Gd'} & Gd && {Gd'}
	\arrow[""{name=0, anchor=center, inner sep=0}, "\shortmid"{marking}, equals, from=1-1, to=1-2]
	\arrow["{\alpha_d}"', from=1-1, to=2-1]
	\arrow[""{name=1, anchor=center, inner sep=0}, "Fm"{inner sep=.8ex}, "\shortmid"{marking}, from=1-2, to=1-3]
	\arrow["{\beta_d}"', from=1-2, to=2-2]
	\arrow["{\beta_{d'}}", from=1-3, to=2-3]
	\arrow[""{name=2, anchor=center, inner sep=0}, "Fm"{inner sep=.8ex}, "\shortmid"{marking}, from=1-4, to=1-5]
	\arrow["{\alpha_d}"', from=1-4, to=2-4]
	\arrow[""{name=3, anchor=center, inner sep=0}, "\shortmid"{marking}, equals, from=1-5, to=1-6]
	\arrow["{\alpha_{d'}}"', from=1-5, to=2-5]
	\arrow["{\beta_{d'}}", from=1-6, to=2-6]
	\arrow[""{name=4, anchor=center, inner sep=0}, "{G(\mathrm{id}_d)}"'{inner sep=.8ex}, "\shortmid"{marking}, from=2-1, to=2-2]
	\arrow[equals, from=2-1, to=3-1]
	\arrow[""{name=5, anchor=center, inner sep=0}, "Gm"'{inner sep=.8ex}, "\shortmid"{marking}, from=2-2, to=2-3]
	\arrow[between={0.1}{0.9}, equals, from=2-3, to=2-4]
	\arrow[equals, from=2-3, to=3-3]
	\arrow[""{name=6, anchor=center, inner sep=0}, "Gm"'{inner sep=.8ex}, "\shortmid"{marking}, from=2-4, to=2-5]
	\arrow[equals, from=2-4, to=3-4]
	\arrow[""{name=7, anchor=center, inner sep=0}, "{G(\mathrm{id}_{d'})}"'{inner sep=.8ex}, "\shortmid"{marking}, from=2-5, to=2-6]
	\arrow[equals, from=2-6, to=3-6]
	\arrow[""{name=8, anchor=center, inner sep=0}, "Gm"'{inner sep=.8ex}, "\shortmid"{marking}, from=3-1, to=3-3]
	\arrow[""{name=9, anchor=center, inner sep=0}, "Gm"'{inner sep=.8ex}, "\shortmid"{marking}, from=3-4, to=3-6]
	\arrow["{\mu_d}"{description}, shift right, draw=none, from=0, to=4]
	\arrow["{\beta_m}"{description}, draw=none, from=1, to=5]
	\arrow["{\alpha_m}"{description}, shift right, draw=none, from=2, to=6]
	\arrow["{\mu_{d'}}"{description}, draw=none, from=3, to=7]
	\arrow["{G_{d,m}}"{description}, draw=none, from=2-2, to=8]
	\arrow["{G_{m,d'}}"{description}, draw=none, from=2-5, to=9]
\end{tikzcd}\]
Either side of this equation can thus be denoted $\mu_m,$ which shows how one recovers a modulation as in 
\cref{def:modulation} from a modulation as in the definition above--for the other direction, one whiskers 
$\mu_{\id_d}$ with $F_d.$ 
\item $\mu_m$ is approximately natural in $m$ in the sense that for every cell 
$\xinlinecell{d_0}{d_1}{d_2}{d_3}{f}{g}{m}{n}{\gamma}$ in $\dbl{D},$ the two pastings below coincide: 
\[\begin{tikzcd}
	{Fd_0} & {Fd_0} & {Fd_1} & {Fd_0} & {Fd_1} & {Fd_1} \\
	{Gd_0} & {Fd_2} & {Fd_3} & {Gd_0} & {Gd_1} & {Fd_3} \\
	{Gd_2} & {Gd_2} & {Gd_3} & {Gd_2} & {Gd_3} & {Gd_3} \\
	{Gd_2} && {Gd_3} & {Gd_2} && {Gd_3}
	\arrow[""{name=0, anchor=center, inner sep=0}, "\shortmid"{marking}, equals, from=1-1, to=1-2]
	\arrow["{\alpha_{d_0}}"', from=1-1, to=2-1]
	\arrow[""{name=1, anchor=center, inner sep=0}, "Fm"{inner sep=.8ex}, "\shortmid"{marking}, from=1-2, to=1-3]
	\arrow["Ff"', from=1-2, to=2-2]
	\arrow["Fg", from=1-3, to=2-3]
	\arrow[""{name=2, anchor=center, inner sep=0}, "Fm"{inner sep=.8ex}, "\shortmid"{marking}, from=1-4, to=1-5]
	\arrow["{\alpha_{d_1}}"', from=1-4, to=2-4]
	\arrow[""{name=3, anchor=center, inner sep=0}, "\shortmid"{marking}, equals, from=1-5, to=1-6]
	\arrow["{\beta_{d_1}}"', from=1-5, to=2-5]
	\arrow["Fg", from=1-6, to=2-6]
	\arrow["Gf"', from=2-1, to=3-1]
	\arrow[""{name=4, anchor=center, inner sep=0}, "Fn"'{inner sep=.8ex}, "\shortmid"{marking}, from=2-2, to=2-3]
	\arrow["{\alpha_{d_2}}"', from=2-2, to=3-2]
	\arrow[""{name=5, anchor=center, inner sep=0}, "{\beta_{d_3}}", from=2-3, to=3-3]
	\arrow[""{name=6, anchor=center, inner sep=0}, "Gm"'{inner sep=.8ex}, "\shortmid"{marking}, from=2-4, to=2-5]
	\arrow[""{name=7, anchor=center, inner sep=0}, "Gf"', from=2-4, to=3-4]
	\arrow["Gg"', from=2-5, to=3-5]
	\arrow["{\beta_{d_3}}", from=2-6, to=3-6]
	\arrow[""{name=8, anchor=center, inner sep=0}, "{G(\mathrm{id}_{d_2})}"'{inner sep=.8ex}, "\shortmid"{marking}, from=3-1, to=3-2]
	\arrow[equals, from=3-1, to=4-1]
	\arrow[""{name=9, anchor=center, inner sep=0}, "Gn"'{inner sep=.8ex}, "\shortmid"{marking}, from=3-2, to=3-3]
	\arrow[equals, from=3-3, to=4-3]
	\arrow[""{name=10, anchor=center, inner sep=0}, "Gn"'{inner sep=.8ex}, "\shortmid"{marking}, from=3-4, to=3-5]
	\arrow[equals, from=3-4, to=4-4]
	\arrow[""{name=11, anchor=center, inner sep=0}, "{G(\mathrm{id}_{d_3})}"'{inner sep=.8ex}, "\shortmid"{marking}, from=3-5, to=3-6]
	\arrow[equals, from=3-6, to=4-6]
	\arrow[""{name=12, anchor=center, inner sep=0}, "Gn"'{inner sep=.8ex}, "\shortmid"{marking}, from=4-1, to=4-3]
	\arrow[""{name=13, anchor=center, inner sep=0}, "Gn"'{inner sep=.8ex}, "\shortmid"{marking}, from=4-4, to=4-6]
	\arrow["{\alpha_f}"{description}, draw=none, from=0, to=8]
	\arrow["{F\gamma}"{description}, draw=none, from=1, to=4]
	\arrow["{\mu_m}"{description}, shift right, draw=none, from=2, to=6]
	\arrow["{\beta_g}"{description}, draw=none, from=3, to=11]
	\arrow["{\mu_n}"{description}, draw=none, from=4, to=9]
	\arrow[between={0.4}{0.6}, equals, from=5, to=7]
	\arrow["{G\gamma}"{description}, shift right, draw=none, from=6, to=10]
	\arrow["{G_{d_2,n}}"{description}, draw=none, from=3-2, to=12]
	\arrow["{G_{n,d_3}}"{description}, draw=none, from=3-5, to=13]
\end{tikzcd}\]
Note that, in this context, one can read the approximate naturality of $\alpha_m$ in $m,$ for the definition of 
a lax transformation $\alpha,$ as the requirement that the components $\alpha_m$ themselves assemble into an identity
endomodulation of $\alpha.$
\end{itemize}
\end{definition}

Note that for any small category $C,$ we can define the constant lax functor $\Delta C:\dbl{D}\to \Span$ 
valued at $C$ by sending every object of $\dbl{D}$ to $C_0$, every tight arrow of $\dbl{D}$ to 
$\id_{C_0},$ every loose arrow of $\dbl{D}$ to the span of arrows $C_0\leftarrow C_1\to C_0,$ 
every 2-cell of $\dbl{D}$ to the identity map on this span, with laxators given by the composition 
of $C$ and unitors given by the identity arrow-assigning function of $C.$ We can use these 
constant functors to establish $\kappa X$ as a lax colimit of $X.$

\begin{proposition}[Universal property of collage]\label{prop:collage_universal_property}
Given a lax functor $X:\dbl{D}\to \Span$,
there is a canonical lax transformation $\eta: X\to \Delta\kappa X,$ composition with which 
induces an equivalence of categories $\Cat(\kappa X,C)\simeq \LaxTrans(X,\Delta C),$ where the codomain 
is the category of lax transformations and modulations between them. 
\end{proposition}
\begin{proof}
We first construct $\eta.$ We define $\eta_d:Xd\to \kappa X_0\cong \sum_{d'} Xd'$ to be the 
inclusion of the $d$ component. Given $f:d\to d',$ for every $x\in Xd,$ in $\eta_f$ we must 
specify a morphism in $\kappa X(x,Xf(x)),$ as illustrated below. We take $[f]_x.$
\[\begin{tikzcd}
	Xd & Xd & Xd \\
	{\sum_{d}Xd} && {Xd'} \\
	{\sum_{d}Xd} & {\kappa X_1} & {\sum_{d}Xd}
	\arrow["{\mathrm{in}_d}"', from=1-1, to=2-1]
	\arrow[equals, from=1-2, to=1-1]
	\arrow[equals, from=1-2, to=1-3]
	\arrow["{\eta_f}"{description}, from=1-2, to=3-2]
	\arrow["Xf", from=1-3, to=2-3]
	\arrow[equals, from=2-1, to=3-1]
	\arrow["{\mathrm{in}_{d'}}", from=2-3, to=3-3]
	\arrow["s"', from=3-2, to=3-1]
	\arrow["t", from=3-2, to=3-3]
\end{tikzcd}\]
Given $m:d\to d',$ for every $h:x\proto x'$ in $Xm,$ we must specify a morphism in 
$\kappa X(x,x'),$ as illustrated below. We take $[h].$
\[\begin{tikzcd}
	Xd & Xm & {Xd'} \\
	{\sum_d X_d} & {\kappa X_1} & {\sum_d X_d}
	\arrow["{\mathrm{in}_d}"', from=1-1, to=2-1]
	\arrow[from=1-2, to=1-1]
	\arrow[from=1-2, to=1-3]
	\arrow["{\eta_m}"', from=1-2, to=2-2]
	\arrow["{\mathrm{in}_{d'}}"', from=1-3, to=2-3]
	\arrow[from=2-2, to=2-1]
	\arrow[from=2-2, to=2-3]
\end{tikzcd}\]

We check the axioms of a lax transformation. The respect of $\alpha_m$ for laxators arises from 
the relations of type 2 in $\kappa X,$ while that for unitors arises from relations of type 4.

Approximate naturality of $\alpha_m$ says that, given a cell $\xinlinecell{d_0}{d_1}{d_2}{d_3}{f}{g}{m}{n}{\gamma}$ 
in $\dbl{D}$ as in the definition of lax transformation and an element $h:x_0\proto x_1$ in $Xm,$ 
then the $x_0\xto{[f]_{x_0}} Xf(x_0)\xto{[X\gamma(h)]} Xg(x_1)$ 
must coincide with the composite 
$x_0\xto{[h]} x_1\xto{[g]_{x_1}} Xg(x_1),$ which is a relation of type 3 in $\kappa X.$ 
Approximate functoriality of $\alpha_f$ says that, given $x\in Xd_0$ and $d_0\xto{f} d_1\xto{g} d_2,$ 
$[x]\xto{[f]_x} [Xf(x)]\xto{[g]_{Xf(x)}} [X(f\cdot g)(x)]$ must coincide with 
$[x]\xto{[g\cdot f]_x} [X(g\cdot f)(x)],$ while furthermore $\id_{[x]}$ must coincide with 
$[\id_d]_x,$ which are both relations of type 1. Thus we have a lax transformation $\eta$ as claimed. 

Since lax transformations and modulations between categories, seen as lax functors out of $\dbl{1},$ 
are readily seen to correspond to functors and natural transformations, we get a functor 
$\Cat(\kappa X,C)\to \Lax(X,\Delta C)$ as desired. Along the same line as the proof that $\eta$ is a lax
transformation, one checks that this functor is a bijection on objects. As for morphisms, 
given lax transformations $\alpha,\beta:\kappa X\to C$ (corresponding to functors) and a 
modulation $\mu:\alpha\to \beta,$ the 
induced modulation $\eta *\mu :\eta * \alpha\to \eta * \beta$ has as components the summands of the unique component 
$\mu_{\bullet}:\kappa X_0\to C_0,$ recalling $\kappa X_0=\sum_d Xd,$ which establishes faithfulness. 
As for fullness, given a modulation $\mu:\eta *\alpha \To \eta *\beta: X\to \Delta C,$ 
we get an unnatural transformation $\alpha\to \beta$ whose component at $[x],$ for $x\in Xd,$ is nothing but 
$\mu_d(x),$ while the respect of $\mu$ for loose composition ensures this unnatural transformation 
is in fact natural at morphisms in $\kappa X$ of form $[h],$ and approximate naturality, 
at morphisms in $\kappa X$ of form $[f]_x.$ Thus we have an equivalence of categories as claimed.
\end{proof}

\section{Discrete opfibrations of models of double theories}\label{sec:discrete_opfibrations}

We move toward our main theorem, showing that instances of models of double
theories can be equivalently described as discrete opfibrations over these
models.

A discrete opfibration of ordinary categories can be cast in abstract terms as a functor $p:E\to B$ such that the following square of sets is a pullback:
\[\begin{tikzcd}
	{E^\WalkingArrow} & {\Ob(E)} \\
	{B^\WalkingArrow} & {\Ob(B)}
	\arrow["\dom", from=1-1, to=1-2]
	\arrow["{p^\WalkingArrow}"', from=1-1, to=2-1]
	\arrow["\lrcorner"{anchor=center, pos=0.125}, draw=none, from=1-1, to=2-2]
	\arrow["{\Ob(p)}", from=1-2, to=2-2]
	\arrow["\dom"', from=2-1, to=2-2]
\end{tikzcd}.\]
Now if $E$ and $B$ are viewed as models of the terminal double
theory $\WalkingOb=\{\ob\xproto{\mor}\ob\}$ in $\Span,$
then we can equally well consider the possibility that the following
square be a pullback, which specializes to the case above:
\[\begin{tikzcd}
	{\top(E(\mor))} & {E(\ob)} \\
	{\top(B(\mor))} & {B(\ob)}
	\arrow["s", from=1-1, to=1-2]
	\arrow["{\top(p_{\mor})}"', from=1-1, to=2-1]
	\arrow["\lrcorner"{anchor=center, pos=0.125}, draw=none, from=1-1, to=2-2]
	\arrow["{p_{\ob}}", from=1-2, to=2-2]
	\arrow["s"', from=2-1, to=2-2]
\end{tikzcd}.\]
Here, $\top$ denotes the tabulator, recalled in \cref{def:tabulator}.

We should like to generalize this further to the case of models of 
an arbitrary double theory. The main issue to consider is
what should be done with a nontrivial loose morphism $m:x\proto y$
in the theory $\dbl{D}$. For instance, the base model
$B$ might be the schema of a weighted graph:
\[Bx=\{s,t:E\rightrightarrows V\},\qquad By=\{W\},\qquad Bm=\{w:E\proto W\}.\]
A discrete opfibration $p:F\to B$ over $B$ 
must presumably contain discrete opfibrations over $B_x$ and $B_y$,
that is, an actual graph $G \coloneqq F_x$ and a set $R \coloneqq B_y.$
Over the loose morphism $m,$ the most natural
guess is that for every $e$ over $E,$ we want a unique heteromorphism out of 
$e$ over $w,$ which we'll denote $\bar w:e\proto r(e).$

In other words, we want $p$ to be a discrete opfibration on the 
collages (or barrels) of $B$ and $F$; or, taking the 2-category of profunctors 
to be the strict slice 2-category $\twocat{K} \coloneqq (\Cat/{\WalkingArrow})$ of categories over the
walking arrow, we are taking the discrete opfibrations in $\twocat{K}$ to be those mapped
to a discrete opfibration by the forgetful 2-functor to $\Cat$. 

This motivates the following definition of discrete opfibration, applicable to
models of an arbitrary double theory $\dbl{D}$:

\begin{definition}\label{def:discrete-opfibration}
  A morphism $p:E\to B$ of models of a double theory $\dbl{D}$,
  as in \cref{def:morphism-lax-functors}, is a
  \define{discrete opfibration} if for each loose morphism $m:x\proto y$ in $\dbl{D}$, the square
\[\begin{tikzcd}
	{\top(Em)} & {Ex} \\
	{\top(Bm)} & {Bx}
	\arrow["s", from=1-1, to=1-2]
	\arrow["{\top(p_m)}"', from=1-1, to=2-1]
	\arrow["\lrcorner"{anchor=center, pos=0.125}, draw=none, from=1-1, to=2-2]
	\arrow["{p_x}", from=1-2, to=2-2]
	\arrow["s"', from=2-1, to=2-2]
\end{tikzcd}\]
is a pullback. 
\end{definition}

Note that the case $m=\id_x$ encompasses the condition that each $p_x$ be a discrete opfibration of categories 
in the usual sense, and that this definition applies to both simple and cartesian double theories,
and for any choice of semantics double category.

For models valued in $\Span,$
the definition of discrete opfibration at $m:x\proto y$ says that, for each $h:a\proto b$ in $Bm$ and each $\bar a$ over $a$ in $Ex,$ there is
a unique $\bar b\in Ey$ and $\bar h:\bar a\proto \bar b$ in $Em$ over $h.$ As promised, this is precisely the condition that 
the collage of $p_m$ be a discrete opfibration of ordinary categories.

One might thus wonder whether a discrete opfibration $p:E\to B$ of models can be fully characterized as a 
discrete opfibration of certain categories associated to $E$ and $B,$ by ``flattening out'' the 
profunctors $Em$ and $Bm$ into their collages and gluing appropriately. Presumably, this flattening 
would have to be via the collage construction above (\cref{cons:collage}). Indeed:

\begin{proposition}[Collage creates discrete opfibrations]
A morphism $p:E\to B$ of span-valued models of a simple double theory $\dbl{D}$ is a discrete opfibration if and only if
$\kappa p:\kappa E\to \kappa B$ is a discrete opfibration of categories.
\end{proposition}
\begin{proof}
Recall that $\kappa B$ has morphisms of two kinds: if $f:x\to y$ is a tight morphism in $\dbl{D}$ and $b\in Bx,$ then
we have a morphism $[f]_b:[b]\to [Bf(b)]$ in $\kappa B;$ and if $m:x\proto x'$ is a loose morphism in 
$\dbl{D},$ $b\in Bx,$ $b'\in By,$ and $h:b\proto b'$ in $Bm(x,x'),$ then we have a morphism $[h]:[b]\to [b']$ in $\kappa B.$

If $\kappa(p)$ is a discrete opfibration, then certainly so is $p,$ for (saying essentially the same thing once again)
the discrete opfibration condition on $p$ at $m$ says precisely that the collage of $p_m$ is a discrete opfibration,
and we have seen that this collage maps into $\kappa(p).$ Conversely, we have only to observe that 
$\kappa(p)$ is \emph{always} a discrete opfibration at morphisms of the $f$-type, because we can 
restrict $\kappa(p)$ to $\int p_0,$ the map of discrete opfibrations obtained by taking the 
categories of elements of $B_0$ and $E_0.$ And maps of discrete opfibrations are, themselves, 
automatically discrete opfibrations.
\end{proof}

We can thus easily compute the discrete opfibrations between models of discrete double 
theories, combining the above result with \cref{prop:functors-out-of-discrete}: 
\begin{corollary}[Discrete opfibrations over discrete theories]\label{cor:discrete-opfibrations-discrete-theories}
For a category $C,$ the discrete opfibrations between models of the loosely discrete 
double theory $\dbl{L}C$ coincide with the discrete opfibrations of categories over $C$ 
under the equivalence $\Lax(\dbl{L}C,\Span)\simeq \Cat/C.$ In a similar way, the tightly discrete 
double theory $\dbl{T}C,$ the 
discrete opfibrations between models coincide with maps of functors $C\to \Cat$ which 
become discrete opfibrations on taking the Grothendieck construction.
\end{corollary}

\subsection{Grothendieck construction for instances}

As our main theorem, we show that instances and discrete opfibrations are
equivalent notions, generalizing the well-known equivalence between actions of a
category $C$ and discrete opfibrations over $C.$

\begin{definition}\label{def:cat-of-dopfs}
Let $\Dopf(B)$ denote the full subcategory of $\Lax(\dbl{D},\Span)/B$ spanned by the discrete
opfibrations over $B.$ 
\end{definition}

We remark that since the discrete opfibrations form the right class of a 
factorization system (see \cref{prop:comprehensive-factorization-simple} below),
they are cancellable on the right, so that every morphism between
discrete opfibrations over $B$ is itself a discrete opfibration, just as for
ordinary discrete opfibrations of categories.

\begin{theorem}\label{thm:instances-discrete-opfibrations}
Let $B$ be a span-valued model of a double theory $\dbl{D}$.
There is an equivalence $\nabla:\Dopf(B)\leftrightarrows \Inst(B):\int$
between the category of discrete opfibrations over $B$
and the category of instances of $B,$ which restricts to an equivalence 
$\CartDopf(B)\leftrightarrows \CartInst(B)$ in the case that $\dbl{D}$ and $B$ are cartesian.
\end{theorem}

For an instance $H$, we call $\int H$ the \define{model of elements} or just the
\define{elements} of $H.$

\begin{proof}[Proof (sketch)]
Here we give only the bare structure of the correspondence on objects, deferring all details 
to be checked in Appendix \ref{app:proof-of-discrete-opfibrations}.

Consider a discrete opfibration $p:E\to B.$ We shall define an instance $\nabla p:1\proto B$
giving the indexed view of $p.$ The material of $\nabla p$ is as follows:
\begin{itemize}
\item On objects, define $\nabla p_x \coloneqq (1\leftarrow Ex\xto{p_x} Bx).$
\item Given the tight morphism $f:x\to y$ in $\dbl{D}$, we define
\begin{equation*}
  \begin{tikzcd}
	{\nabla p_x} \\
	{\nabla p_y}
	\arrow["{\nabla p_f}"', from=1-1, to=2-1]
  \end{tikzcd}
  \quad\coloneqq\quad
  \begin{tikzcd}
	1 & Ex & Bx \\
	1 & Ey & By
	\arrow[equals, from=1-1, to=2-1]
	\arrow[from=1-2, to=1-1]
	\arrow["{p_x}", from=1-2, to=1-3]
	\arrow["Ef"', from=1-2, to=2-2]
	\arrow["Bf", from=1-3, to=2-3]
	\arrow[from=2-2, to=2-1]
	\arrow["{p_y}"', from=2-2, to=2-3]
  \end{tikzcd}.
\end{equation*}
\item Given the loose morphism $m:x\proto y$ in $\dbl{D}$, we define the cell $\nabla p_m$ to be:
\[\begin{tikzcd}
	1 & Ex & Bx & Bm & By \\
	&& Em \\
	1 && Ey && By
	\arrow[equals, from=1-1, to=3-1]
	\arrow[from=1-2, to=1-1]
	\arrow["p_x", from=1-2, to=1-3]
	\arrow[from=1-4, to=1-3]
	\arrow[from=1-4, to=1-5]
	\arrow[equals, from=1-5, to=3-5]
	\arrow[from=2-3, to=1-2]
	\arrow["\lrcorner"{anchor=center, pos=0.125, rotate=135}, draw=none, from=2-3, to=1-3]
	\arrow["{p_m}"{description}, from=2-3, to=1-4]
	\arrow[from=2-3, to=3-3]
	\arrow[from=3-3, to=3-1]
	\arrow[from=3-3, to=3-5]
\end{tikzcd}\]
\end{itemize}

In the other direction, consider an instance $P:1\proto B.$ We shall define
the corresponding discrete opfibration $\pi:\int P\to B.$ 

The material of $\int P$ and $\pi$ is as follows:
$P$ has the component $P_0:\dbl{D}_0\to \Span_1$ satisfying 
$P_0\cdot t=B_0.$ To get the functor $(\int P)_0:\dbl{D}_0\to \Set$ and the 
natural transformation $\int P_0\to B_0,$  
we may thus whisker $P_0$ with the natural transformation $r:c\To t:\Span_1\to \Set,$ 
where $r$ picks out the right-hand morphism of a span, $c$ the summit, and $t$ the 
codomain of $r.$ Note that there is no substantial ``collection of elements''
process necessary here, 
since an instance is already defined in a fibered manner, rather than 
in terms of actual functors from the values of $B$ to $\Set.$

\begin{itemize}
\item (On loose morphisms) For $m:x\proto y$ in $\dbl{D},$ we are given the cell $Pm$ below,
where again we define $\int Pm$ and $\pi_m$ by naming appropriate components of $P.$
\[\begin{tikzcd}
	1 & {\int P x} & Bx & Bm & By \\
	&& {\int Pm} \\
	1 && {\int P_y} && By
	\arrow[equals, from=1-1, to=3-1]
	\arrow[from=1-2, to=1-1]
	\arrow["{\pi_x}", from=1-2, to=1-3]
	\arrow[from=1-4, to=1-3]
	\arrow[from=1-4, to=1-5]
	\arrow[equals, from=1-5, to=3-5]
	\arrow[from=2-3, to=1-2]
	\arrow["\lrcorner"{anchor=center, pos=0.125, rotate=135}, draw=none, from=2-3, to=1-3]
	\arrow["{\pi_m}"{description}, from=2-3, to=1-4]
	\arrow[from=2-3, to=3-3]
	\arrow[from=3-3, to=3-1]
	\arrow["{\pi_y}", from=3-3, to=3-5]
\end{tikzcd}\]

\item (On cells) Given a cell $\xinlinecell{x}{y}{z}{w}{f}{g}{m}{n}{\alpha}$ in $\dbl{D},$
we have available so far the data below.
\[\begin{tikzcd}
	{\int P x} & {\int P m} & {\int P y} \\
	{\int P z} & {\int P n} & {\int P w} \\
	Bx & Bm & By \\
	Bz & Bn & Bw
	\arrow["{\int P f}", from=1-1, to=2-1]
	\arrow["{\pi_x}"{description}, curve={height=24pt}, dotted, from=1-1, to=3-1]
	\arrow["{s_m}"', from=1-2, to=1-1]
	\arrow["{t_m}", from=1-2, to=1-3]
	\arrow["{\int P\alpha}", dashed, from=1-2, to=2-2]
	\arrow["{\pi_m}"{description, pos=0.2}, curve={height=12pt}, dotted, from=1-2, to=3-2]
	\arrow["{\int P g}", from=1-3, to=2-3]
	\arrow["{\pi_y}"{description}, curve={height=-24pt}, dotted, from=1-3, to=3-3]
	\arrow["{\pi_z}"{description}, curve={height=24pt}, dotted, from=2-1, to=4-1]
	\arrow["{s_n}", from=2-2, to=2-1]
	\arrow["{t_n}"', from=2-2, to=2-3]
	\arrow["{\pi_n}"{description, pos=0.3}, curve={height=-24pt}, dotted, from=2-2, to=4-2]
	\arrow["{\pi_w}"{description}, curve={height=-24pt}, dotted, from=2-3, to=4-3]
	\arrow["Bf", from=3-1, to=4-1]
	\arrow[from=3-2, to=3-1]
	\arrow[from=3-2, to=3-3]
	\arrow["{B\alpha}", from=3-2, to=4-2]
	\arrow["Bg", from=3-3, to=4-3]
	\arrow[from=4-2, to=4-1]
	\arrow[from=4-2, to=4-3]
\end{tikzcd}\]

Now recall that we defined $\int P n \coloneqq \int Pz \times_{Bz} Bn$, so that
by the universal property of the pullback we may define
$\int P\alpha$ to be the unique map such that $(\int P\alpha)\cdot \pi_n=\pi_m\cdot B\alpha$
and $(\int P\alpha)\cdot s_n=s_m\cdot (\int P f).$ 

\item (Laxators)
Given composable loose morphisms $x\xproto{m} y\xproto{n} z$ in $\dbl{D},$ we have already defined
everything except $q$ in the diagram below. This map $q$ is defined
by the universal property of the pullback at $\bullet_B$ in such a
way as to produce a span map between the spans with apexes $\bullet_E$ and $\bullet_B.$

\[\begin{tikzcd}
	{\int P x} & {\int P m} & {\int P y} & {\int P n} & {\int P z} \\
	&& {\bullet_E} \\
	{\int P x} && {\int P (m\odot n)} && {\int P z} \\
	Bx & Bm & By & Bn & Bz \\
	&& {\bullet_B} \\
	Bx && {B(m\odot n)} && Bz
	\arrow[equals, from=1-1, to=3-1]
	\arrow[from=1-2, to=1-1]
	\arrow[from=1-2, to=1-3]
	\arrow["\pi"{description}, dotted, from=1-2, to=4-2]
	\arrow[from=1-4, to=1-3]
	\arrow[from=1-4, to=1-5]
	\arrow["\pi"{description}, dotted, from=1-4, to=4-4]
	\arrow[equals, from=1-5, to=3-5]
	\arrow[from=2-3, to=1-2]
	\arrow["\lrcorner"{anchor=center, pos=0.125, rotate=135}, draw=none, from=2-3, to=1-3]
	\arrow[from=2-3, to=1-4]
	\arrow["{\int P_{m,n}}"', dashed, from=2-3, to=3-3]
	\arrow["q"{description}, curve={height=-24pt}, dotted, from=2-3, to=5-3]
	\arrow["{s_{m\odot n}}"{description, pos=0.7}, from=3-3, to=3-1]
	\arrow[from=3-3, to=3-5]
	\arrow["\pi"{description}, curve={height=12pt}, dotted, from=3-3, to=6-3]
	\arrow[equals, from=4-1, to=6-1]
	\arrow[from=4-2, to=4-1]
	\arrow[from=4-2, to=4-3]
	\arrow[from=4-4, to=4-3]
	\arrow[from=4-4, to=4-5]
	\arrow[equals, from=4-5, to=6-5]
	\arrow[from=5-3, to=4-2]
	\arrow["\lrcorner"{anchor=center, pos=0.125, rotate=135}, draw=none, from=5-3, to=4-3]
	\arrow[from=5-3, to=4-4]
	\arrow["{B_{m,n}}", from=5-3, to=6-3]
	\arrow[from=6-3, to=6-1]
	\arrow[from=6-3, to=6-5]
\end{tikzcd}\]

This permits us to define the laxator
\begin{equation*} \textstyle
  \int P_{m,n}:\bullet_E \to \int P(m\odot n)=B(m\odot n)\times_{Bx} \int P x
\end{equation*}
using the universal property of the pullback. Specifically, we require that
$\int P_{m,n}\cdot \pi_{m\odot n}=q\cdot B_{m,n}$
and that $\int P_{m,n}\cdot s_{m\odot n}$ makes the left-hand pentagon commute.

\item (Unitors)
Similarly, we are given the data to uniquely define a unitor
from the universal property of the pullback $\int P \id_x=\int P x \times_{Bx} B\id_x$
by filling the upper boundary below:
\[\begin{tikzcd}
	{\int Px} & {\int Px} & {\int Px} \\
	{\int Px} & {\int P \id_x} & {\int Px} \\
	Bx & Bx & Bx \\
	Bx & {B\id_x} & Bx
	\arrow[equals, from=1-1, to=2-1]
	\arrow[curve={height=18pt}, dotted, from=1-1, to=3-1]
	\arrow[from=1-2, to=1-1]
	\arrow[from=1-2, to=1-3]
	\arrow["{P_x}"', dashed, from=1-2, to=2-2]
	\arrow[curve={height=-18pt}, dotted, from=1-2, to=3-2]
	\arrow[equals, from=1-3, to=2-3]
	\arrow[curve={height=18pt}, dotted, from=2-1, to=4-1]
	\arrow[from=2-2, to=2-1]
	\arrow[from=2-2, to=2-3]
	\arrow[curve={height=-18pt}, dotted, from=2-2, to=4-2]
	\arrow[equals, from=3-1, to=4-1]
	\arrow[from=3-2, to=3-1]
	\arrow[from=3-2, to=3-3]
	\arrow["{B_x}"{description}, from=3-2, to=4-2]
	\arrow[equals, from=3-3, to=4-3]
	\arrow[from=4-2, to=4-1]
	\arrow[from=4-2, to=4-3]
\end{tikzcd}\]

Thus $\int P_x\cdot \pi_{\id_x}=\pi_x\cdot B_x,$ 
while $\int P_x \cdot s_x=\id_{\int Px}.$ 
\end{itemize}

It remains to check that $\nabla p$ is an instance; that $\int P$ is a model and $\pi$
is a discrete opfibration into $B$; and to show that $\nabla$ and $\int$ extend to mutually
quasi-inverse functors. All this is routine, if lengthy, verification, and we defer the details
to Appendix \ref{app:proof-of-discrete-opfibrations}. Here, we finish with the 
additional argument required for the cartesian case: 

Suppose that $P$ is a cartesian instance of $B.$ We must show that
$\int P:\dbl{D}\to \Span$ preserves finite products on
objects and on loose morphisms. On objects, this is precisely
the definition of a cartesian instance. As for loose morphisms,
the pullback square characterizing $\int P$ as a discrete opfibration
makes $\int P(\id_I)= B(\id_I)\times_{B(I)} P(I)=I\times_I I=I,$
handling the nullary case. For the binary case, this follows
from commutation of limits with limits since
$\int P(m\times m')$ is the pullback of the span
$B(m\times m')\times_{B(x\times x')} \int P(x\times x'),$ which is
isomorphic to the product of the spans for $m$ and $m'.$

Conversely, if we suppose that $\int P$ is a cartesian model, then by definition
we have the canonical isomorphisms
$\int P(x\times y)\cong \int P(x)\times \int P(y)$ and
$\int P(I)\cong I,$ which give rise to the cartesian instance isomorphisms for
$\nabla \int P$ and thus for the isomorphic instance $P.$
\end{proof}

\subsection{Comprehensive factorization systems for models of double theories}

One of the most important properties of discrete opfibrations of categories is that they 
are precisely the functors right orthogonal to the initial functors \cite{street-walters-1973}.
We cannot recover such a result for models of double theories directly from the proposition above,
because we know that the collage functor $\kappa$
is neither full nor faithful (\cref{rem:kappa-bad}). Nonetheless, we can still establish the existence of such 
a factorization system: 

\begin{proposition}[Comprehensive factorization for models of simple theories]\label{prop:comprehensive-factorization-simple}
  The discrete opfibrations of span-valued models of a simple double theory
  $\dbl{D}$ are the right class of an orthogonal factorization system on
  $\Lax(\dbl{D},\Span).$
\end{proposition}
\begin{proof}
  Recall from \cref{prop:models_locally_presentable} that $\Lax(\dbl{D},\Span)$ is locally presentable,
  with limits and filtered colimits created by the forgetful functor
  \begin{equation*}
    U:\Lax(\dbl{D},\Span)\to \Set^{\Loose(\dbl{D})+\Ob(\dbl{D})},
    \quad
    X\mapsto(m\mapsto X(m);d\mapsto X(d)).
  \end{equation*}
  Thus, by Theorem 11.3 of Kelly's \cite{transfinite-constructions-1980}, for any small set $\mathcal{G}$
  of maps in $\Lax(\dbl{D},\Span)$, there is a factorization system whose right class 
  $\mathcal{G}^\perp$ is the class of maps right orthogonal to $\mathcal{G}$.

  Now consider the forgetful functor
  \begin{equation*}
    V:\Lax(\dbl{D},\Span)\to \Span_1^{\Loose(\dbl{D})}
    \quad
    X\mapsto (m:x\proto y)\mapsto (Xx\leftarrow Xm\to Xy),
  \end{equation*}
  where $\Span_1$ is the loose category of $\Span$,
  hence its \emph{objects} are spans of sets and we consider the mere \emph{set} $\Loose(\dbl{D})$
  of loose arrows. Since $V$ factors through $U$ via the limit- and filtered-colimit-creating
  functor which evaluates every span on its summit, and also the span at $\id_x$ at its domain,
  $V$ preserves limits and filtered colimits. Thus, $V$ has a left adjoint $F.$
  We will define the desired factorization system to be generated by $F(\mathcal{G})$ for a set 
  $\mathcal{G}$ generating a natural notion of ``discrete opfibration'' in $\Span_1^{\Loose(\dbl{D})}.$
  
  To wit, given spans $E_0\leftarrow E\to E_1$ and $B_0\leftarrow B\to B_1$ 
  in $\Set,$ we say that a span map $(p_0,p,p_1): (E_0\leftarrow E\to E_1) \to (B_0\leftarrow B\to B_1)$
  is a discrete opfibration if and only if the square
  \[\begin{tikzcd}
    {E} & {E_0} \\
    {B} & {B_0}
    \arrow[to=1-1, from=1-2]
    \arrow[to=2-1, from=2-2]
    \arrow["p_0", from=1-2, to=2-2]
    \arrow["p"', from=1-1, to=2-1]
  \end{tikzcd}\]
  is a pullback; in other words, the discrete opfibrations in $\Span_1$ are the class right orthogonal 
  to the span map $f:(1\leftarrow 0\to 0)\to (1\leftarrow 1\to 1).$ Then we can define the discrete opfibrations
  in a power $\Span_1^A$ levelwise, and this class will consist of the morphisms right orthogonal
  to the class $\mathcal{G}=\{f^a\}$ of maps $f^a$ in $\Span_1^A$ given by $f$ at some particular $a\in A$ and the identity of the 
  empty span elsewhere.
  Thus lifting along $f^a$ adjoins a lift of a heteromorphism at the $a$ component, given a lift of its domain.
  
  One sees immediately that the discrete opfibrations in $\Lax(\dbl{D},\Span)$ are those mapped to discrete
  opfibrations in $\Span_1^{\Loose(\dbl{D})}$ under $V$, which means that they are precisely the right orthogonality 
  class $F(\mathcal{G})^\perp,$ as desired.
\end{proof}

As in the simple case, we have
\begin{proposition}[Cartesian models are locally presentable]\label{prop:cartesian-models-locally-presentable}
If $\dbl{D}$ is a cartesian double theory, then the category 
$\Cart(\dbl{D},\Span)$, consisting of (cartesian) span-valued models of $\dbl{D}$
and tight natural transformations, is locally finitely presentable. The full inclusion 
$\Cart(\dbl{D},\Span)\to \Lax(\dbl{D},\Span)$ is a finitely accessible right adjoint.
\end{proposition}
\begin{proof}
Recall from \cref{prop:models_locally_presentable} that the simple models
$\Lax(\dbl{D},\Span)$ are locally presentable, being models of the finite limit 
sketch $\varphi(\dbl{D})$ constructed there. Now, a cartesian model is precisely 
one sending products of objects in $\dbl{D}$ to products of sets, and similarly for loose morphisms.
For any two objects $d_1,d_2\in \dbl{D},$ we have a cone $[d_1]\leftarrow [d_1\times d_2]\to [d_2]$
in $\varphi(\dbl{D}),$ and similarly for pairs of loose morphisms, as well as for 
the nullary cones on the terminal object and loose morphism. Thus, cartesian models are precisely
the models of the finite limit sketch $\varphi(X)$ extended with the product cones just
listed.
\end{proof}

We can also construct a comprehensive factorization for
morphisms of cartesian models much as in the simple case.

\begin{corollary}[Comprehensive factorization of cartesian models]
There is an orthogonal factorization system on the category of cartesian models
of a cartesian double theory $\dbl{D}$ whose right class consists of those 
morphisms $p:X\to Y$ of cartesian models which are discrete opfibrations
when considered as morphisms of simple models.  
\end{corollary}
\begin{proof}
Consider a small class $\mathcal{G}$ generating the 
initial morphisms of $\Lax(\dbl{D},\Span).$ 
Since by \cref{prop:cartesian-models-locally-presentable} we have 
a left adjoint $L\dashv i,$ to say that given a morphism $p:X\to Y$ of cartesian models,
$ip$ is a discrete opfibration, i.e. $ip$ is right orthogonal to $\mathcal{G},$ 
it is equivalent to say that $p$ is right orthogonal to $L(\mathcal{G}).$ 
This being a small class in a locally presentable category, as above we 
may obtain a factorization system as desired.
\end{proof}

We call the left class of this factorization system the \emph{initial morphisms} of models of $\dbl{D}$ in $\Span,$ by
analogy with the case of initial functors between categories, i.e., models of the terminal double theory. 
\begin{example}[Initial morphisms of profunctors]
Let us now consider the case of profunctors, that is, models of the walking loose morphism, $\WalkingLoose=(\vdash\xproto{\ell}\dashv).$
Consider a profunctor as a barrel, that is, a category $C\to \WalkingArrow$ over the walking arrow
$0\xto{a} 1,$ and recall (see e.g. Joyal's notes \cite[Theorem 5.2]{joyal-factorisation}) that any initial functor may be written as a transfinite composition 
of pushouts of this form:
\[\begin{tikzcd}
	{\{0\}} & C \\
	\WalkingArrow & C'
	\arrow[from=1-1, to=1-2]
	\arrow[from=1-1, to=2-1]
	\arrow["{C'}", from=1-2, to=2-2]
	\arrow[from=2-1, to=2-2]
\end{tikzcd}\]
To lay out such a sequence of pushouts over $\WalkingArrow$ is precisely to choose, for each pushout
such that the image of $\{0\}$ is in $C_\vdash,$ whether the new object in $C'$ is to lie over
$0$ or $1$ in $\WalkingArrow.$ Thus there are three kinds of generating initial morphisms which 
lie over $\WalkingArrow$: those which add a morphism over $1,$ those which add a morphism over $0,$ 
and those which add a heteromorphism over $a.$ Thus the initial functors lying over $\WalkingArrow$ 
are generated, as the left class of a factorization system, by the same three arrows as 
the initial morphisms of models of $\WalkingLoose,$ so they are the same class. 
\end{example}

\begin{proposition}[Initial morphisms over degenerate double theories]
Let $\dbl{D}$ be a simple double theory. 
\begin{enumerate}
\item If $\dbl{D}=\dbl{T}C$ is tightly discrete, 
then under the equivalence $\Lax(\dbl{D},\Span)\simeq \CAT(C,\Cat)$ (\cref{prop:span_models_of_discrete_theories}) 
the initial morphisms of models of $\dbl{D}$ are precisely the levelwise-initial natural transformations.
\item If $\dbl{D}=\dbl{L}C$ is loosely discrete, then under the equivalence 
$\Lax(\dbl{D},\Span)\simeq \Cat/C,$ the initial morphisms of models of $\dbl{D}$ are precisely the initial functors over $C.$
\end{enumerate}
\end{proposition}
\begin{proof}

\begin{enumerate}
\item We have already seen that the discrete opfibrations of models of $\dbl{T}C$ 
are detected levelwise, and the levelwise initial functors are precisely the left orthogonality 
class of the levelwise discrete opfibrations.

\item Similarly, we have already seen that the discrete opfibrations of models correspond 
discrete opfibrations over $C,$ so this reduces to the lifting of factorization systems to slices.
\end{enumerate} 

\end{proof}

\begin{example}[Initial morphisms of natural transformations]
Next let $\dbl{D}$ be the simplest double theory generated by a nontrivial 2-category. That is, 
$\dbl{D}$ is generated by a parallel pair of tight morphisms $c\xto{f,g}d$ and a 2-cell $\xi:f\To g.$
Then models of $\dbl{D}$ are natural transformations $\Xi:F\To G:\cat{C}\to \cat{D}$
between parallel pairs of functors. 

There are two generating initial morphisms over 
$\dbl{D},$ one per object. The generator corresponding to $\id_d$ is uninteresting, just generating 
$\dbl{D}$-model morphisms that are identity on the $c$ component and arbitrary inital functors on 
the $d$ components. The generating initial morphism corresponding to $\id_c$ is more interesting:
it includes the canonical natural transformation between the two functors $s,t:\bullet\to \WalkingArrow$
into its levelwise product with another copy of $\WalkingArrow.$ Thus pushing out along a copy of this morphism
adjoints to a model $X$ of $\dbl{D}$ a new morphism in the $c$ component, with given domain, and a whole new square,
with one edge given by a component of the natural transformation constituting $X,$ to the $d$ component.

Unlike the cases above, we do not see how to describe these initial morphisms of natural transformations
in terms of ordinary initial functors.
\end{example}

As with any orthogonal factorization system, we can
immediately derive the following by factoring a model morphism into an initial map followed by a discrete opfibration.

 \begin{corollary}
  The category $\Dopf(B)$ of discrete opfibrations over a model $B$ of a simple double theory is reflective in 
  the category of models $\Lax(\dbl{D},\Span)/B$ sliced over $B.$
 \end{corollary}

This allows us to \emph{present} a discrete opfibration over $B$ via an arbitrary
morphism $f:X\to B$ of models, which is the implementation approach currently taken in the 
CatColab tool.

\bibliographystyle{alpha}
\bibliography{double-instances}

\begin{appendices}

\section{Loose transformations versus modules}\label{app:loose-transformations}

In this appendix, we prove results relating the two notions of loose morphism between lax functors,
namely loose transformations and modules, as well as their respective higher morphisms, namely
modifications and modulations. The main constructions here are given, and the main results are
asserted, for bicategories in \cite{modules-2003}. We take the opportunity to make a correction
(\cref{rem:modifications-not-faithful-in-modulations}) on the relationship between
modifications and modulations and to prove the result (\cref{thm:loose-transformations-into-modules})
that the category of loose transformations and modulations is a full subcategory of that of
modules and modulations in the more general setting of double categories.

We first recall the definition of a loose transformation, essentially as
given in \cite{transposing-2024}.\footnote{We reserve the name ``loose transformation'' 
for the lax case of an object that also appears in colax and pseudo variants, 
because that is the only version of interest here.} Such transformations generalize the
transformations of internal category theory to pseudocategories internal to any
2-category; but let us caution that they are \emph{not} the transformations between
pseudofunctors of pseudocategories studied in the subject's most detailed
reference \cite{pseudo-categories-2006}, which instead generalize the tight
transformations used above. Furthermore, there can be no ``packed'' definition of a loose
transformation in terms of multivariate lax functors, since lax functors
out of a product with the walking loose arrow already serve to define modules,
which are more general. Thus we must be content with an unpacked definition:

\begin{definition}[Loose transformation] \label{def:lax-loose-transformation}
Let $F,G: \dbl{D} \to \dbl{E}$ be lax functors between double categories. 
A \define{loose transformation} $\tau: F \to G$ consists of:
\begin{itemize}
  \item A functor $\tau_0:\dbl{D}_0\to \dbl{E}_1,$ whose values $\tau_x \coloneqq \tau_0(x)$
    and $\tau_f \coloneqq \tau_0(f)$ at objects and tight morphisms of $\dbl{D}$
    are called the \define{components} of $\tau$ and look like this:
      \[\begin{tikzcd}
        Fx & Gx \\
        Fy & Gy
        \arrow["Ff"', from=1-1, to=2-1]
        \arrow[""{name=0, anchor=center, inner sep=0}, "{\tau_x}", "\shortmid"{marking}, from=1-1, to=1-2]
        \arrow["Gf", from=1-2, to=2-2]
        \arrow[""{name=1, anchor=center, inner sep=0}, "{\tau_y}"', "\shortmid"{marking}, from=2-1, to=2-2]
        \arrow["{\tau_f}"{description}, draw=none, from=0, to=1]
      \end{tikzcd}\]
  \item A globular natural transformation, also denoted $\tau,$ with signature
    \[\tau:\tau_0\odot G_1 \To F_1\odot \tau_0 :\dbl{D}_1 \to \dbl{E}_1,\] whose component at a
    loose morphism $m:x\proto y,$ called a \define{naturality comparison}, looks
    like this:
      \[\begin{tikzcd}
        Fx & Gx & Gy \\
        Fx & Fy & Gy
        \arrow["Fm"', "\shortmid"{marking}, from=2-1, to=2-2]
        \arrow["{\tau_y}"', "\shortmid"{marking}, from=2-2, to=2-3]
        \arrow["{\tau_x}", "\shortmid"{marking}, from=1-1, to=1-2]
        \arrow["Gm", "\shortmid"{marking}, from=1-2, to=1-3]
        \arrow[Rightarrow, no head, from=1-1, to=2-1]
        \arrow[Rightarrow, no head, from=1-3, to=2-3]
        \arrow["{\tau_m}"{description}, draw=none, from=1-2, to=2-2]
      \end{tikzcd}\]
\end{itemize}
This data must be coherent with the unitors and laxators of $F$ and $G$ in the
sense that the following diagrams commute, in the functor categories indicated
below each diagram:
\[\begin{tikzcd}[column sep=small]
  {\tau_0\odot G_1\odot G_1} && {F_1\odot \tau_0\odot G} && {\tau_0\odot\mathrm{id}_G} & {\tau_0\odot G(\mathrm{id})} \\
  {\tau_0\odot G(-\odot -)} && {F_1\odot F_1\odot \tau_0} && {\mathrm{id}_F\odot \tau_0} & {F(\mathrm{id})\odot \tau_0} \\
  & {F(-\odot -)\odot \tau} \\
  {\dbl{D}_1\times_{\dbl{D}_0}\dbl{D}_1} && {\dbl{E}_1} && {\dbl{D}_0} & {\dbl{E}_1}
  \arrow["{\tau\odot 1_G}", from=1-1, to=1-3]
  \arrow["{1_\tau\odot G_{-,-}}"', from=1-1, to=2-1]
  \arrow["{1_F\odot \tau}", from=1-3, to=2-3]
  \arrow["{1_\tau\odot G_{-}}", from=1-5, to=1-6]
  \arrow[equals, from=1-5, to=2-5]
  \arrow["{\tau_{\mathrm{id}}}", from=1-6, to=2-6]
  \arrow["{\tau_{-\odot -}}"{description}, from=2-1, to=3-2]
  \arrow["{F_{-,-}\odot 1_{\tau}}"{description}, from=2-3, to=3-2]
  \arrow["{F_{-}\odot 1_\tau}"', from=2-5, to=2-6]
  \arrow[from=4-1, to=4-3]
  \arrow[from=4-5, to=4-6]
\end{tikzcd}\]
\end{definition}

There are two kinds of morphisms between loose transformations: modifications and modulations. 

\begin{definition}[Modifications and modulations of loose transformations]\label{def:modmodofpro}
Let $\sigma,\tau: F\To G:\dbl{D}\to \dbl{E}$ be loose transformations between lax double functors. 
A \define{modification} $\sigma \to \tau$ is a globular natural transformation
$\mu:\sigma_0\Rightarrow \tau_0:\dbl{D}_0\to \dbl{E}_1,$ thus with components as below left, such that 
the square below right commutes in $\Cat(\dbl{D}_1,\dbl{E}_1):$
\[\begin{tikzcd}
	Fx & Gx & {\sigma_0\odot G} & {\tau_0\odot G} \\
	Fx & Gx & {F\odot\sigma_0} & {F\odot \tau_0}
	\arrow[""{name=0, anchor=center, inner sep=0}, "{\sigma_x}"{inner sep=.8ex}, "\shortmid"{marking}, from=1-1, to=1-2]
	\arrow[equals, from=1-1, to=2-1]
	\arrow[equals, from=1-2, to=2-2]
	\arrow["{\mu\odot 1_G}", from=1-3, to=1-4]
	\arrow["\sigma"', from=1-3, to=2-3]
	\arrow["\tau", from=1-4, to=2-4]
	\arrow[""{name=1, anchor=center, inner sep=0}, "{\tau_x}"'{inner sep=.8ex}, "\shortmid"{marking}, from=2-1, to=2-2]
	\arrow["{1_F\odot \mu}"', from=2-3, to=2-4]
	\arrow["{\mu_x}"{description}, draw=none, from=0, to=1]
\end{tikzcd}\]

In contrast, a \define{modulation} $\sigma \to \tau$ between loose transformations is a globular natural
transformation $\mu: \sigma_0\to F(\id)\odot \tau_0$, thus with components as below left, 
such that the rectangle below right commutes, also in $\Cat(\dbl{D}_1,\dbl{E}_1):$
\[\begin{tikzcd}
	& Fx & Gx & {\sigma_0\odot G} & {F(\mathrm{id})\odot \tau_0\odot G} & {F(\mathrm{id})\odot F\odot \tau_0} \\
	Fx & Fx & Gx & {F\odot\sigma_0} & {F\odot F(\mathrm{id})\odot \tau_0} & {F\odot \tau_0}
	\arrow[""{name=0, anchor=center, inner sep=0}, "{\sigma_x}", from=1-2, to=1-3]
	\arrow[equals, from=1-2, to=2-1]
	\arrow[equals, from=1-3, to=2-3]
	\arrow["{\mu\odot 1_G}", from=1-4, to=1-5]
	\arrow["\sigma"', from=1-4, to=2-4]
	\arrow["{1\odot \tau}", from=1-5, to=1-6]
	\arrow["{F_{\mathrm{id},-}\odot 1}", from=1-6, to=2-6]
	\arrow["{F(\mathrm{id}_x)}"'{inner sep=.8ex}, "\shortmid"{marking}, from=2-1, to=2-2]
	\arrow["{\tau_x}"'{inner sep=.8ex}, "\shortmid"{marking}, from=2-2, to=2-3]
	\arrow["{1_F\odot \mu}"', from=2-4, to=2-5]
	\arrow["{F_{-,\mathrm{id}}\odot 1}"', from=2-5, to=2-6]
	\arrow["{\mu_x}"{description}, draw=none, from=0, to=2-2]
\end{tikzcd}\]
\end{definition}

One motivation for the non-obvious definition of modulation is that modifications are sometimes too strict.
For instance, a lax functor $F:\dbl{1}\to \Span$ is a category; a loose transformation 
$\tau:F\to G$ between two such lax functors can be a functor, in the case that it is carried by a span
whose left leg is an identity; but the category of such transformations and modifications
between them is discrete, while the category of such transformations and modulations reproduces
the usual functor category. In the double-categorical world, this example loses some of its sheen,
since we prefer to recover functors via tight transformations, but it is still illustrative. 

A more substantive reason to consider modulations between loose transformations is that they are
in bijection with modulations between the induced modules, as we shall show below. This is 
particularly important in that a modulation of loose transformations has components indexed by 
objects of $\dbl{D},$ while a modulation of modules has components indexed by loose arrows of $\dbl{D},$ 
so the former is a lighter-weight structure, as we saw throughout the body of the paper in the 
special case of instances. 

\begin{remark}[Modifications not faithful in modulations]\label{rem:modifications-not-faithful-in-modulations}
One can induce a modulation from a modification $\mu:\sigma \to \tau$
of loose transformations by pasting the unitor of $F$ to the left of each component cell
of $\mu.$ The modulation axiom follows from the modification axiom together with unitality of $F.$ 

\emph{Contra} Cockett, Koslowski, Seely, and Wood in \cite[4.3]{modules-2003}, where the 
induced 2-functor 
from the bicategory of lax functors, lax transformations, and modifications to that
of lax functors, lax transformations, and modulations is claimed to be an sub-2-category inclusion, 
we note that since the mapping of modifications to modulations involves whiskering with a unitor, 
which may not be monic, this mapping may not be faithful. 

For example, consider lax functors $S,T:\dbl{1}\to \dbl{C}\mathsf{at},$ where for the codomain we view the 2-category
of categories, functors, and natural transformations as a double category in the \emph{loose} direction.
Thus such a lax functor $T$ is the same as an ordinary monad on a category, as has been 
known since lax functors were first defined (\cite[Definition 5.4.1]{benabou-bicategories}) 
If $T$ corresponds to the monad $T:\cat{C}_T\to \cat{C}_T,$ then 
a lax transformation $q:S\to T$ amounts to a functor $F_q:\cat{C}_S\to \cat{C}_T$ 
together with a natural transformation $\alpha_q:F_q S\to TF_q,$ appropriately respecting the 
monad units and multiplications in a sense on which our counterexample will depend only insofar
as we need to let $T$ have its identity endomorphism, with $F_{\id}=\id_{\cat{C}_T}$ and $\alpha_{\id}=\id.$ 
Finally, a modification $\mu:q\To q':S\to T$ is constituted by a natural transformation 
$\mu:F_q\to F_{q'}$ such that a certain square commutes. 

Now consider the constant $0$ monad on the category $\cat{Ab}$ of abelian groups, with its identity monad
endomorphism. The endo-modifications of $\id_0$ coincide with the endo-natural transformations of $\id_{\cat{Ab}},$ 
since the modification constraint is trivial in this case, and one knows that there are thus infinitely 
many distinct such endo-modifications, one per integer. 
However, for any modification $\mu:\id_0\To \id_0,$ the corresponding 
modulation $\bar\mu$ is constituted by a natural transformation $\id_{\cat{Ab}}\to \id_{\cat{Ab}}\cdot 0,$
of which there is, of course, only one. 

The upshot is that, in order that the embedding of modifications in modulations be faithful, it is 
sufficient that, what's quite common in practice, the unitors be monic; but this 
faithfulness does not hold at all in general.
\end{remark}

We aim next to functorially produce a module from any loose transformation. The construction below is in essence
identical to that given by Cockett, Koslowski, Seely, and Wood in \cite{modules-2003} for the 
case of bicategories. Here we cover a more general case by adding tight morphisms, but this does 
not complicate the arguments.

\begin{construction}[Modules and modulations from loose transformations]\label{cons:loose-transformation-to-module}
  We give four constructions: (1) a module from a loose transformation, (2) a
  loose transformation from a module with invertible left actions, (3) a
  modulation between modules from a modulation between loose transformations,
  and (4) a converse construction to (3).
  \begin{enumerate}
  \item \emph{Given a lax loose transformation $\tau: F \To G:\dbl{D}\to \dbl{E}$ between lax double functors,
  we construct a module $M_\tau:F\proTo G.$}
  
    First, we given the induced functor $M_\tau:\dbl{D}_1\to \dbl{E}_1.$ Given 
    $m:x\proto y$ in $\dbl{D},$ we define $M_\tau(m):Fx\proto Gy$ as the composite $Fx\xproto{Fm} Fy \xproto{\tau_y} Gy,$ and 
    for a cell $\xinlinecell{x}{y}{x'}{y'}{f}{g}{m}{m'}{\alpha}$ in $\dbl{D},$ we define $M_\tau(\alpha)$ as the analogous composite:
    \[\begin{tikzcd}
      Fx & Fy & Gy \\
      {Fx'} & {Fy'} & {Gy'}
      \arrow[""{name=0, anchor=center, inner sep=0}, "Fm"{inner sep=.8ex}, "\shortmid"{marking}, from=1-1, to=1-2]
      \arrow["Ff"', from=1-1, to=2-1]
      \arrow[""{name=1, anchor=center, inner sep=0}, "{\tau_y}"{inner sep=.8ex}, "\shortmid"{marking}, from=1-2, to=1-3]
      \arrow["Fg"{description}, from=1-2, to=2-2]
      \arrow["Gg", from=1-3, to=2-3]
      \arrow[""{name=2, anchor=center, inner sep=0}, "{Fm'}"'{inner sep=.8ex}, "\shortmid"{marking}, from=2-1, to=2-2]
      \arrow[""{name=3, anchor=center, inner sep=0}, "{\tau_{y'}}"', from=2-2, to=2-3]
      \arrow["{F\alpha}"{description}, draw=none, from=0, to=2]
      \arrow["{\tau_g}"{description}, draw=none, from=1, to=3]
    \end{tikzcd}\]
    That this provides a functor $\dbl{D}_1\to \dbl{E}_1$ follows from functoriality of $\tau_f$ in $f$ 
    together with that of $F\alpha$ in $\alpha.$

    Next, we give the left and right actions making $M_\tau$ into a module. Here are $M^\ell_{m,n}$ and $M^r_{m,n}$ for
    $m:x\proto y$ and $n:y\proto z$:\footnote{Note that the left action would not have been definable had we attempted to set 
    $M_\tau(m)=\tau_x\cdot Gm$, since, $\tau$ being lax, we have no way to move a $\tau$ leftward. 
    This definition is instead appropriate for the module associated to an \emph{oplax} loose transformation,
    a notion we do not consider here.} 
    \[\begin{tikzcd}
      Fx & Fy & Fz & Gz & Fx & Fy & Gy & Gz \\
      Fx && Fz & Gz & Fx & Fy & Fz & Gz \\
      &&&& Fx && Fz & Gz
      \arrow["Fm"{inner sep=.8ex}, "\shortmid"{marking}, from=1-1, to=1-2]
      \arrow[equals, from=1-1, to=2-1]
      \arrow["Fn"{inner sep=.8ex}, "\shortmid"{marking}, from=1-2, to=1-3]
      \arrow["{M_\tau n}"{inner sep=.8ex}, "\shortmid"{marking}, curve={height=-18pt}, from=1-2, to=1-4]
      \arrow["{\tau_z}"{inner sep=.8ex}, "\shortmid"{marking}, from=1-3, to=1-4]
      \arrow[equals, from=1-3, to=2-3]
      \arrow[equals, from=1-4, to=2-4]
      \arrow["Fm"{inner sep=.8ex}, "\shortmid"{marking}, from=1-5, to=1-6]
      \arrow["{M_\tau m}"{inner sep=.8ex}, "\shortmid"{marking}, curve={height=-18pt}, from=1-5, to=1-7]
      \arrow[equals, from=1-5, to=2-5]
      \arrow["{\tau_y}"{inner sep=.8ex}, "\shortmid"{marking}, from=1-6, to=1-7]
      \arrow[equals, from=1-6, to=2-6]
      \arrow["Gn"{inner sep=.8ex}, "\shortmid"{marking}, from=1-7, to=1-8]
      \arrow["{\tau_n}"{description}, draw=none, from=1-7, to=2-7]
      \arrow[equals, from=1-8, to=2-8]
      \arrow[""{name=0, anchor=center, inner sep=0}, "{F(m\odot n)}"'{inner sep=.8ex}, "\shortmid"{marking}, from=2-1, to=2-3]
      \arrow["{M_\tau (m\odot n)}"'{inner sep=.8ex}, "\shortmid"{marking}, curve={height=24pt}, from=2-1, to=2-4]
      \arrow["{\tau_z}"'{inner sep=.8ex}, "\shortmid"{marking}, from=2-3, to=2-4]
      \arrow["Fm"{inner sep=.8ex}, "\shortmid"{marking}, from=2-5, to=2-6]
      \arrow[equals, from=2-5, to=3-5]
      \arrow["Fn"{inner sep=.8ex}, "\shortmid"{marking}, from=2-6, to=2-7]
      \arrow["{\tau_z}"{inner sep=.8ex}, "\shortmid"{marking}, from=2-7, to=2-8]
      \arrow[equals, from=2-7, to=3-7]
      \arrow[equals, from=2-8, to=3-8]
      \arrow[""{name=1, anchor=center, inner sep=0}, "{F(m\odot n)}"'{inner sep=.8ex}, "\shortmid"{marking}, from=3-5, to=3-7]
      \arrow["{M_\tau(m\odot n)}"{description}, curve={height=24pt}, from=3-5, to=3-8]
      \arrow["{\tau_z}"{inner sep=.8ex}, "\shortmid"{marking}, from=3-7, to=3-8]
      \arrow["{F_{m,n}}"{description}, draw=none, from=1-2, to=0]
      \arrow["{F_{m,n}}"{description}, draw=none, from=2-6, to=1]
    \end{tikzcd}\]
  
  \item \emph{Let $M:F\proto G:\dbl{D}\to\dbl{E}$ be a module such that $M^\ell$ is invertible. We construct a loose transformation $\tau^M:F \to G.$}

   First, let $\tau^M_x \coloneqq M_{\id_x}:Fx\proto Gx,$ and similarly $\tau^M_f \coloneqq M_{\id_f}.$
   This provides a functor $\tau_0$ by functoriality of $M.$
   For the naturality comparison $\tau_m,$ we take the following pasting: 
  \[\begin{tikzcd}
    Fx & Gx & Gy \\
    Fx && Gy \\
    Fx & Fy & Gy
    \arrow["{M\mathrm{id}_x}"{inner sep=.8ex}, "\shortmid"{marking}, from=1-1, to=1-2]
    \arrow[equals, from=1-1, to=2-1]
    \arrow["Gm"{inner sep=.8ex}, "\shortmid"{marking}, from=1-2, to=1-3]
    \arrow[equals, from=1-3, to=2-3]
    \arrow[""{name=0, anchor=center, inner sep=0}, "Mm"'{inner sep=.8ex}, "\shortmid"{marking}, from=2-1, to=2-3]
    \arrow[equals, from=2-1, to=3-1]
    \arrow[equals, from=2-3, to=3-3]
    \arrow["Fm"'{inner sep=.8ex}, "\shortmid"{marking}, from=3-1, to=3-2]
    \arrow["{M\mathrm{id}_y}"'{inner sep=.8ex}, "\shortmid"{marking}, from=3-2, to=3-3]
    \arrow["{M^r_m}"{description}, draw=none, from=1-2, to=0]
    \arrow["{(M^\ell_m)^{-1}}"{description, pos=0.7}, draw=none, from=0, to=3-2]
  \end{tikzcd}\]
  Naturality in $m$ follows from naturality of the actions of $M.$

  \item \emph{Given a modulation $\mu:\sigma\To \tau:F\to G$ between lax loose transformations, 
  we construct a modulation of modules $\bar\mu:M_\sigma \To M_\tau$.}

  We define the components $\bar\mu_m$ of the intended modulation $\bar\mu$ as follows:
  \[\begin{tikzcd}
    Fx & Fy & Gy & \\
    Fx & Fy & Fy & Gy \\
    Fx & Fy & Gy
    \arrow["Fm"{inner sep=.8ex}, "\shortmid"{marking}, from=1-1, to=1-2]
    \arrow["{M_\sigma m}"{inner sep=.8ex}, "\shortmid"{marking}, curve={height=-18pt}, from=1-1, to=1-3]
    \arrow[equals, from=1-1, to=2-1]
    \arrow[""{name=0, anchor=center, inner sep=0}, "{\sigma_y}"{inner sep=.8ex}, "\shortmid"{marking}, from=1-2, to=1-3]
    \arrow[equals, from=1-2, to=2-2]
    \arrow[equals, from=1-3, to=2-4]
    \arrow["Fm"'{inner sep=.8ex}, "\shortmid"{marking}, from=2-1, to=2-2]
    \arrow[equals, from=2-1, to=3-1]
    \arrow["{F(\mathrm{id}_y)}"'{inner sep=.8ex}, "\shortmid"{marking}, from=2-2, to=2-3]
    \arrow["{\tau_y}"'{inner sep=.8ex}, "\shortmid"{marking}, from=2-3, to=2-4]
    \arrow[curve={height=-6pt}, equals, from=2-3, to=3-2]
    \arrow[equals, from=2-4, to=3-3]
    \arrow[""{name=1, anchor=center, inner sep=0}, "Fm"'{inner sep=.8ex}, "\shortmid"{marking}, from=3-1, to=3-2]
    \arrow["{M_\tau m}"'{inner sep=.8ex}, "\shortmid"{marking}, curve={height=18pt}, from=3-1, to=3-3]
    \arrow["{\tau_y}"'{inner sep=.8ex}, "\shortmid"{marking}, from=3-2, to=3-3]
    \arrow["{\mu_y}"{description}, draw=none, from=0, to=2-3]
    \arrow["{F_{m,\mathrm{id}_y}}"{description}, draw=none, from=2-2, to=1]
  \end{tikzcd}\]

  \item \emph{Given a modulation $\mu:M_\sigma\To M_\tau:F\proto G$ between modules, we construct a modulation
  of loose transformations $\underline{\mu}:\sigma\To \tau.$}

  We define the components $\underline{\mu}_x$ of the intended modulation
  $\underline{\mu}:\sigma\To \tau$ by whiskering $\mu_{\mathrm{id}_x}$ with a unitor of $F$ as follows:
  \[\begin{tikzcd}
    Fx & Fx & Gx \\
    Fx & Fx & Gx \\
    Fx & Fx & Gx
    \arrow[""{name=0, anchor=center, inner sep=0}, "\shortmid"{marking}, equals, from=1-1, to=1-2]
    \arrow[equals, from=1-1, to=2-1]
    \arrow["{\sigma_x}"'{inner sep=.8ex}, "\shortmid"{marking}, from=1-2, to=1-3]
    \arrow[equals, from=1-2, to=2-2]
    \arrow[equals, from=1-3, to=2-3]
    \arrow[""{name=1, anchor=center, inner sep=0}, "{F\mathrm{id}_x}"'{inner sep=.8ex}, "\shortmid"{marking}, from=2-1, to=2-2]
    \arrow[equals, from=2-1, to=3-1]
    \arrow["{\sigma_x}"'{inner sep=.8ex}, "\shortmid"{marking}, from=2-2, to=2-3]
    \arrow["{\mu_{\mathrm{id}_x}}"{description}, draw=none, from=2-2, to=3-2]
    \arrow[equals, from=2-3, to=3-3]
    \arrow["{F\mathrm{id}_x}"'{inner sep=.8ex}, "\shortmid"{marking}, from=3-1, to=3-2]
    \arrow["{\tau_x}"'{inner sep=.8ex}, "\shortmid"{marking}, from=3-2, to=3-3]
    \arrow["{F_x}"{description}, draw=none, from=0, to=1]
  \end{tikzcd}\]
  \end{enumerate}
\end{construction}

\begin{lemma}\label{lemma:constructions-work}
  The preceding constructions are well-defined. That is:
  \begin{enumerate}
  \item Given a lax loose transformation $\tau: F \To G:\dbl{D}\to \dbl{E}$ between lax double functors, 
  the construction $M_\tau$ above indeed gives a module.
  \item Given a module $M:F\proTo G:\dbl{D}\to \dbl{E}$ such that $M^\ell$ is invertible,
  the construction $\tau^M$ above indeed gives a loose transformation.
  \item Given a modulation $\mu:\sigma\To \tau:F\to G$ between lax loose transformations, 
  the construction $\bar\mu$ above indeed gives a modulation of modules.
  \item Given a modulation $\mu:M_\sigma\To M_\tau:F\proto G$ of modules, the construction
  $\underline{\mu}$ above indeed gives a modulation of loose transformations.
  \end{enumerate}
\end{lemma}
\begin{proof}

  \begin{enumerate}
  \item 
  Given a cell $\xinlinecell{x}{y}{x'}{y'}{f}{g}{m}{m'}{\alpha}$ in $\dbl{D}$
    consecutive with $\xinlinecell{y}{z}{y'}{z'}{g}{h}{n}{n'}{\beta},$ the naturality of the left action uses only the naturality of laxators of $F:$ 
    \[\begin{tikzcd}[column sep=small]
      Fx & Fy & Fz & Gz && Fx & Fy & Fz & Gz \\
      {Fx'} & {Fy'} & {Fz'} & {Gz'} && Fx && Fz & Gz \\
      {Fx'} && {Fz'} & {Gz'} && {Fx'} && {Fz'} & {Gz'}
      \arrow[""{name=0, anchor=center, inner sep=0}, "Fm"{inner sep=.8ex}, "\shortmid"{marking}, from=1-1, to=1-2]
      \arrow["Ff"', from=1-1, to=2-1]
      \arrow[""{name=1, anchor=center, inner sep=0}, "Fn"{inner sep=.8ex}, "\shortmid"{marking}, from=1-2, to=1-3]
      \arrow["Fg"{description}, from=1-2, to=2-2]
      \arrow[""{name=2, anchor=center, inner sep=0}, "{\tau_z}"{inner sep=.8ex}, "\shortmid"{marking}, from=1-3, to=1-4]
      \arrow["Fh"{description}, from=1-3, to=2-3]
      \arrow["Gh", from=1-4, to=2-4]
      \arrow["Fm"{inner sep=.8ex}, "\shortmid"{marking}, from=1-6, to=1-7]
      \arrow[equals, from=1-6, to=2-6]
      \arrow["Fn"{inner sep=.8ex}, "\shortmid"{marking}, from=1-7, to=1-8]
      \arrow["{\tau_z}"{inner sep=.8ex}, "\shortmid"{marking}, from=1-8, to=1-9]
      \arrow[equals, from=1-8, to=2-8]
      \arrow[equals, from=1-9, to=2-9]
      \arrow[""{name=3, anchor=center, inner sep=0}, "{Fm'}"{inner sep=.8ex}, "\shortmid"{marking}, from=2-1, to=2-2]
      \arrow[equals, from=2-1, to=3-1]
      \arrow[""{name=4, anchor=center, inner sep=0}, "{Fn'}"{inner sep=.8ex}, "\shortmid"{marking}, from=2-2, to=2-3]
      \arrow[""{name=5, anchor=center, inner sep=0}, "{\tau_{z'}}"{inner sep=.8ex}, "\shortmid"{marking}, from=2-3, to=2-4]
      \arrow[equals, from=2-3, to=3-3]
      \arrow[between={0.4}{0.6}, equals, from=2-4, to=2-6]
      \arrow[equals, from=2-4, to=3-4]
      \arrow[""{name=6, anchor=center, inner sep=0}, "{F(m\odot n)}"', from=2-6, to=2-8]
      \arrow["Ff"', from=2-6, to=3-6]
      \arrow[""{name=7, anchor=center, inner sep=0}, "{\tau_z}", from=2-8, to=2-9]
      \arrow["Fh"{description}, from=2-8, to=3-8]
      \arrow["Gh", from=2-9, to=3-9]
      \arrow[""{name=8, anchor=center, inner sep=0}, "{F(m'\odot n')}"', from=3-1, to=3-3]
      \arrow["{\tau_{z'}}", from=3-3, to=3-4]
      \arrow[""{name=9, anchor=center, inner sep=0}, "{F(m'\odot n')}"', from=3-6, to=3-8]
      \arrow[""{name=10, anchor=center, inner sep=0}, "{\tau_{z'}}"', from=3-8, to=3-9]
      \arrow["{F\alpha}"{description}, draw=none, from=0, to=3]
      \arrow["{F\beta}"{description}, draw=none, from=1, to=4]
      \arrow["{\tau_h}"{description}, draw=none, from=2, to=5]
      \arrow["{F_{m,n}}"{description}, draw=none, from=1-7, to=6]
      \arrow["{F_{m',n'}}"{description}, draw=none, from=2-2, to=8]
      \arrow["{F(\alpha\odot \beta)}"{description}, draw=none, from=6, to=9]
      \arrow["{\tau_h}"{description}, draw=none, from=7, to=10]
    \end{tikzcd}\]
    The naturality of the right action requires one extra step, applying the naturality of $\tau_n$ with 
    respect to the cell $\beta:$
\[\begin{tikzcd}
	Fx & Fy & Gy & Gz & Fx & Fy & Gy & Gz \\
	{Fx'} & {Fy'} & {Gy'} & {Gz'} & Fx & Fy & Fz & Gz \\
	{Fx'} & {Fy'} & {Fz'} & {Gz'} & Fx && Fz & Gz \\
	{Fx'} && {Fz'} & {Gz'} & {Fx'} && {Fz'} & {Gz'} \\
	&& Fx & Fy & Gy & Gz \\
	&& Fx & Fy & Fz & Gz \\
	&& {Fx'} & {Fy'} & {Fz'} & {Gz'} \\
	&& {Fx'} && {Fz'} & {Gz'}
	\arrow[""{name=0, anchor=center, inner sep=0}, "Fm"{inner sep=.8ex}, "\shortmid"{marking}, from=1-1, to=1-2]
	\arrow["Ff"', from=1-1, to=2-1]
	\arrow[""{name=1, anchor=center, inner sep=0}, "{{\tau_y}}"{inner sep=.8ex}, "\shortmid"{marking}, from=1-2, to=1-3]
	\arrow["Fg"{description}, from=1-2, to=2-2]
	\arrow[""{name=2, anchor=center, inner sep=0}, "Gn"{inner sep=.8ex}, "\shortmid"{marking}, from=1-3, to=1-4]
	\arrow["Gg"{description}, from=1-3, to=2-3]
	\arrow["Gh", from=1-4, to=2-4]
	\arrow["Fm"{inner sep=.8ex}, "\shortmid"{marking}, from=1-5, to=1-6]
	\arrow[equals, from=1-5, to=2-5]
	\arrow["{{{\tau_y}}}"{inner sep=.8ex}, "\shortmid"{marking}, from=1-6, to=1-7]
	\arrow[equals, from=1-6, to=2-6]
	\arrow["Gn"{inner sep=.8ex}, "\shortmid"{marking}, from=1-7, to=1-8]
	\arrow["{{{\tau_n}}}"{description}, draw=none, from=1-7, to=2-7]
	\arrow[equals, from=1-8, to=2-8]
	\arrow[""{name=3, anchor=center, inner sep=0}, "{{Fm'}}"'{inner sep=.8ex}, "\shortmid"{marking}, from=2-1, to=2-2]
	\arrow[equals, from=2-1, to=3-1]
	\arrow[""{name=4, anchor=center, inner sep=0}, "{{\tau_{y'}}}"'{inner sep=.8ex}, "\shortmid"{marking}, from=2-2, to=2-3]
	\arrow[equals, from=2-2, to=3-2]
	\arrow[""{name=5, anchor=center, inner sep=0}, "{{Gn'}}"'{inner sep=.8ex}, "\shortmid"{marking}, from=2-3, to=2-4]
	\arrow["{{{\tau_{n'}}}}"{description}, draw=none, from=2-3, to=3-3]
	\arrow[equals, from=2-4, to=3-4]
	\arrow["Fm"{inner sep=.8ex}, "\shortmid"{marking}, from=2-5, to=2-6]
	\arrow[equals, from=2-5, to=3-5]
	\arrow["Fn"{inner sep=.8ex}, "\shortmid"{marking}, from=2-6, to=2-7]
	\arrow["{{{\tau_z}}}"{inner sep=.8ex}, "\shortmid"{marking}, from=2-7, to=2-8]
	\arrow[equals, from=2-7, to=3-7]
	\arrow[equals, from=2-8, to=3-8]
	\arrow["{{Fm'}}"{inner sep=.8ex}, "\shortmid"{marking}, from=3-1, to=3-2]
	\arrow[equals, from=3-1, to=4-1]
	\arrow["{{Fn'}}"{inner sep=.8ex}, "\shortmid"{marking}, from=3-2, to=3-3]
	\arrow["{{{\tau_{z'}}}}"{inner sep=.8ex}, "\shortmid"{marking}, from=3-3, to=3-4]
	\arrow[equals, from=3-3, to=4-3]
	\arrow[equals, from=3-4, to=4-4]
	\arrow[""{name=6, anchor=center, inner sep=0}, "{{{F(m\odot n)}}}"'{inner sep=.8ex}, "\shortmid"{marking}, from=3-5, to=3-7]
	\arrow["Ff"', from=3-5, to=4-5]
	\arrow["{{{\tau_z}}}", from=3-7, to=3-8]
	\arrow["Fh"{description}, from=3-7, to=4-7]
	\arrow["Gh", from=3-8, to=4-8]
	\arrow[""{name=7, anchor=center, inner sep=0}, "{{{F(m'\odot n')}}}"'{inner sep=.8ex}, "\shortmid"{marking}, from=4-1, to=4-3]
	\arrow["{{{\tau_{z'}}}}", from=4-3, to=4-4]
	\arrow[""{name=8, anchor=center, inner sep=0}, "{{F(m'\odot n')}}"', from=4-5, to=4-7]
	\arrow["{{\tau_{z'}}}"', from=4-7, to=4-8]
	\arrow["Fm"{inner sep=.8ex}, "\shortmid"{marking}, from=5-3, to=5-4]
	\arrow[equals, from=5-3, to=6-3]
	\arrow[""{name=9, anchor=center, inner sep=0}, "{{\tau_y}}"{inner sep=.8ex}, "\shortmid"{marking}, from=5-4, to=5-5]
	\arrow[equals, from=5-4, to=6-4]
	\arrow["Gn"{inner sep=.8ex}, "\shortmid"{marking}, from=5-5, to=5-6]
	\arrow["{{\tau_n}}"{description}, draw=none, from=5-5, to=6-5]
	\arrow[equals, from=5-6, to=6-6]
	\arrow[""{name=10, anchor=center, inner sep=0}, "Fm"{inner sep=.8ex}, "\shortmid"{marking}, from=6-3, to=6-4]
	\arrow["Ff"', from=6-3, to=7-3]
	\arrow[""{name=11, anchor=center, inner sep=0}, "Fn"{inner sep=.8ex}, "\shortmid"{marking}, from=6-4, to=6-5]
	\arrow["Fg"{description}, from=6-4, to=7-4]
	\arrow[""{name=12, anchor=center, inner sep=0}, "{{\tau_z}}"{inner sep=.8ex}, "\shortmid"{marking}, from=6-5, to=6-6]
	\arrow["Fh"{description}, from=6-5, to=7-5]
	\arrow["Gh", from=6-6, to=7-6]
	\arrow[""{name=13, anchor=center, inner sep=0}, "{{Fm'}}"{inner sep=.8ex}, "\shortmid"{marking}, from=7-3, to=7-4]
	\arrow[equals, from=7-3, to=8-3]
	\arrow[""{name=14, anchor=center, inner sep=0}, "{{Fn'}}"{inner sep=.8ex}, "\shortmid"{marking}, from=7-4, to=7-5]
	\arrow[""{name=15, anchor=center, inner sep=0}, "{{{\tau_{z'}}}}"{inner sep=.8ex}, "\shortmid"{marking}, from=7-5, to=7-6]
	\arrow[equals, from=7-5, to=8-5]
	\arrow[equals, from=7-6, to=8-6]
	\arrow[""{name=16, anchor=center, inner sep=0}, "{{{F(m'\odot n')}}}"'{inner sep=.8ex}, "\shortmid"{marking}, from=8-3, to=8-5]
	\arrow["{{{\tau_{z'}}}}", from=8-5, to=8-6]
	\arrow["{{F\alpha}}"{description}, draw=none, from=0, to=3]
	\arrow["{{\tau_g}}"{description}, draw=none, from=1, to=4]
	\arrow["{{G\beta}}"{description}, draw=none, from=2, to=5]
	\arrow["{{{F_{m,n}}}}"{description}, draw=none, from=2-6, to=6]
	\arrow["{{{F_{m',n'}}}}"{description}, draw=none, from=3-2, to=7]
	\arrow["{{F(\alpha\odot\beta)}}"{description}, draw=none, from=6, to=8]
	\arrow[between={0.3}{0.7}, squiggly, from=4-3, to=9]
	\arrow[between={0.3}{0.7}, squiggly, from=9, to=8]
	\arrow["{{F\alpha}}"{description}, draw=none, from=10, to=13]
	\arrow["{{F\beta}}"{description}, draw=none, from=11, to=14]
	\arrow["{{\tau_h}}"{description}, draw=none, from=12, to=15]
	\arrow["{{{F_{m',n'}}}}"{description}, draw=none, from=7-4, to=16]
\end{tikzcd}\]
  
    For associativity of the actions, we consider a composable triple $x\xproto{m}y\xproto{n}z\xproto{p}w$ in $\dbl{D}.$ Associativity of the
  left action, again, reduces to associativity of $F$'s laxators.
  \[\begin{tikzcd}[column sep=small]
    Fx & Fy & Fz & Fw & Gw && Fx & Fy & Fz & Fw & Gw \\
    Fx & Fy && Fw & Gw & {=} & Fx && Fz & Fw & Gw \\
    Fx &&& Fw & Gw && Fx &&& Fw & Gw
    \arrow["Fm"{inner sep=.8ex}, "\shortmid"{marking}, from=1-1, to=1-2]
    \arrow[equals, from=1-1, to=2-1]
    \arrow["Fn"{inner sep=.8ex}, "\shortmid"{marking}, from=1-2, to=1-3]
    \arrow[equals, from=1-2, to=2-2]
    \arrow["Fp"{inner sep=.8ex}, "\shortmid"{marking}, from=1-3, to=1-4]
    \arrow["{\tau_w}"{inner sep=.8ex}, "\shortmid"{marking}, from=1-4, to=1-5]
    \arrow[equals, from=1-4, to=2-4]
    \arrow[equals, from=1-5, to=2-5]
    \arrow["Fm"{inner sep=.8ex}, "\shortmid"{marking}, from=1-7, to=1-8]
    \arrow[equals, from=1-7, to=2-7]
    \arrow["Fn"{inner sep=.8ex}, "\shortmid"{marking}, from=1-8, to=1-9]
    \arrow["Fp"{inner sep=.8ex}, "\shortmid"{marking}, from=1-9, to=1-10]
    \arrow[equals, from=1-9, to=2-9]
    \arrow["{\tau_w}"{inner sep=.8ex}, "\shortmid"{marking}, from=1-10, to=1-11]
    \arrow[equals, from=1-10, to=2-10]
    \arrow[equals, from=1-11, to=2-11]
    \arrow["Fm"'{inner sep=.8ex}, "\shortmid"{marking}, from=2-1, to=2-2]
    \arrow[equals, from=2-1, to=3-1]
    \arrow[""{name=0, anchor=center, inner sep=0}, "{F(n\odot p)}"'{inner sep=.8ex}, "\shortmid"{marking}, from=2-2, to=2-4]
    \arrow["{\tau_w}"{inner sep=.8ex}, "\shortmid"{marking}, from=2-4, to=2-5]
    \arrow[equals, from=2-4, to=3-4]
    \arrow[equals, from=2-5, to=3-5]
    \arrow[""{name=1, anchor=center, inner sep=0}, "{F(m\odot n)}"'{inner sep=.8ex}, "\shortmid"{marking}, from=2-7, to=2-9]
    \arrow[equals, from=2-7, to=3-7]
    \arrow["Fp"'{inner sep=.8ex}, "\shortmid"{marking}, from=2-9, to=2-10]
    \arrow["{\tau_w}", "\shortmid"{marking}, from=2-10, to=2-11]
    \arrow[equals, from=2-10, to=3-10]
    \arrow[equals, from=2-11, to=3-11]
    \arrow[""{name=2, anchor=center, inner sep=0}, "{F(m\odot n\odot p)}"'{inner sep=.8ex}, "\shortmid"{marking}, from=3-1, to=3-4]
    \arrow["{\tau_w}"'{inner sep=.8ex}, "\shortmid"{marking}, from=3-4, to=3-5]
    \arrow[""{name=3, anchor=center, inner sep=0}, "{F(m\odot n\odot p)}"', from=3-7, to=3-10]
    \arrow["{\tau_w}"', from=3-10, to=3-11]
    \arrow["{F_{n,p}}"{description}, draw=none, from=1-3, to=0]
    \arrow["{F_{m,n}}"{description}, draw=none, from=1-8, to=1]
    \arrow["{F_{m,n\odot p}}"{description}, draw=none, from=0, to=2]
    \arrow["{F_{m\odot n,p}}"{description}, draw=none, from=2-9, to=3]
  \end{tikzcd}\]
    
  Associativity of the right action requires us to also 
    slide in a pasting before applying associativity of $F$'s laxtors and to apply coherence of $\tau$ with $G$'s laxators
    afterward: 
\[\begin{tikzcd}[column sep=small]
	Fx & Fy & Gy & Gz & Gw &&&&&& \\
	Fx & Fy & Fz & Gz & Gw && Fx & Fy & Gy & Gz & Gw \\
	Fx && Fz & Gz & Gw && Fx & Fy & Gy && Gw \\
	Fx && Fz & Fw & Gw && Fx & Fy && Fw & Gw \\
	Fx &&& Fw & Gw && Fx &&& Fw & Gw \\
	\\
	Fx & Fy & Gy & Gz & Gw && Fx & Fy & Gy & Gz & Gw \\
	Fx & Fy & Fz & Gz & Gw && Fx & Fy & Fz & Gz & Gw \\
	Fx & Fy & Fz & Fw & Gw && Fx & Fy & Fz & Fw & Gw \\
	Fx && Fz & Fw & Gw && Fx & Fy && Fw & Gw \\
	Fx &&& Fw & Gw && Fx &&& Fw & Gw
	\arrow["Fm"{inner sep=.8ex}, "\shortmid"{marking}, from=1-1, to=1-2]
	\arrow[equals, from=1-1, to=2-1]
	\arrow["{{\tau_y}}"{inner sep=.8ex}, "\shortmid"{marking}, from=1-2, to=1-3]
	\arrow[equals, from=1-2, to=2-2]
	\arrow["Gn"{inner sep=.8ex}, "\shortmid"{marking}, from=1-3, to=1-4]
	\arrow["{{\tau_n}}"{description}, draw=none, from=1-3, to=2-3]
	\arrow["Gp"{inner sep=.8ex}, "\shortmid"{marking}, from=1-4, to=1-5]
	\arrow[equals, from=1-4, to=2-4]
	\arrow[equals, from=1-5, to=2-5]
	\arrow["Fm"'{inner sep=.8ex}, "\shortmid"{marking}, from=2-1, to=2-2]
	\arrow[equals, from=2-1, to=3-1]
	\arrow["Fn"'{inner sep=.8ex}, "\shortmid"{marking}, from=2-2, to=2-3]
	\arrow["{{\tau_z}}"'{inner sep=.8ex}, "\shortmid"{marking}, from=2-3, to=2-4]
	\arrow[equals, from=2-3, to=3-3]
	\arrow["Gp"'{inner sep=.8ex}, "\shortmid"{marking}, from=2-4, to=2-5]
	\arrow[equals, from=2-4, to=3-4]
	\arrow[equals, from=2-5, to=3-5]
	\arrow["Fm"{inner sep=.8ex}, "\shortmid"{marking}, from=2-7, to=2-8]
	\arrow[equals, from=2-7, to=3-7]
	\arrow["{{\tau_y}}"{inner sep=.8ex}, "\shortmid"{marking}, from=2-8, to=2-9]
	\arrow[equals, from=2-8, to=3-8]
	\arrow["Gn"{inner sep=.8ex}, "\shortmid"{marking}, from=2-9, to=2-10]
	\arrow[equals, from=2-9, to=3-9]
	\arrow["Gp"{inner sep=.8ex}, "\shortmid"{marking}, from=2-10, to=2-11]
	\arrow[equals, from=2-11, to=3-11]
	\arrow[""{name=0, anchor=center, inner sep=0}, "{{F(m\odot n)}}"'{inner sep=.8ex}, "\shortmid"{marking}, from=3-1, to=3-3]
	\arrow[equals, from=3-1, to=4-1]
	\arrow["{{\tau_z}}"'{inner sep=.8ex}, "\shortmid"{marking}, from=3-3, to=3-4]
	\arrow[equals, from=3-3, to=4-3]
	\arrow["Gp"'{inner sep=.8ex}, "\shortmid"{marking}, from=3-4, to=3-5]
	\arrow["{{\tau_p}}"{description}, draw=none, from=3-4, to=4-4]
	\arrow[equals, from=3-5, to=4-5]
	\arrow["Fm"'{inner sep=.8ex}, "\shortmid"{marking}, from=3-7, to=3-8]
	\arrow[equals, from=3-7, to=4-7]
	\arrow["{{\tau_y}}"'{inner sep=.8ex}, "\shortmid"{marking}, from=3-8, to=3-9]
	\arrow[equals, from=3-8, to=4-8]
	\arrow[""{name=1, anchor=center, inner sep=0}, "{{G(n\odot p)}}"'{inner sep=.8ex}, "\shortmid"{marking}, from=3-9, to=3-11]
	\arrow["{{\tau_{n\odot p}}}"{description}, draw=none, from=3-9, to=4-10]
	\arrow[equals, from=3-11, to=4-11]
	\arrow["{{F(m\odot n)}}"'{inner sep=.8ex}, "\shortmid"{marking}, from=4-1, to=4-3]
	\arrow[equals, from=4-1, to=5-1]
	\arrow["Fp"'{inner sep=.8ex}, "\shortmid"{marking}, from=4-3, to=4-4]
	\arrow["{{\tau_w}}"'{inner sep=.8ex}, "\shortmid"{marking}, from=4-4, to=4-5]
	\arrow[equals, from=4-4, to=5-4]
	\arrow[equals, from=4-5, to=5-5]
	\arrow["Fm"'{inner sep=.8ex}, "\shortmid"{marking}, from=4-7, to=4-8]
	\arrow[equals, from=4-7, to=5-7]
	\arrow["{{F(n\odot p)}}"'{inner sep=.8ex}, "\shortmid"{marking}, from=4-8, to=4-10]
	\arrow["{{\tau_w}}"'{inner sep=.8ex}, "\shortmid"{marking}, from=4-10, to=4-11]
	\arrow[equals, from=4-10, to=5-10]
	\arrow[equals, from=4-11, to=5-11]
	\arrow[""{name=2, anchor=center, inner sep=0}, "{{F(m\odot n\odot p)}}"'{inner sep=.8ex}, "\shortmid"{marking}, from=5-1, to=5-4]
	\arrow["{{\tau_w}}"'{inner sep=.8ex}, "\shortmid"{marking}, from=5-4, to=5-5]
	\arrow[""{name=3, anchor=center, inner sep=0}, "{{F(m\odot n\odot p)}}"'{inner sep=.8ex}, "\shortmid"{marking}, from=5-7, to=5-10]
	\arrow["{{\tau_w}}"'{inner sep=.8ex}, "\shortmid"{marking}, from=5-10, to=5-11]
	\arrow["Fm"{inner sep=.8ex}, "\shortmid"{marking}, from=7-1, to=7-2]
	\arrow[equals, from=7-1, to=8-1]
	\arrow["{{\tau_y}}"{inner sep=.8ex}, "\shortmid"{marking}, from=7-2, to=7-3]
	\arrow[equals, from=7-2, to=8-2]
	\arrow["Gn"{inner sep=.8ex}, "\shortmid"{marking}, from=7-3, to=7-4]
	\arrow["{{\tau_n}}"{description}, draw=none, from=7-3, to=8-3]
	\arrow["Gp"{inner sep=.8ex}, "\shortmid"{marking}, from=7-4, to=7-5]
	\arrow[equals, from=7-4, to=8-4]
	\arrow[equals, from=7-5, to=8-5]
	\arrow["Fm"{inner sep=.8ex}, "\shortmid"{marking}, from=7-7, to=7-8]
	\arrow[equals, from=7-7, to=8-7]
	\arrow["{{\tau_y}}"{inner sep=.8ex}, "\shortmid"{marking}, from=7-8, to=7-9]
	\arrow[equals, from=7-8, to=8-8]
	\arrow[between={0.3}{0.7}, squiggly, from=7-9, to=5-10]
	\arrow["Gn"{inner sep=.8ex}, "\shortmid"{marking}, from=7-9, to=7-10]
	\arrow["{{\tau_n}}"{description}, draw=none, from=7-9, to=8-9]
	\arrow["Gp"{inner sep=.8ex}, "\shortmid"{marking}, from=7-10, to=7-11]
	\arrow[equals, from=7-10, to=8-10]
	\arrow[equals, from=7-11, to=8-11]
	\arrow["Fm"'{inner sep=.8ex}, "\shortmid"{marking}, from=8-1, to=8-2]
	\arrow[equals, from=8-1, to=9-1]
	\arrow["Fn"'{inner sep=.8ex}, "\shortmid"{marking}, from=8-2, to=8-3]
	\arrow[equals, from=8-2, to=9-2]
	\arrow["{{\tau_z}}"'{inner sep=.8ex}, "\shortmid"{marking}, from=8-3, to=8-4]
	\arrow[equals, from=8-3, to=9-3]
	\arrow["Gp"'{inner sep=.8ex}, "\shortmid"{marking}, from=8-4, to=8-5]
	\arrow["{{\tau_p}}"{description}, draw=none, from=8-4, to=9-4]
	\arrow[equals, from=8-5, to=9-5]
	\arrow["Fm"'{inner sep=.8ex}, "\shortmid"{marking}, from=8-7, to=8-8]
	\arrow[equals, from=8-7, to=9-7]
	\arrow["Fn"'{inner sep=.8ex}, "\shortmid"{marking}, from=8-8, to=8-9]
	\arrow[equals, from=8-8, to=9-8]
	\arrow["{{\tau_z}}"'{inner sep=.8ex}, "\shortmid"{marking}, from=8-9, to=8-10]
	\arrow[equals, from=8-9, to=9-9]
	\arrow["Gp"'{inner sep=.8ex}, "\shortmid"{marking}, from=8-10, to=8-11]
	\arrow["{{\tau_p}}"{description}, draw=none, from=8-10, to=9-10]
	\arrow[equals, from=8-11, to=9-11]
	\arrow["Fm"'{inner sep=.8ex}, "\shortmid"{marking}, from=9-1, to=9-2]
	\arrow[equals, from=9-1, to=10-1]
	\arrow["Fn"'{inner sep=.8ex}, "\shortmid"{marking}, from=9-2, to=9-3]
	\arrow["Fp"'{inner sep=.8ex}, "\shortmid"{marking}, from=9-3, to=9-4]
	\arrow[equals, from=9-3, to=10-3]
	\arrow["{{\tau_w}}"'{inner sep=.8ex}, "\shortmid"{marking}, from=9-4, to=9-5]
	\arrow[equals, from=9-4, to=10-4]
	\arrow[between={0.3}{0.7}, squiggly, from=9-5, to=9-7]
	\arrow[equals, from=9-5, to=10-5]
	\arrow["Fm"'{inner sep=.8ex}, "\shortmid"{marking}, from=9-7, to=9-8]
	\arrow[equals, from=9-7, to=10-7]
	\arrow["Fn"'{inner sep=.8ex}, "\shortmid"{marking}, from=9-8, to=9-9]
	\arrow[equals, from=9-8, to=10-8]
	\arrow["Fp"'{inner sep=.8ex}, "\shortmid"{marking}, from=9-9, to=9-10]
	\arrow["{{\tau_w}}"'{inner sep=.8ex}, "\shortmid"{marking}, from=9-10, to=9-11]
	\arrow[equals, from=9-10, to=10-10]
	\arrow[equals, from=9-11, to=10-11]
	\arrow[""{name=4, anchor=center, inner sep=0}, "{{F(m\odot n)}}"'{inner sep=.8ex}, "\shortmid"{marking}, from=10-1, to=10-3]
	\arrow[equals, from=10-1, to=11-1]
	\arrow["Fp"'{inner sep=.8ex}, "\shortmid"{marking}, from=10-3, to=10-4]
	\arrow["{{\tau_w}}"'{inner sep=.8ex}, "\shortmid"{marking}, from=10-4, to=10-5]
	\arrow[equals, from=10-4, to=11-4]
	\arrow[equals, from=10-5, to=11-5]
	\arrow["Fm"'{inner sep=.8ex}, "\shortmid"{marking}, from=10-7, to=10-8]
	\arrow[equals, from=10-7, to=11-7]
	\arrow[""{name=5, anchor=center, inner sep=0}, "{{F(n\odot p)}}"'{inner sep=.8ex}, "\shortmid"{marking}, from=10-8, to=10-10]
	\arrow["{{\tau_w}}"'{inner sep=.8ex}, "\shortmid"{marking}, from=10-10, to=10-11]
	\arrow[equals, from=10-10, to=11-10]
	\arrow[equals, from=10-11, to=11-11]
	\arrow[""{name=6, anchor=center, inner sep=0}, "{{F(m\odot n\odot p)}}"'{inner sep=.8ex}, "\shortmid"{marking}, from=11-1, to=11-4]
	\arrow["{{\tau_w}}"'{inner sep=.8ex}, "\shortmid"{marking}, from=11-4, to=11-5]
	\arrow[""{name=7, anchor=center, inner sep=0}, "{{F(m\odot n\odot p)}}"'{inner sep=.8ex}, "\shortmid"{marking}, from=11-7, to=11-10]
	\arrow["{{\tau_w}}"'{inner sep=.8ex}, "\shortmid"{marking}, from=11-10, to=11-11]
	\arrow["{{F_{m,n}}}"{description}, draw=none, from=2-2, to=0]
	\arrow["{{G_{n,p}}}"{description}, draw=none, from=2-10, to=1]
	\arrow["{{F_{m\odot n,p}}}"{description}, draw=none, from=4-3, to=2]
	\arrow["{{F_{m,n\odot p}}}"{description}, draw=none, from=4-8, to=3]
	\arrow[between={0.4}{0.7}, squiggly, from=2, to=7-3]
	\arrow["{{F_{m,n}}}"{description}, draw=none, from=9-2, to=4]
	\arrow["{{F_{n,p}}}"{description}, draw=none, from=9-9, to=5]
	\arrow["{{F_{m\odot n,p}}}"{description}, draw=none, from=10-3, to=6]
	\arrow["{{F_{m,n\odot p}}}"{description}, draw=none, from=5, to=7]
\end{tikzcd}\]
  
  For unitality of the actions, we rely on unitality of $F,$ along with, for the right action, the coherence of $\tau$
    with unitors: 
    \[\begin{tikzcd}
      Fx & Fx & Fy & Gy && Fx & Fy & Gy \\
      Fx & Fx & Fy & Gy & {=} \\
      Fx && Fy & Gy && Fx & Fy & Gy \\
      Fx & Fy & Gy & Gy && Fx & Fy & Gy \\
      Fx & Fy & Gy & Gy & {=} \\
      Fx & Fy & Fy & Gy && Fx & Fy & Gy \\
      Fx && Fy & Gy
      \arrow[""{name=0, anchor=center, inner sep=0}, "\shortmid"{marking}, equals, from=1-1, to=1-2]
      \arrow[equals, from=1-1, to=2-1]
      \arrow["Fm"{inner sep=.8ex}, "\shortmid"{marking}, from=1-2, to=1-3]
      \arrow[equals, from=1-2, to=2-2]
      \arrow["{\tau_y}", from=1-3, to=1-4]
      \arrow[equals, from=1-3, to=2-3]
      \arrow[equals, from=1-4, to=2-4]
      \arrow["Fm"{inner sep=.8ex}, "\shortmid"{marking}, from=1-6, to=1-7]
      \arrow[equals, from=1-6, to=3-6]
      \arrow["{\tau_y}", from=1-7, to=1-8]
      \arrow[equals, from=1-8, to=3-8]
      \arrow[""{name=1, anchor=center, inner sep=0}, "{F(\mathrm{id}_x)}"'{inner sep=.8ex}, "\shortmid"{marking}, from=2-1, to=2-2]
      \arrow[equals, from=2-1, to=3-1]
      \arrow["Fm"'{inner sep=.8ex}, "\shortmid"{marking}, from=2-2, to=2-3]
      \arrow["{\tau_y}"'{inner sep=.8ex}, "\shortmid"{marking}, from=2-3, to=2-4]
      \arrow[equals, from=2-3, to=3-3]
      \arrow[equals, from=2-4, to=3-4]
      \arrow[""{name=2, anchor=center, inner sep=0}, "Fm"'{inner sep=.8ex}, "\shortmid"{marking}, from=3-1, to=3-3]
      \arrow["{\tau_y}"'{inner sep=.8ex}, "\shortmid"{marking}, from=3-3, to=3-4]
      \arrow["Fm"{inner sep=.8ex}, "\shortmid"{marking}, from=3-6, to=3-7]
      \arrow["{\tau_y}", from=3-7, to=3-8]
      \arrow["Fm"{inner sep=.8ex}, "\shortmid"{marking}, from=4-1, to=4-2]
      \arrow[equals, from=4-1, to=5-1]
      \arrow["{\tau_y}"{inner sep=.8ex}, "\shortmid"{marking}, from=4-2, to=4-3]
      \arrow[equals, from=4-2, to=5-2]
      \arrow["\shortmid"{marking}, equals, from=4-3, to=4-4]
      \arrow[equals, from=4-3, to=5-3]
      \arrow[equals, from=4-4, to=5-4]
      \arrow["Fm"{inner sep=.8ex}, "\shortmid"{marking}, from=4-6, to=4-7]
      \arrow[equals, from=4-6, to=6-6]
      \arrow["{\tau_y}", from=4-7, to=4-8]
      \arrow[equals, from=4-8, to=6-8]
      \arrow["Fm"'{inner sep=.8ex}, "\shortmid"{marking}, from=5-1, to=5-2]
      \arrow[equals, from=5-1, to=6-1]
      \arrow["{\tau_y}"'{inner sep=.8ex}, "\shortmid"{marking}, from=5-2, to=5-3]
      \arrow[equals, from=5-2, to=6-2]
      \arrow["{G(\mathrm{id}_y)}"'{inner sep=.8ex}, "\shortmid"{marking}, from=5-3, to=5-4]
      \arrow["{\tau_{\mathrm{id}_y}}"{description}, draw=none, from=5-3, to=6-3]
      \arrow[equals, from=5-4, to=6-4]
      \arrow["Fm"'{inner sep=.8ex}, "\shortmid"{marking}, from=6-1, to=6-2]
      \arrow[equals, from=6-1, to=7-1]
      \arrow["{F(\mathrm{id}_y)}"'{inner sep=.8ex}, "\shortmid"{marking}, from=6-2, to=6-3]
      \arrow["{\tau_y}"'{inner sep=.8ex}, "\shortmid"{marking}, from=6-3, to=6-4]
      \arrow[equals, from=6-3, to=7-3]
      \arrow[equals, from=6-4, to=7-4]
      \arrow["Fm"{inner sep=.8ex}, "\shortmid"{marking}, from=6-6, to=6-7]
      \arrow["{\tau_y}", from=6-7, to=6-8]
      \arrow[""{name=3, anchor=center, inner sep=0}, "Fm"'{inner sep=.8ex}, "\shortmid"{marking}, from=7-1, to=7-3]
      \arrow["{\tau_y}"'{inner sep=.8ex}, "\shortmid"{marking}, from=7-3, to=7-4]
      \arrow["{F_x}"{description}, draw=none, from=0, to=1]
      \arrow["{F_{\mathrm{id}_x,m}}"{description}, draw=none, from=2-2, to=2]
      \arrow["{F_{m,\mathrm{id}_y}}"{description}, draw=none, from=6-2, to=3]
    \end{tikzcd}\]

  \item We must check that the proposed $\tau$ intertwines the unitors and laxators of $F$ and $G$ 
  correctly. The pentagon for the laxators is shown to commute via this diagram, where we 
  make the first step by reversing the square for associativity of the left action of $M,$ 
  the second by using the commutativity of the left and right actions of $M,$ the t
  third by cancelling $M^\ell$ with $(M^\ell)^{-1},$ and the last using the associativity of the
  right action of M.
\[\begin{tikzcd}[column sep=small]
	Fx & Gx & Gy & Gz && Fx & Gx & Gy & Gz \\
	Fx && Gy &&& Fx & Gx && Gz \\
	Fx & Fy & Gy & Gz && Fx &&& Gz \\
	Fx & Fy && Gz && Fx && Fz & Gz \\
	Fx & Fy & Fz & Gz && Fx & Gx & Gy & Gz \\
	Fx && Fz & Gz && Fx && Gy & Gz \\
	Fx & Gx & Gy & Gz && Fx &&& Gz \\
	Fx && Gy &&& Fx && Fz & Gz \\
	Fx & Fy & Gy & Gz && Fx & Gx & Gy & Gz \\
	Fx & Fy && Gz && Fx && Gy \\
	Fx &&& Gz && Fx & Fy & Gy & Gz \\
	Fx && Fz & Gz && Fx && Gy & Gz \\
	&&&&& Fx &&& Gz \\
	&&&&& Fx && Fz & Gz
	\arrow["{{M\mathrm{id}_x}}"{inner sep=.8ex}, "\shortmid"{marking}, from=1-1, to=1-2]
	\arrow[equals, from=1-1, to=2-1]
	\arrow["Gm"{inner sep=.8ex}, "\shortmid"{marking}, from=1-2, to=1-3]
	\arrow["Gn"{inner sep=.8ex}, "\shortmid"{marking}, from=1-3, to=1-4]
	\arrow[equals, from=1-3, to=2-3]
	\arrow[equals, from=1-4, to=3-4]
	\arrow["{{M\mathrm{id}_x}}"{inner sep=.8ex}, "\shortmid"{marking}, from=1-6, to=1-7]
	\arrow[equals, from=1-6, to=2-6]
	\arrow["Gm"{inner sep=.8ex}, "\shortmid"{marking}, from=1-7, to=1-8]
	\arrow[equals, from=1-7, to=2-7]
	\arrow["Gn"{inner sep=.8ex}, "\shortmid"{marking}, from=1-8, to=1-9]
	\arrow[equals, from=1-9, to=2-9]
	\arrow[""{name=0, anchor=center, inner sep=0}, "Mm"'{inner sep=.8ex}, "\shortmid"{marking}, from=2-1, to=2-3]
	\arrow[equals, from=2-1, to=3-1]
	\arrow[equals, from=2-3, to=3-3]
	\arrow["{{M\mathrm{id}_x}}"'{inner sep=.8ex}, "\shortmid"{marking}, from=2-6, to=2-7]
	\arrow[equals, from=2-6, to=3-6]
	\arrow[""{name=1, anchor=center, inner sep=0}, "{{G(m\odot n)}}"'{inner sep=.8ex}, "\shortmid"{marking}, from=2-7, to=2-9]
	\arrow[equals, from=2-9, to=3-9]
	\arrow["Fm"'{inner sep=.8ex}, "\shortmid"{marking}, from=3-1, to=3-2]
	\arrow[equals, from=3-1, to=4-1]
	\arrow["{{M\mathrm{id}_y}}"'{inner sep=.8ex}, "\shortmid"{marking}, from=3-2, to=3-3]
	\arrow[equals, from=3-2, to=4-2]
	\arrow["Gn"'{inner sep=.8ex}, "\shortmid"{marking}, from=3-3, to=3-4]
	\arrow[equals, from=3-4, to=4-4]
	\arrow[""{name=2, anchor=center, inner sep=0}, "{{M(m\odot n)}}"'{inner sep=.8ex}, "\shortmid"{marking}, from=3-6, to=3-9]
	\arrow[equals, from=3-6, to=4-6]
	\arrow[equals, from=3-9, to=4-9]
	\arrow["Fm"'{inner sep=.8ex}, "\shortmid"{marking}, from=4-1, to=4-2]
	\arrow[equals, from=4-1, to=5-1]
	\arrow[""{name=3, anchor=center, inner sep=0}, "Mn"'{inner sep=.8ex}, "\shortmid"{marking}, from=4-2, to=4-4]
	\arrow[equals, from=4-2, to=5-2]
	\arrow[equals, from=4-4, to=5-4]
	\arrow["{{F(m\odot n)}}"', from=4-6, to=4-8]
	\arrow["{{M\mathrm{id}_z}}"', from=4-8, to=4-9]
	\arrow["Fm"'{inner sep=.8ex}, "\shortmid"{marking}, from=5-1, to=5-2]
	\arrow[equals, from=5-1, to=6-1]
	\arrow["Fn"'{inner sep=.8ex}, "\shortmid"{marking}, from=5-2, to=5-3]
	\arrow["{{M\mathrm{id}_z}}"', from=5-3, to=5-4]
	\arrow[equals, from=5-3, to=6-3]
	\arrow[equals, from=5-4, to=6-4]
	\arrow["{{M\mathrm{id}_x}}"{inner sep=.8ex}, "\shortmid"{marking}, from=5-6, to=5-7]
	\arrow[equals, from=5-6, to=6-6]
	\arrow["Gm"{inner sep=.8ex}, "\shortmid"{marking}, from=5-7, to=5-8]
	\arrow[squiggly, from=5-8, to=4-8]
	\arrow["Gn"{inner sep=.8ex}, "\shortmid"{marking}, from=5-8, to=5-9]
	\arrow[equals, from=5-8, to=6-8]
	\arrow[equals, from=5-9, to=6-9]
	\arrow[""{name=4, anchor=center, inner sep=0}, "{{F(m\odot n)}}"', from=6-1, to=6-3]
	\arrow["{{M\mathrm{id}_z}}"', from=6-3, to=6-4]
	\arrow[squiggly, from=6-3, to=7-3]
	\arrow[""{name=5, anchor=center, inner sep=0}, "Mm"'{inner sep=.8ex}, "\shortmid"{marking}, from=6-6, to=6-8]
	\arrow[equals, from=6-6, to=7-6]
	\arrow["Gn"'{inner sep=.8ex}, "\shortmid"{marking}, from=6-8, to=6-9]
	\arrow[equals, from=6-9, to=7-9]
	\arrow["{{M\mathrm{id}_x}}"{inner sep=.8ex}, "\shortmid"{marking}, from=7-1, to=7-2]
	\arrow[equals, from=7-1, to=8-1]
	\arrow["Gm"{inner sep=.8ex}, "\shortmid"{marking}, from=7-2, to=7-3]
	\arrow["Gn"{inner sep=.8ex}, "\shortmid"{marking}, from=7-3, to=7-4]
	\arrow[equals, from=7-3, to=8-3]
	\arrow[equals, from=7-4, to=9-4]
	\arrow[""{name=6, anchor=center, inner sep=0}, "{{M(m\odot n)}}"'{inner sep=.8ex}, "\shortmid"{marking}, from=7-6, to=7-9]
	\arrow[equals, from=7-6, to=8-6]
	\arrow[equals, from=7-9, to=8-9]
	\arrow[""{name=7, anchor=center, inner sep=0}, "Mm"'{inner sep=.8ex}, "\shortmid"{marking}, from=8-1, to=8-3]
	\arrow[equals, from=8-1, to=9-1]
	\arrow[equals, from=8-3, to=9-3]
	\arrow["{{F(m\odot n)}}"', from=8-6, to=8-8]
	\arrow["{{M\mathrm{id}_z}}"', from=8-8, to=8-9]
	\arrow["Fm"'{inner sep=.8ex}, "\shortmid"{marking}, from=9-1, to=9-2]
	\arrow[equals, from=9-1, to=10-1]
	\arrow["{{M\mathrm{id}_y}}"'{inner sep=.8ex}, "\shortmid"{marking}, from=9-2, to=9-3]
	\arrow[equals, from=9-2, to=10-2]
	\arrow["Gn"'{inner sep=.8ex}, "\shortmid"{marking}, from=9-3, to=9-4]
	\arrow[equals, from=9-4, to=10-4]
	\arrow["{{M\mathrm{id}_x}}"{inner sep=.8ex}, "\shortmid"{marking}, from=9-6, to=9-7]
	\arrow[equals, from=9-6, to=10-6]
	\arrow["Gm"{inner sep=.8ex}, "\shortmid"{marking}, from=9-7, to=9-8]
	\arrow[squiggly, from=9-8, to=8-8]
	\arrow["Gn"{inner sep=.8ex}, "\shortmid"{marking}, from=9-8, to=9-9]
	\arrow[equals, from=9-8, to=10-8]
	\arrow[equals, from=9-9, to=11-9]
	\arrow["Fm"'{inner sep=.8ex}, "\shortmid"{marking}, from=10-1, to=10-2]
	\arrow[equals, from=10-1, to=11-1]
	\arrow[""{name=8, anchor=center, inner sep=0}, "Mn"'{inner sep=.8ex}, "\shortmid"{marking}, from=10-2, to=10-4]
	\arrow[squiggly, from=10-4, to=10-6]
	\arrow[equals, from=10-4, to=11-4]
	\arrow[""{name=9, anchor=center, inner sep=0}, "Mm"'{inner sep=.8ex}, "\shortmid"{marking}, from=10-6, to=10-8]
	\arrow[equals, from=10-6, to=11-6]
	\arrow[equals, from=10-8, to=11-8]
	\arrow[""{name=10, anchor=center, inner sep=0}, "{{M(m\odot n)}}"'{inner sep=.8ex}, "\shortmid"{marking}, from=11-1, to=11-4]
	\arrow[equals, from=11-1, to=12-1]
	\arrow[equals, from=11-4, to=12-4]
	\arrow["Fm"'{inner sep=.8ex}, "\shortmid"{marking}, from=11-6, to=11-7]
	\arrow[equals, from=11-6, to=12-6]
	\arrow["{{M\mathrm{id}_y}}"'{inner sep=.8ex}, "\shortmid"{marking}, from=11-7, to=11-8]
	\arrow["Gn"'{inner sep=.8ex}, "\shortmid"{marking}, from=11-8, to=11-9]
	\arrow[equals, from=11-8, to=12-8]
	\arrow[equals, from=11-9, to=12-9]
	\arrow["{{F(m\odot n)}}"', from=12-1, to=12-3]
	\arrow["{{M\mathrm{id}_z}}"', from=12-3, to=12-4]
	\arrow[""{name=11, anchor=center, inner sep=0}, "Mm"'{inner sep=.8ex}, "\shortmid"{marking}, from=12-6, to=12-8]
	\arrow[equals, from=12-6, to=13-6]
	\arrow["Gn"'{inner sep=.8ex}, "\shortmid"{marking}, from=12-8, to=12-9]
	\arrow[equals, from=12-9, to=13-9]
	\arrow[""{name=12, anchor=center, inner sep=0}, "{{M(m\odot n)}}"'{inner sep=.8ex}, "\shortmid"{marking}, from=13-6, to=13-9]
	\arrow[equals, from=13-6, to=14-6]
	\arrow[equals, from=13-9, to=14-9]
	\arrow["{{F(m\odot n)}}"', from=14-6, to=14-8]
	\arrow["{{M\mathrm{id}_z}}"', from=14-8, to=14-9]
	\arrow["{{M^r_m}}"{description}, draw=none, from=1-2, to=0]
	\arrow["{{G_{m,n}}}"{description}, draw=none, from=1-8, to=1]
	\arrow["{{(M^\ell_m)^{-1}}}"{description, pos=0.7}, draw=none, from=0, to=3-2]
	\arrow["{{M^r_{m\odot n}}}"{description}, draw=none, from=2-7, to=2]
	\arrow["{{M^r_n}}"{description}, draw=none, from=3-3, to=3]
	\arrow["{{(M^\ell_{m\odot n})^{-1}}}"{description, pos=0.7}, draw=none, from=2, to=4-8]
	\arrow["{{(M^\ell_n)^{-1}}}"{description, pos=0.7}, draw=none, from=3, to=5-3]
	\arrow["{{F_{m,n}}}"{description}, draw=none, from=5-2, to=4]
	\arrow["{{M^r_m}}"{description}, draw=none, from=5-7, to=5]
	\arrow["{{M^r_{m,n}}}"{description}, draw=none, from=6-8, to=6]
	\arrow["{{M^r_m}}"{description}, draw=none, from=7-2, to=7]
	\arrow["{{(M^\ell_{m\odot n})^{-1}}}"{description, pos=0.7}, draw=none, from=6, to=8-8]
	\arrow["{{(M^\ell_m)^{-1}}}"{description, pos=0.7}, draw=none, from=7, to=9-2]
	\arrow["{{M^r_n}}"{description}, draw=none, from=9-3, to=8]
	\arrow["{{M^r_m}}"{description}, draw=none, from=9-7, to=9]
	\arrow["{{M^\ell_{m,n}}}"{description}, draw=none, from=10-2, to=10]
	\arrow["{{(M^\ell_m)^{-1}}}"{description, pos=0.7}, draw=none, from=9, to=11-7]
	\arrow["{{(M^\ell_{m\odot n})^{-1}}}"{description, pos=0.7}, draw=none, from=10, to=12-3]
	\arrow["{{M^\ell_m}}"{description}, draw=none, from=11-7, to=11]
	\arrow["{{M^r_{m,n}}}"{description}, draw=none, from=12-8, to=12]
	\arrow["{{(M^\ell_{m\odot n})^{-1}}}"{description, pos=0.7}, draw=none, from=12, to=14-8]
\end{tikzcd}\]
  For the unitors, we compose the desired square with the isomorphism 
  $M^\ell:F(\mathrm{id})\odot \tau^M_0\to \tau^M_0.$ This produces the following 
  list of equal pastings, with the moves being made by unitality of the actions of $M.$
\[\begin{tikzcd}
	Fx & Gx & Gx & Fx & Fx & Gx \\
	Fx & Gx & Gx & Fx & Fx & Gx \\
	Fx && Gx & Fx && Gx \\
	Fx & Fx & Gx & Fx & Gx \\
	Fx && Gx & Fx & Gx
	\arrow["{M\mathrm{id}_x}"{inner sep=.8ex}, "\shortmid"{marking}, from=1-1, to=1-2]
	\arrow[equals, from=1-1, to=2-1]
	\arrow[""{name=0, anchor=center, inner sep=0}, "\shortmid"{marking}, equals, from=1-2, to=1-3]
	\arrow[equals, from=1-2, to=2-2]
	\arrow[equals, from=1-3, to=2-3]
	\arrow[""{name=1, anchor=center, inner sep=0}, "\shortmid"{marking}, equals, from=1-4, to=1-5]
	\arrow[equals, from=1-4, to=2-4]
	\arrow["{M\mathrm{id}_x}"{inner sep=.8ex}, "\shortmid"{marking}, from=1-5, to=1-6]
	\arrow[equals, from=1-5, to=2-5]
	\arrow[equals, from=1-6, to=2-6]
	\arrow["{M\mathrm{id}_x}"'{inner sep=.8ex}, "\shortmid"{marking}, from=2-1, to=2-2]
	\arrow[equals, from=2-1, to=3-1]
	\arrow[""{name=2, anchor=center, inner sep=0}, "{G(\mathrm{id}_x)}"', from=2-2, to=2-3]
	\arrow[equals, from=2-3, to=3-3]
	\arrow[""{name=3, anchor=center, inner sep=0}, "{F(\mathrm{id}_x)}"'{inner sep=.8ex}, "\shortmid"{marking}, from=2-4, to=2-5]
	\arrow[equals, from=2-4, to=3-4]
	\arrow["{M\mathrm{id}_x}"'{inner sep=.8ex}, "\shortmid"{marking}, from=2-5, to=2-6]
	\arrow[equals, from=2-6, to=3-6]
	\arrow[""{name=4, anchor=center, inner sep=0}, "{M\mathrm{id_x}}"'{inner sep=.8ex}, "\shortmid"{marking}, from=3-1, to=3-3]
	\arrow[equals, from=3-1, to=4-1]
	\arrow[equals, from=3-3, to=4-3]
	\arrow[squiggly, from=3-3, to=4-4]
	\arrow[""{name=5, anchor=center, inner sep=0}, "{M\mathrm{id}_x}"'{inner sep=.8ex}, "\shortmid"{marking}, from=3-4, to=3-6]
	\arrow["{F(\mathrm{id}_x)}"'{inner sep=.8ex}, "\shortmid"{marking}, from=4-1, to=4-2]
	\arrow[equals, from=4-1, to=5-1]
	\arrow["{M\mathrm{id}_x}"'{inner sep=.8ex}, "\shortmid"{marking}, from=4-2, to=4-3]
	\arrow[equals, from=4-3, to=5-3]
	\arrow[""{name=6, anchor=center, inner sep=0}, "{M\mathrm{id_x}}"'{inner sep=.8ex}, "\shortmid"{marking}, from=4-4, to=4-5]
	\arrow[equals, from=4-4, to=5-4]
	\arrow[equals, from=4-5, to=5-5]
	\arrow[""{name=7, anchor=center, inner sep=0}, "{M\mathrm{id}_x}"'{inner sep=.8ex}, "\shortmid"{marking}, from=5-1, to=5-3]
	\arrow["{M\mathrm{id}_x}"'{inner sep=.8ex}, "\shortmid"{marking}, from=5-4, to=5-5]
	\arrow["{G_x}"{description}, draw=none, from=0, to=2]
	\arrow["{F_x}"{description}, draw=none, from=1, to=3]
	\arrow["{M^r}"{description}, draw=none, from=2-2, to=4]
	\arrow["{M^\ell}"{description}, draw=none, from=2-5, to=5]
	\arrow["{(M^\ell)^{-1}}"{description, pos=0.7}, draw=none, from=4, to=4-2]
	\arrow[between={0.2}{0.8}, squiggly, from=3-4, to=6]
	\arrow["{M^\ell}"{description}, draw=none, from=4-2, to=7]
\end{tikzcd}\]

  \item We must check that the proposed modulation correctly intertwines the right and left 
  actions of $M_\sigma$ and $M_\tau.$ For the left action, we need only associativity of laxators for 
  $F$ and some sliding: 
\[\begin{tikzcd}[column sep=small]
	& Fx & Fy & Fz & Gz &&& Fx & Fy & Fz & Gz \\
	Fx & Fy & Fz & Fz & Gz &&& Fx && Fz & Gz \\
	& Fx & Fy & Fz & Gz && Fx && Fz & Fz & Gz \\
	& Fx && Fz & Gz &&& Fx && Fz & Gz
	\arrow["Fm"{inner sep=.8ex}, "\shortmid"{marking}, from=1-2, to=1-3]
	\arrow[equals, from=1-2, to=2-1]
	\arrow["Fn"{inner sep=.8ex}, "\shortmid"{marking}, from=1-3, to=1-4]
	\arrow[equals, from=1-3, to=2-2]
	\arrow[""{name=0, anchor=center, inner sep=0}, "{{\sigma_z}}"{inner sep=.8ex}, "\shortmid"{marking}, from=1-4, to=1-5]
	\arrow[equals, from=1-4, to=2-3]
	\arrow[equals, from=1-5, to=2-5]
	\arrow["Fm"{inner sep=.8ex}, "\shortmid"{marking}, from=1-8, to=1-9]
	\arrow[equals, from=1-8, to=2-8]
	\arrow["Fn"{inner sep=.8ex}, "\shortmid"{marking}, from=1-9, to=1-10]
	\arrow["{{\sigma_z}}"{inner sep=.8ex}, "\shortmid"{marking}, from=1-10, to=1-11]
	\arrow[equals, from=1-10, to=2-10]
	\arrow[equals, from=1-11, to=2-11]
	\arrow["Fm"'{inner sep=.8ex}, "\shortmid"{marking}, from=2-1, to=2-2]
	\arrow[equals, from=2-1, to=3-2]
	\arrow["Fn"'{inner sep=.8ex}, "\shortmid"{marking}, from=2-2, to=2-3]
	\arrow[equals, from=2-2, to=3-3]
	\arrow["{{F(\mathrm{id}_z)}}"'{inner sep=.8ex}, "\shortmid"{marking}, from=2-3, to=2-4]
	\arrow["{{\tau_z}}"'{inner sep=.8ex}, "\shortmid"{marking}, from=2-4, to=2-5]
	\arrow[equals, from=2-4, to=3-4]
	\arrow[""{name=1, anchor=center, inner sep=0}, equals, from=2-5, to=3-5]
	\arrow[""{name=2, anchor=center, inner sep=0}, "{{F(m\odot n)}}"'{inner sep=.8ex}, "\shortmid"{marking}, from=2-8, to=2-10]
	\arrow[""{name=3, anchor=center, inner sep=0}, equals, from=2-8, to=3-7]
	\arrow[""{name=4, anchor=center, inner sep=0}, "{{\sigma_z}}"'{inner sep=.8ex}, "\shortmid"{marking}, from=2-10, to=2-11]
	\arrow[equals, from=2-10, to=3-9]
	\arrow[equals, from=2-11, to=3-11]
	\arrow["Fm"'{inner sep=.8ex}, "\shortmid"{marking}, from=3-2, to=3-3]
	\arrow[equals, from=3-2, to=4-2]
	\arrow[""{name=5, anchor=center, inner sep=0}, "Fn"'{inner sep=.8ex}, "\shortmid"{marking}, from=3-3, to=3-4]
	\arrow["{{\tau_z}}"'{inner sep=.8ex}, "\shortmid"{marking}, from=3-4, to=3-5]
	\arrow[equals, from=3-4, to=4-4]
	\arrow[equals, from=3-5, to=4-5]
	\arrow["{{F(m\odot n)}}"'{inner sep=.8ex}, "\shortmid"{marking}, from=3-7, to=3-9]
	\arrow[equals, from=3-7, to=4-8]
	\arrow["{{F(\mathrm{id}_z)}}"'{inner sep=.8ex}, "\shortmid"{marking}, from=3-9, to=3-10]
	\arrow["{{\tau_z}}"'{inner sep=.8ex}, "\shortmid"{marking}, from=3-10, to=3-11]
	\arrow[equals, from=3-10, to=4-10]
	\arrow[equals, from=3-11, to=4-11]
	\arrow[""{name=6, anchor=center, inner sep=0}, "{{F(m\odot n)}}"'{inner sep=.8ex}, "\shortmid"{marking}, from=4-2, to=4-4]
	\arrow["{{\tau_z}}"', from=4-4, to=4-5]
	\arrow[""{name=7, anchor=center, inner sep=0}, "{{F(m\odot n)}}"'{inner sep=.8ex}, "\shortmid"{marking}, from=4-8, to=4-10]
	\arrow["{{\tau_z}}"', from=4-10, to=4-11]
	\arrow["{{\mu_z}}"{description}, draw=none, from=0, to=2-4]
	\arrow["{{F_{m,n}}}"{description}, draw=none, from=1-9, to=2]
	\arrow["{{F_{n,\mathrm{id}_z}}}"{description}, draw=none, from=2-3, to=5]
	\arrow[between={0.3}{0.6}, squiggly, from=1, to=3]
	\arrow["{{\mu_z}}"{description, pos=0.6}, draw=none, from=4, to=3-10]
	\arrow["{{F_{m,n}}}"{description}, draw=none, from=3-3, to=6]
	\arrow["{{F_{m\odot n,\mathrm{id}_z}}}"{description}, draw=none, from=3-9, to=7]
\end{tikzcd}\]
  For the right action, we rearrange the pasting 
  for ``$\mu$, then $M_\tau^r$'', apply associativity of laxators, apply the modulation 
  axiom for $\mu$, and then reverse the process: 
\[\begin{tikzcd}[column sep=small]
	& Fx & Fy & Gy & Gz && Fx & Fy & Gy & Gz & \\
	Fx & Fy & Fy & Gy & Gz && Fx & Fy & Fz & Gz \\
	& Fx & Fy & Gy & Gz && Fx && Fz & Gz \\
	& Fx & Fy & Fz & Gz && Fx && Fz & Fz & Gz \\
	& Fx && Fz & Gz && Fx && Fz & Gz \\
	\\
	& Fx & Fy & Gy & Gz &&& Fx & Fy & Gy & Gz \\
	Fx & Fy & Fy & Gy & Gz &&& Fx & Fy & Fz & Gz \\
	Fx & Fy & Fy & Fz & Gz && Fx & Fy & Fz & Fz & Gz \\
	& Fx & Fy & Fz & Gz &&& Fx & Fy & Fz & Gz \\
	& Fx && Fz & Gz &&& Fx && Fz & Gz
	\arrow["Fm"{inner sep=.8ex}, "\shortmid"{marking}, from=1-2, to=1-3]
	\arrow[equals, from=1-2, to=2-1]
	\arrow[""{name=0, anchor=center, inner sep=0}, "{{{{\sigma_y}}}}"{inner sep=.8ex}, "\shortmid"{marking}, from=1-3, to=1-4]
	\arrow[equals, from=1-3, to=2-2]
	\arrow["Gn"{inner sep=.8ex}, "\shortmid"{marking}, from=1-4, to=1-5]
	\arrow[equals, from=1-4, to=2-4]
	\arrow[equals, from=1-5, to=2-5]
	\arrow["Fm"{inner sep=.8ex}, "\shortmid"{marking}, from=1-7, to=1-8]
	\arrow[equals, from=1-7, to=2-7]
	\arrow["{{{{\sigma_y}}}}"{inner sep=.8ex}, "\shortmid"{marking}, from=1-8, to=1-9]
	\arrow[equals, from=1-8, to=2-8]
	\arrow["Gn"{inner sep=.8ex}, "\shortmid"{marking}, from=1-9, to=1-10]
	\arrow["{{{{\sigma_n}}}}"{description}, draw=none, from=1-9, to=2-9]
	\arrow[equals, from=1-10, to=2-10]
	\arrow["Fm"'{inner sep=.8ex}, "\shortmid"{marking}, from=2-1, to=2-2]
	\arrow[equals, from=2-1, to=3-2]
	\arrow["{{{{F(\mathrm{id}_y)}}}}"'{inner sep=.8ex}, "\shortmid"{marking}, from=2-2, to=2-3]
	\arrow["{{{{\tau_y}}}}"'{inner sep=.8ex}, "\shortmid"{marking}, from=2-3, to=2-4]
	\arrow[equals, from=2-3, to=3-3]
	\arrow["Gn"'{inner sep=.8ex}, "\shortmid"{marking}, from=2-4, to=2-5]
	\arrow[equals, from=2-4, to=3-4]
	\arrow[equals, from=2-5, to=3-5]
	\arrow["Fm"'{inner sep=.8ex}, "\shortmid"{marking}, from=2-7, to=2-8]
	\arrow[equals, from=2-7, to=3-7]
	\arrow["Fn"'{inner sep=.8ex}, "\shortmid"{marking}, from=2-8, to=2-9]
	\arrow["{{\sigma_z}}"'{inner sep=.8ex}, "\shortmid"{marking}, from=2-9, to=2-10]
	\arrow[equals, from=2-9, to=3-9]
	\arrow[equals, from=2-10, to=3-10]
	\arrow[""{name=1, anchor=center, inner sep=0}, "Fm"'{inner sep=.8ex}, "\shortmid"{marking}, from=3-2, to=3-3]
	\arrow[equals, from=3-2, to=4-2]
	\arrow["{{{{\tau_y}}}}"'{inner sep=.8ex}, "\shortmid"{marking}, from=3-3, to=3-4]
	\arrow[equals, from=3-3, to=4-3]
	\arrow["Gn"'{inner sep=.8ex}, "\shortmid"{marking}, from=3-4, to=3-5]
	\arrow["{{{{\tau_n}}}}"{description}, draw=none, from=3-4, to=4-4]
	\arrow[equals, from=3-5, to=4-5]
	\arrow[""{name=2, anchor=center, inner sep=0}, "{{F(m\odot n)}}"'{inner sep=.8ex}, "\shortmid"{marking}, from=3-7, to=3-9]
	\arrow[equals, from=3-7, to=4-7]
	\arrow[""{name=3, anchor=center, inner sep=0}, "{{\sigma_z}}"'{inner sep=.8ex}, "\shortmid"{marking}, from=3-9, to=3-10]
	\arrow[equals, from=3-9, to=4-9]
	\arrow[equals, from=3-10, to=4-11]
	\arrow["Fm"'{inner sep=.8ex}, "\shortmid"{marking}, from=4-2, to=4-3]
	\arrow[equals, from=4-2, to=5-2]
	\arrow["Fn"'{inner sep=.8ex}, "\shortmid"{marking}, from=4-3, to=4-4]
	\arrow["{{{{\tau_z}}}}"'{inner sep=.8ex}, "\shortmid"{marking}, from=4-4, to=4-5]
	\arrow[equals, from=4-4, to=5-4]
	\arrow[equals, from=4-5, to=5-5]
	\arrow["{{F(m\odot n)}}"'{inner sep=.8ex}, "\shortmid"{marking}, from=4-7, to=4-9]
	\arrow[equals, from=4-7, to=5-7]
	\arrow["{{F(\mathrm{id}_z)}}"'{inner sep=.8ex}, "\shortmid"{marking}, from=4-9, to=4-10]
	\arrow["{{\tau_z}}"'{inner sep=.8ex}, "\shortmid"{marking}, from=4-10, to=4-11]
	\arrow[curve={height=-6pt}, equals, from=4-10, to=5-9]
	\arrow[equals, from=4-11, to=5-10]
	\arrow[""{name=4, anchor=center, inner sep=0}, "{{{{F(m\odot n)}}}}"'{inner sep=.8ex}, "\shortmid"{marking}, from=5-2, to=5-4]
	\arrow["{{{{\tau_z}}}}"'{inner sep=.8ex}, "\shortmid"{marking}, from=5-4, to=5-5]
	\arrow[between={0.3}{0.6}, squiggly, from=5-4, to=7-4]
	\arrow[""{name=5, anchor=center, inner sep=0}, "{{F(m\odot n)}}"'{inner sep=.8ex}, "\shortmid"{marking}, from=5-7, to=5-9]
	\arrow["{{\tau_z}}"'{inner sep=.8ex}, "\shortmid"{marking}, from=5-9, to=5-10]
	\arrow["Fm"{inner sep=.8ex}, "\shortmid"{marking}, from=7-2, to=7-3]
	\arrow[equals, from=7-2, to=8-1]
	\arrow[""{name=6, anchor=center, inner sep=0}, "{{{{\sigma_y}}}}"{inner sep=.8ex}, "\shortmid"{marking}, from=7-3, to=7-4]
	\arrow[equals, from=7-3, to=8-2]
	\arrow["Gn"{inner sep=.8ex}, "\shortmid"{marking}, from=7-4, to=7-5]
	\arrow[equals, from=7-4, to=8-4]
	\arrow[equals, from=7-5, to=8-5]
	\arrow["Fm"{inner sep=.8ex}, "\shortmid"{marking}, from=7-8, to=7-9]
	\arrow[equals, from=7-8, to=8-8]
	\arrow[between={0.3}{0.6}, squiggly, from=7-9, to=5-9]
	\arrow["{{{{\sigma_y}}}}"{inner sep=.8ex}, "\shortmid"{marking}, from=7-9, to=7-10]
	\arrow[equals, from=7-9, to=8-9]
	\arrow["Gn"{inner sep=.8ex}, "\shortmid"{marking}, from=7-10, to=7-11]
	\arrow["{{\sigma_n}}"{description}, draw=none, from=7-10, to=8-10]
	\arrow[equals, from=7-11, to=8-11]
	\arrow["Fm"'{inner sep=.8ex}, "\shortmid"{marking}, from=8-1, to=8-2]
	\arrow[equals, from=8-1, to=9-1]
	\arrow["{{{{F(\mathrm{id}_y)}}}}"'{inner sep=.8ex}, "\shortmid"{marking}, from=8-2, to=8-3]
	\arrow[equals, from=8-2, to=9-2]
	\arrow["{{{{\tau_y}}}}"'{inner sep=.8ex}, "\shortmid"{marking}, from=8-3, to=8-4]
	\arrow[equals, from=8-3, to=9-3]
	\arrow["Gn"'{inner sep=.8ex}, "\shortmid"{marking}, from=8-4, to=8-5]
	\arrow["{{\tau_n}}"{description}, draw=none, from=8-4, to=9-4]
	\arrow[equals, from=8-5, to=9-5]
	\arrow["Fm"'{inner sep=.8ex}, "\shortmid"{marking}, from=8-8, to=8-9]
	\arrow[equals, from=8-8, to=9-7]
	\arrow["Fn"'{inner sep=.8ex}, "\shortmid"{marking}, from=8-9, to=8-10]
	\arrow[equals, from=8-9, to=9-8]
	\arrow[""{name=7, anchor=center, inner sep=0}, "{{\sigma_z}}"'{inner sep=.8ex}, "\shortmid"{marking}, from=8-10, to=8-11]
	\arrow[equals, from=8-10, to=9-9]
	\arrow[equals, from=8-11, to=9-11]
	\arrow["Fm"'{inner sep=.8ex}, "\shortmid"{marking}, from=9-1, to=9-2]
	\arrow[equals, from=9-1, to=10-2]
	\arrow["{{F(\mathrm{id}_y)}}"'{inner sep=.8ex}, "\shortmid"{marking}, from=9-2, to=9-3]
	\arrow[curve={height=6pt}, equals, from=9-2, to=10-3]
	\arrow["Fn"'{inner sep=.8ex}, "\shortmid"{marking}, from=9-3, to=9-4]
	\arrow["{{\tau_z}}"'{inner sep=.8ex}, "\shortmid"{marking}, from=9-4, to=9-5]
	\arrow[equals, from=9-4, to=10-4]
	\arrow[between={0.3}{0.7}, squiggly, from=9-5, to=9-7]
	\arrow[equals, from=9-5, to=10-5]
	\arrow["Fm"'{inner sep=.8ex}, "\shortmid"{marking}, from=9-7, to=9-8]
	\arrow[equals, from=9-7, to=10-8]
	\arrow["{{F(n)}}"'{inner sep=.8ex}, "\shortmid"{marking}, from=9-8, to=9-9]
	\arrow[curve={height=6pt}, equals, from=9-8, to=10-9]
	\arrow["{{F(\mathrm{id}_z)}}"'{inner sep=.8ex}, "\shortmid"{marking}, from=9-9, to=9-10]
	\arrow["{{\tau_z}}"'{inner sep=.8ex}, "\shortmid"{marking}, from=9-10, to=9-11]
	\arrow[equals, from=9-10, to=10-10]
	\arrow[equals, from=9-11, to=10-11]
	\arrow["Fm"'{inner sep=.8ex}, "\shortmid"{marking}, from=10-2, to=10-3]
	\arrow[equals, from=10-2, to=11-2]
	\arrow[""{name=8, anchor=center, inner sep=0}, "Fn"'{inner sep=.8ex}, "\shortmid"{marking}, from=10-3, to=10-4]
	\arrow["{{{{\tau_z}}}}"'{inner sep=.8ex}, "\shortmid"{marking}, from=10-4, to=10-5]
	\arrow[equals, from=10-4, to=11-4]
	\arrow[equals, from=10-5, to=11-5]
	\arrow["Fm"'{inner sep=.8ex}, "\shortmid"{marking}, from=10-8, to=10-9]
	\arrow[equals, from=10-8, to=11-8]
	\arrow[""{name=9, anchor=center, inner sep=0}, "Fn"'{inner sep=.8ex}, "\shortmid"{marking}, from=10-9, to=10-10]
	\arrow["{{{{\tau_z}}}}"'{inner sep=.8ex}, "\shortmid"{marking}, from=10-10, to=10-11]
	\arrow[equals, from=10-10, to=11-10]
	\arrow[equals, from=10-11, to=11-11]
	\arrow[""{name=10, anchor=center, inner sep=0}, "{{{{F(m\odot n)}}}}"'{inner sep=.8ex}, "\shortmid"{marking}, from=11-2, to=11-4]
	\arrow["{{{{\tau_z}}}}"'{inner sep=.8ex}, "\shortmid"{marking}, from=11-4, to=11-5]
	\arrow[""{name=11, anchor=center, inner sep=0}, "{{{{F(m\odot n)}}}}"'{inner sep=.8ex}, "\shortmid"{marking}, from=11-8, to=11-10]
	\arrow["{{{{\tau_z}}}}"'{inner sep=.8ex}, "\shortmid"{marking}, from=11-10, to=11-11]
	\arrow["{{{{\mu_y}}}}"{description}, draw=none, from=0, to=2-3]
	\arrow["{{{{F_{m,\mathrm{id}_y}}}}}"{description}, draw=none, from=2-2, to=1]
	\arrow["{{F_{m,n}}}"{description}, draw=none, from=2-8, to=2]
	\arrow["{{\mu_z}}"{description, pos=0.6}, draw=none, from=3, to=4-10]
	\arrow["{{{{F_{m,n}}}}}"{description}, draw=none, from=4-3, to=4]
	\arrow["{{F_{m\odot n,\mathrm{id}_z}}}"{description}, draw=none, from=4-9, to=5]
	\arrow["{{{{\mu_y}}}}"{description}, draw=none, from=6, to=8-3]
	\arrow["{{\mu_z}}"{description, pos=0.6}, draw=none, from=7, to=9-10]
	\arrow["{{F_{\mathrm{id}_y,n}}}"{description}, draw=none, from=9-3, to=8]
	\arrow["{{F_{n,\mathrm{id}_z}}}"{description}, draw=none, from=9-9, to=9]
	\arrow["{{{{F_{m,n}}}}}"{description}, draw=none, from=10-3, to=10]
	\arrow["{{{{F_{m,n}}}}}"{description}, draw=none, from=10-9, to=11]
\end{tikzcd}\]

  \item We must check that the modulation hexagon commutes. Its opposite sides comprise the upper-left
  and upper-right pastings in the diagram below. Starting from the upper-left, we move down 
  by applying the respect of $\mu$ for the right actions of $M_\sigma$ and $M_\tau,$ 
  then to the center by applying unitality of $F$'s laxators. Similarly, starting from the 
  upper-right, we move down by applying the respect of $\mu$ for the left actions of $M_\sigma$ and $M_\tau,$ 
  then to the center by applying unitality of $F$'s laxators. 
\[\begin{tikzcd}[column sep=small]
	& Fx & Gx & Gy &&&&& Fx & Gx & Gy \\
	Fx & Fx & Gx &&&&&& Fx & Fy & Gy \\
	Fx & Fx & Gx & Gy &&&& Fx & Fy & Fy & Gy \\
	Fx & Fx & Fy & Gy &&&& Fx & Fy & Fy & Gy \\
	Fx && Fy & Fy & Fx & Gx & Gy & Fx && Fy & Fy \\
	&&&& Fx & Fy & Gy \\
	& Fx & Gx & Gy & Fx & Fy & Gy && Fx & Gx & Gy \\
	Fx & Fx & Gx & Gy &&&&& Fx & Fy & Gy \\
	Fx & Fx & Fy & Gy &&&& Fx & Fy & Fy & Gy \\
	Fx && Fy & Gy &&&& Fx && Fy & Gy \\
	Fx && Fy & Fy &&&& Fx && Fy & Fy
	\arrow["{\sigma_x}"{inner sep=.8ex}, "\shortmid"{marking}, from=1-2, to=1-3]
	\arrow[equals, from=1-2, to=2-1]
	\arrow[equals, from=1-2, to=2-2]
	\arrow["Gm"{inner sep=.8ex}, "\shortmid"{marking}, from=1-3, to=1-4]
	\arrow[equals, from=1-3, to=2-3]
	\arrow[equals, from=1-4, to=3-4]
	\arrow["{\sigma_x}"{inner sep=.8ex}, "\shortmid"{marking}, from=1-9, to=1-10]
	\arrow[equals, from=1-9, to=2-9]
	\arrow["Gm"{inner sep=.8ex}, "\shortmid"{marking}, from=1-10, to=1-11]
	\arrow["{\sigma_m}"{description}, draw=none, from=1-10, to=2-10]
	\arrow[equals, from=1-11, to=2-11]
	\arrow[""{name=0, anchor=center, inner sep=0}, "{F(\mathrm{id}_x)}"'{inner sep=.8ex}, "\shortmid"{marking}, from=2-1, to=2-2]
	\arrow[equals, from=2-1, to=3-1]
	\arrow["{\sigma_x}"'{inner sep=.8ex}, "\shortmid"{marking}, from=2-2, to=2-3]
	\arrow["{\mu_{\mathrm{id}_x}}"{description}, draw=none, from=2-2, to=3-2]
	\arrow[equals, from=2-3, to=3-3]
	\arrow["Fm"'{inner sep=.8ex}, "\shortmid"{marking}, from=2-9, to=2-10]
	\arrow[equals, from=2-9, to=3-8]
	\arrow["{\sigma_y}"'{inner sep=.8ex}, "\shortmid"{marking}, from=2-10, to=2-11]
	\arrow[equals, from=2-10, to=3-9]
	\arrow[equals, from=2-10, to=3-10]
	\arrow[equals, from=2-11, to=3-11]
	\arrow["{F(\mathrm{id}_x)}"'{inner sep=.8ex}, "\shortmid"{marking}, from=3-1, to=3-2]
	\arrow[equals, from=3-1, to=4-1]
	\arrow["{\tau_x}"'{inner sep=.8ex}, "\shortmid"{marking}, from=3-2, to=3-3]
	\arrow[equals, from=3-2, to=4-2]
	\arrow["Gm"'{inner sep=.8ex}, "\shortmid"{marking}, from=3-3, to=3-4]
	\arrow["{\tau_m}"{description}, draw=none, from=3-3, to=4-3]
	\arrow[equals, from=3-4, to=4-4]
	\arrow["Fm"'{inner sep=.8ex}, "\shortmid"{marking}, from=3-8, to=3-9]
	\arrow[equals, from=3-8, to=4-8]
	\arrow[""{name=1, anchor=center, inner sep=0}, "{F(\mathrm{id}_y)}"'{inner sep=.8ex}, "\shortmid"{marking}, from=3-9, to=3-10]
	\arrow[equals, from=3-9, to=4-9]
	\arrow["{\sigma_y}"'{inner sep=.8ex}, "\shortmid"{marking}, from=3-10, to=3-11]
	\arrow["{\mu_{\mathrm{id}_y}}"{description}, draw=none, from=3-10, to=4-10]
	\arrow[equals, from=3-11, to=4-11]
	\arrow["{F(\mathrm{id}_x)}"'{inner sep=.8ex}, "\shortmid"{marking}, from=4-1, to=4-2]
	\arrow[equals, from=4-1, to=5-1]
	\arrow["Fm"'{inner sep=.8ex}, "\shortmid"{marking}, from=4-2, to=4-3]
	\arrow["{\tau_y}"'{inner sep=.8ex}, "\shortmid"{marking}, from=4-3, to=4-4]
	\arrow[equals, from=4-3, to=5-3]
	\arrow[equals, from=4-4, to=5-4]
	\arrow["Fm"'{inner sep=.8ex}, "\shortmid"{marking}, from=4-8, to=4-9]
	\arrow[equals, from=4-8, to=5-8]
	\arrow["{F(\mathrm{id}_y)}"'{inner sep=.8ex}, "\shortmid"{marking}, from=4-9, to=4-10]
	\arrow["{\tau_y}"'{inner sep=.8ex}, "\shortmid"{marking}, from=4-10, to=4-11]
	\arrow[equals, from=4-10, to=5-10]
	\arrow[equals, from=4-11, to=5-11]
	\arrow[""{name=2, anchor=center, inner sep=0}, "Fm"'{inner sep=.8ex}, "\shortmid"{marking}, from=5-1, to=5-3]
	\arrow["{\tau_y}"'{inner sep=.8ex}, "\shortmid"{marking}, from=5-3, to=5-4]
	\arrow[squiggly, from=5-3, to=7-3]
	\arrow["{\sigma_x}"{inner sep=.8ex}, "\shortmid"{marking}, from=5-5, to=5-6]
	\arrow[equals, from=5-5, to=6-5]
	\arrow["Gm"{inner sep=.8ex}, "\shortmid"{marking}, from=5-6, to=5-7]
	\arrow["{\sigma_m}"{description}, draw=none, from=5-6, to=6-6]
	\arrow[equals, from=5-7, to=6-7]
	\arrow[""{name=3, anchor=center, inner sep=0}, "Fm"'{inner sep=.8ex}, "\shortmid"{marking}, from=5-8, to=5-10]
	\arrow["{\tau_y}"'{inner sep=.8ex}, "\shortmid"{marking}, from=5-10, to=5-11]
	\arrow[squiggly, from=5-10, to=7-10]
	\arrow["Fm"'{inner sep=.8ex}, "\shortmid"{marking}, from=6-5, to=6-6]
	\arrow[equals, from=6-5, to=7-5]
	\arrow["{\sigma_y}"'{inner sep=.8ex}, "\shortmid"{marking}, from=6-6, to=6-7]
	\arrow["{\mu_m}"{description}, draw=none, from=6-6, to=7-6]
	\arrow[equals, from=6-7, to=7-7]
	\arrow["{\sigma_x}"{inner sep=.8ex}, "\shortmid"{marking}, from=7-2, to=7-3]
	\arrow[equals, from=7-2, to=8-1]
	\arrow[equals, from=7-2, to=8-2]
	\arrow["Gm"{inner sep=.8ex}, "\shortmid"{marking}, from=7-3, to=7-4]
	\arrow[equals, from=7-3, to=8-3]
	\arrow[equals, from=7-4, to=8-4]
	\arrow["Fm"'{inner sep=.8ex}, "\shortmid"{marking}, from=7-5, to=7-6]
	\arrow["{\tau_y}"'{inner sep=.8ex}, "\shortmid"{marking}, from=7-6, to=7-7]
	\arrow["{\sigma_x}"{inner sep=.8ex}, "\shortmid"{marking}, from=7-9, to=7-10]
	\arrow[equals, from=7-9, to=8-9]
	\arrow["Gm"{inner sep=.8ex}, "\shortmid"{marking}, from=7-10, to=7-11]
	\arrow["{\sigma_m}"{description}, draw=none, from=7-10, to=8-10]
	\arrow[equals, from=7-11, to=8-11]
	\arrow[""{name=4, anchor=center, inner sep=0}, "{F(\mathrm{id}_x)}"'{inner sep=.8ex}, "\shortmid"{marking}, from=8-1, to=8-2]
	\arrow[equals, from=8-1, to=9-1]
	\arrow["{\sigma_x}"'{inner sep=.8ex}, "\shortmid"{marking}, from=8-2, to=8-3]
	\arrow[equals, from=8-2, to=9-2]
	\arrow["Gm"'{inner sep=.8ex}, "\shortmid"{marking}, from=8-3, to=8-4]
	\arrow["{\sigma_m}"{description}, draw=none, from=8-3, to=9-3]
	\arrow[squiggly, from=8-4, to=6-5]
	\arrow[equals, from=8-4, to=9-4]
	\arrow[squiggly, from=8-9, to=6-7]
	\arrow["Fm"'{inner sep=.8ex}, "\shortmid"{marking}, from=8-9, to=8-10]
	\arrow[equals, from=8-9, to=9-8]
	\arrow["{\sigma_y}"'{inner sep=.8ex}, "\shortmid"{marking}, from=8-10, to=8-11]
	\arrow[equals, from=8-10, to=9-9]
	\arrow[equals, from=8-10, to=9-10]
	\arrow[equals, from=8-11, to=9-11]
	\arrow["{F(\mathrm{id}_x)}"'{inner sep=.8ex}, "\shortmid"{marking}, from=9-1, to=9-2]
	\arrow[equals, from=9-1, to=10-1]
	\arrow["Fm"'{inner sep=.8ex}, "\shortmid"{marking}, from=9-2, to=9-3]
	\arrow["{\sigma_y}"'{inner sep=.8ex}, "\shortmid"{marking}, from=9-3, to=9-4]
	\arrow[equals, from=9-3, to=10-3]
	\arrow[equals, from=9-4, to=10-4]
	\arrow["Fm"'{inner sep=.8ex}, "\shortmid"{marking}, from=9-8, to=9-9]
	\arrow[equals, from=9-8, to=10-8]
	\arrow[""{name=5, anchor=center, inner sep=0}, "{F(\mathrm{id}_y)}"'{inner sep=.8ex}, "\shortmid"{marking}, from=9-9, to=9-10]
	\arrow["{\sigma_y}"'{inner sep=.8ex}, "\shortmid"{marking}, from=9-10, to=9-11]
	\arrow[equals, from=9-10, to=10-10]
	\arrow[equals, from=9-11, to=10-11]
	\arrow[""{name=6, anchor=center, inner sep=0}, "Fm"'{inner sep=.8ex}, "\shortmid"{marking}, from=10-1, to=10-3]
	\arrow[equals, from=10-1, to=11-1]
	\arrow["{\sigma_y}"'{inner sep=.8ex}, "\shortmid"{marking}, from=10-3, to=10-4]
	\arrow["{\mu_m}"{description}, draw=none, from=10-3, to=11-3]
	\arrow[equals, from=10-4, to=11-4]
	\arrow[""{name=7, anchor=center, inner sep=0}, "Fm"'{inner sep=.8ex}, "\shortmid"{marking}, from=10-8, to=10-10]
	\arrow[equals, from=10-8, to=11-8]
	\arrow["{\sigma_y}"'{inner sep=.8ex}, "\shortmid"{marking}, from=10-10, to=10-11]
	\arrow["{\mu_m}"{description}, draw=none, from=10-10, to=11-10]
	\arrow[equals, from=10-11, to=11-11]
	\arrow["Fm"'{inner sep=.8ex}, "\shortmid"{marking}, from=11-1, to=11-3]
	\arrow["{\tau_y}"'{inner sep=.8ex}, "\shortmid"{marking}, from=11-3, to=11-4]
	\arrow["Fm"'{inner sep=.8ex}, "\shortmid"{marking}, from=11-8, to=11-10]
	\arrow["{\tau_y}"'{inner sep=.8ex}, "\shortmid"{marking}, from=11-10, to=11-11]
	\arrow["{F_x}"{description}, draw=none, from=1-2, to=0]
	\arrow["{F_y}"{description}, draw=none, from=2-10, to=1]
	\arrow["{F_{\mathrm{id}_x,m}}"{description}, draw=none, from=4-2, to=2]
	\arrow["{F_{m,\mathrm{id}_y}}"{description}, draw=none, from=4-9, to=3]
	\arrow["{F_x}"{description}, draw=none, from=7-2, to=4]
	\arrow["{F_y}"{description}, draw=none, from=8-10, to=5]
	\arrow["{F_{\mathrm{id}_x,m}}"{description}, draw=none, from=9-2, to=6]
	\arrow["{F_{m,\mathrm{id}_y}}"{description}, draw=none, from=9-9, to=7]
\end{tikzcd}\]
  \end{enumerate}
    
\end{proof}

\begin{theorem}[Embedding loose transformations into modules]\label{thm:loose-transformations-into-modules}
(Cockett, Koslowski, Seely, and Wood)
For any two lax functors $F,G:\dbl{D}\to \dbl{E},$ we have a fully faithful embedding 
of the category $\twocat{LTrans}(F,G)$ of loose transformations and \emph{modulations} 
into the category $\twocat{Mod}(F,G)$ of modules and modulations. 

Suppose $F$ is pseudo. Then this embedding is 
an equivalence onto the category of modules with invertible left actions; 
if furthermore $\dbl{D}=\dbl{T}C$ is tightly discrete, then the embedding is essentially 
surjective. If $F$ is even strict, 
then the embedding is an equivalence onto the category of modules with identity left actions.
\end{theorem}
\begin{proof}
The embedding will of course map $\tau\mapsto M_\tau$ and $\mu\mapsto \bar\mu$ via 
\cref{cons:loose-transformation-to-module} above, which succeeds by \cref{lemma:constructions-work}.

We must still check functoriality of the mapping $\mu\mapsto \bar \mu.$ Preservation of the identity modulation
follows from unitality of $F$'s laxators. Given a composite 
$\sigma\xto{\mu}\tau\xto{\nu}\upsilon,$ we check that $\overline{\nu\circ \mu}=\bar\nu\circ \bar\mu.$
This follows from sliding and naturality of laxators for $F$: 
\[\begin{tikzcd}
	Fx & Fy & Gy && Fx & Fy & Gy && \\
	Fx & Fy & Fy & Gy & Fx & Fy & Fy & Gy \\
	Fx & Fy & Gy && Fx & Fy & Fy & Fy & Gy \\
	Fx & Fy & Fy & Gy & Fx & Fy & Fy & Gy \\
	Fx & Fy & Gy && Fx & Fy & Gy
	\arrow["Fm"{inner sep=.8ex}, "\shortmid"{marking}, from=1-1, to=1-2]
	\arrow[equals, from=1-1, to=2-1]
	\arrow[""{name=0, anchor=center, inner sep=0}, "{\sigma_y}"{inner sep=.8ex}, "\shortmid"{marking}, from=1-2, to=1-3]
	\arrow[equals, from=1-2, to=2-2]
	\arrow[equals, from=1-3, to=2-4]
	\arrow["Fm"{inner sep=.8ex}, "\shortmid"{marking}, from=1-5, to=1-6]
	\arrow[equals, from=1-5, to=2-5]
	\arrow[""{name=1, anchor=center, inner sep=0}, "{\sigma_y}"{inner sep=.8ex}, "\shortmid"{marking}, from=1-6, to=1-7]
	\arrow[equals, from=1-6, to=2-6]
	\arrow[equals, from=1-7, to=2-8]
	\arrow["Fm"'{inner sep=.8ex}, "\shortmid"{marking}, from=2-1, to=2-2]
	\arrow[equals, from=2-1, to=3-1]
	\arrow["{F(\mathrm{id}_y)}"'{inner sep=.8ex}, "\shortmid"{marking}, from=2-2, to=2-3]
	\arrow["{\tau_y}"'{inner sep=.8ex}, "\shortmid"{marking}, from=2-3, to=2-4]
	\arrow[curve={height=-6pt}, equals, from=2-3, to=3-2]
	\arrow[equals, from=2-4, to=3-3]
	\arrow["Fm"'{inner sep=.8ex}, "\shortmid"{marking}, from=2-5, to=2-6]
	\arrow[equals, from=2-5, to=3-5]
	\arrow["{F(\mathrm{id}_y)}"'{inner sep=.8ex}, "\shortmid"{marking}, from=2-6, to=2-7]
	\arrow[equals, from=2-6, to=3-6]
	\arrow[""{name=2, anchor=center, inner sep=0}, "{\tau_y}"'{inner sep=.8ex}, "\shortmid"{marking}, from=2-7, to=2-8]
	\arrow[equals, from=2-7, to=3-7]
	\arrow[equals, from=2-8, to=3-9]
	\arrow[""{name=3, anchor=center, inner sep=0}, "Fm"'{inner sep=.8ex}, "\shortmid"{marking}, from=3-1, to=3-2]
	\arrow[equals, from=3-1, to=4-1]
	\arrow[""{name=4, anchor=center, inner sep=0}, "{\tau_y}"'{inner sep=.8ex}, "\shortmid"{marking}, from=3-2, to=3-3]
	\arrow[equals, from=3-2, to=4-2]
	\arrow[between={0.1}{0.9}, squiggly, from=3-3, to=3-5]
	\arrow[equals, from=3-3, to=4-4]
	\arrow["Fm"'{inner sep=.8ex}, "\shortmid"{marking}, from=3-5, to=3-6]
	\arrow[equals, from=3-5, to=4-5]
	\arrow["{F(\mathrm{id}_y)}"'{inner sep=.8ex}, "\shortmid"{marking}, from=3-6, to=3-7]
	\arrow[equals, from=3-6, to=4-6]
	\arrow["{F(\mathrm{id}_y)}"'{inner sep=.8ex}, "\shortmid"{marking}, from=3-7, to=3-8]
	\arrow["{\upsilon_y}"'{inner sep=.8ex}, "\shortmid"{marking}, from=3-8, to=3-9]
	\arrow[curve={height=-6pt}, equals, from=3-8, to=4-7]
	\arrow[equals, from=3-9, to=4-8]
	\arrow["Fm"'{inner sep=.8ex}, "\shortmid"{marking}, from=4-1, to=4-2]
	\arrow[equals, from=4-1, to=5-1]
	\arrow["{F(\mathrm{id}_y)}"'{inner sep=.8ex}, "\shortmid"{marking}, from=4-2, to=4-3]
	\arrow["{\upsilon_y}"'{inner sep=.8ex}, "\shortmid"{marking}, from=4-3, to=4-4]
	\arrow[curve={height=-6pt}, equals, from=4-3, to=5-2]
	\arrow[equals, from=4-4, to=5-3]
	\arrow["Fm"'{inner sep=.8ex}, "\shortmid"{marking}, from=4-5, to=4-6]
	\arrow[equals, from=4-5, to=5-5]
	\arrow[""{name=5, anchor=center, inner sep=0}, "{F(\mathrm{id}_y)}"'{inner sep=.8ex}, "\shortmid"{marking}, from=4-6, to=4-7]
	\arrow["{\upsilon_y}"'{inner sep=.8ex}, "\shortmid"{marking}, from=4-7, to=4-8]
	\arrow[curve={height=-6pt}, equals, from=4-7, to=5-6]
	\arrow[equals, from=4-8, to=5-7]
	\arrow[""{name=6, anchor=center, inner sep=0}, "Fm"'{inner sep=.8ex}, "\shortmid"{marking}, from=5-1, to=5-2]
	\arrow["{\tau_y}"'{inner sep=.8ex}, "\shortmid"{marking}, from=5-2, to=5-3]
	\arrow[""{name=7, anchor=center, inner sep=0}, "Fm"'{inner sep=.8ex}, "\shortmid"{marking}, from=5-5, to=5-6]
	\arrow["{\tau_y}"'{inner sep=.8ex}, "\shortmid"{marking}, from=5-6, to=5-7]
	\arrow["{\mu_y}"{description}, draw=none, from=0, to=2-3]
	\arrow["{\mu_y}"{description}, draw=none, from=1, to=2-7]
	\arrow["{F_{m,\mathrm{id}_y}}"{description}, draw=none, from=2-2, to=3]
	\arrow["{\nu_y}"{description, pos=0.6}, draw=none, from=2, to=3-8]
	\arrow["{\nu_y}"{description, pos=0.6}, draw=none, from=4, to=4-3]
	\arrow["{F_{\mathrm{id}_y,\mathrm{id}_y}}"{description}, draw=none, from=3-7, to=5]
	\arrow["{F_{m,\mathrm{id}_y}}"{description}, draw=none, from=4-2, to=6]
	\arrow["{F_{m,\mathrm{id}_y}}"{description}, draw=none, from=4-6, to=7]
\end{tikzcd}\]

We can see that $\underline{(\bar \mu)}=\mu$ and $\overline{(\underline{\mu})}=\mu,$ 
respectively by unitality of $F$ and by the same axiom together with respect for $\mu$ for 
the left actions of $M_\sigma,M_\tau.$ This proves full faithfulness of the embedding.

Furthermore, supposing that $F$ is pseudo, then since for a loose transformation $\tau:F\to G,$ 
the associated module $M_\tau$ has as left actions the pasting of a laxator of $F$ with an identity
cell, our embedding lands in modules with invertible left actions. Conversely, given 
a module $M$ with invertible left actions, we have constructed the loose transformation 
$\tau^M$ with components $\tau^M_x=M(\id_x).$ Then $\tau^{M_\tau}_x=F(\id_x)\odot \tau_x,$ 
so that $\tau^{M_\tau}\cong \tau$ as long as $F$ has invertible unitors, which gives the result.
The case of $F$ strict goes the same way. If $\dbl{D}$ lacks nontrivial loose arrows, then the only
left actions for $M$ are indexed by identities, which are invertible by the invertibility of 
$F$'s unitors and unitality of $M^\ell.$ Thus in this case every module arises from a 
loose transformation.  
\end{proof}

\section{Local presentability of categories of models and their instances}\label{app:local-presentability}

In this appendix, we explicitly construct a finite limit theory whose models in $\Set$ 
coincide with the models of a given simple double theory $\dbl{D}$ in $\Span.$ 

\begin{remark}[Comparison with the construction of \cite{limits-in-double-1999}]\label{rem:connection-to-double-limits}
Our construction can be seen as a refinement of that in \cite[5.7 (A)]{limits-in-double-1999}, 
which gives a tightly discrete double category, limits over whose diagrams coincide with double 
limits over a given double category. Our construction is more complex because it is harder 
to get equivalent categories of diagrams than to get a functor between them commuting 
with limits. This phenomenon is familiar from ordinary category theory, where one can even take
the limit of a diagram by considering the free category on its underlying graph. 
That said, Grandis and Par\'e handle a more general codomain; we expect that 
the two arguments could be merged to give an analogous result 
to \cref{prop:cartesian-models-locally-presentable} for models of any double theory 
in any double category with pullbacks and tabulators.
\end{remark}

\begin{construction}[Flattening of a double theory]
  Given a simple double theory $\dbl{D},$ consider the following sketch $\varphi(\dbl{D}),$
  which we call the \emph{flattening} of $\dbl{D}.$

  It will aid the understanding of the construction to look ahead to our intent to prove that 
  models of this sketch are equivalent to
  models of $\dbl{D}$ in $\Span.$
  \begin{itemize}
  \item Objects of $\varphi(\dbl{D})$ come in four sorts:
  \begin{enumerate}
  \item Each object $x$ of $\dbl{D}_0$ provides an object $[x]$ of $\varphi(\dbl{D})$.
  \item Each loose morphism $m:x\proto y$ in $\dbl{D}$ provides an object $[m]$ of $\varphi(\dbl{D})$.
  \item Each composable pair $x\xproto{m}y\xproto{n}z$ in $\dbl{D}$ 
  provides an object $[m,n]$ of $\varphi(\dbl{D})$.
  \item Each composable triple $x\xproto{m}y\xproto{n}z\xproto{p}w$ in $\dbl{D}$ provides
  an object $[m,n,p]$ of $\varphi(\dbl{D})$.
  \end{enumerate}
  \item Morphisms in $\varphi(\dbl{D})$ are generated by the following eight classes:
  \begin{enumerate}
  \item Each tight morphism $f:x\to y$ in $\dbl{D}$ provides a morphism $[f]:[x]\to [y]$ in $\varphi(\dbl{D})$.
  \item Each loose morphism $m:x\proto y$ in $\dbl{D}$ provides a span $[x]\stackrel{s_m}{\leftarrow} [m]\xto{t_m} [y]$ in $\varphi(\dbl{D})$.
  \item Each composable pair $(m,n)$ in $\dbl{D}$ provides morphisms 
    \begin{itemize}
      \item $s_{m,n}:[m,n]\to [m]$
      \item $c_{m,n}:[m,n]\to [m\odot n]$
      \item $t_{m,n}:[m,n]\to [n].$ 
    \end{itemize} 
  \item Each cell $\alpha:m\to m'$ in $\dbl{D}$ provides a morphism $[\alpha]:[m]\to [n]$ in $\varphi(\dbl{D})$.
  \item Each loosely-composable pair $\alpha:m\to m',\beta:n\to n'$ of cells provides a morphism
  $[\alpha,\beta]: [m,n]\to [m',n']$ in $\varphi(\dbl{D})$.
  \item Each object $x$ in $\dbl{D}$ provides a morphism $u_x:[x]\to [\id_x]$ in $\varphi(\dbl{D})$.
  \item Each loose morphism $m:x\proto y$ in $\dbl{D}$ provides morphisms 
    \begin{itemize}
    \item $\ell_m:[m]\to [\id_x,m]$ 
    \item $r_m:[m]\to [m,\id_y].$
    \end{itemize} 
  \item Each composable triple $x\xproto{m}y\xproto{n}z\xproto{p}w$ in $\dbl{D}$ provides morphisms 
    \begin{itemize}
    \item $\ell_{m,n,p}:[m,n,p]\to [m\odot n,p]$
    \item $r_{m,n,p}:[m,n,p]\to [m,n\odot p]$
    \item $s_{m,n,p}:[m,n,p]\to [m,n]$
    \item $t_{m,n,p}:[m,n,p]\to [n,p]$.
    \end{itemize}
  \end{enumerate}

  \item The relations in $\varphi(\dbl{D})$ are as follows:
  \begin{enumerate}
  \item (Functoriality) $[f]\cdot[g]=[f\cdot g],[1_x]=1_{[x]},[\alpha\cdot \beta]=[\alpha]\cdot[\beta],[1_m]=1_{[m]}$ 
  \item (Span maps from cells, laxators, unitors) 
  For a cell $\xinlinecell{x}{y}{z}{w}{f}{g}{m}{m'}{\alpha}$ in 
  $\dbl{D},$ we have the commutative squares: 
  \[\begin{tikzcd}
    {[m]} & {[m']} & {[m]} & {[m']} \\
    {[x]} & {[z]} & {[y]} & {[w]}
    \arrow["{[\alpha]}", from=1-1, to=1-2]
    \arrow["{s_m}"', from=1-1, to=2-1]
    \arrow["{s_{m'}}", from=1-2, to=2-2]
    \arrow["{[\alpha]}", from=1-3, to=1-4]
    \arrow["{t_m}"', from=1-3, to=2-3]
    \arrow["{t_{m'}}", from=1-4, to=2-4]
    \arrow["{[f]}"', from=2-1, to=2-2]
    \arrow["{[g]}"', from=2-3, to=2-4]
  \end{tikzcd}\]
  \[\begin{tikzcd}
    {[m,n]} & {[m]} & {[m,n]} & {[n]} & {[x]} & {[\id_x]} \\
    {[m\odot n]} & {[x]} & {[m\odot n]} & {[z]} & {[\id_x]} & {[x]}
    \arrow["{s_{m,n}}", from=1-1, to=1-2]
    \arrow["{c_{m,n}}"', from=1-1, to=2-1]
    \arrow["{s_m}", from=1-2, to=2-2]
    \arrow["{t_{m,n}}", from=1-3, to=1-4]
    \arrow["{c_{m,n}}"', from=1-3, to=2-3]
    \arrow["{t_n}", from=1-4, to=2-4]
    \arrow["{u_x}", from=1-5, to=1-6]
    \arrow["{u_x}"', from=1-5, to=2-5]
    \arrow[equals, from=1-5, to=2-6]
    \arrow["{s_{\id_x}}", from=1-6, to=2-6]
    \arrow["{s_{m\odot n}}"', from=2-1, to=2-2]
    \arrow["{t_{m\odot n}}"', from=2-3, to=2-4]
    \arrow["{t_{\id_x}}"', from=2-5, to=2-6]
  \end{tikzcd}\]
    \item (Naturality of laxators and unitors) For each composable pair $\alpha,\beta$ of cells as above
  and for each $f:x\to y,$ we have the  commutative squares:
  \[\begin{tikzcd}
    {[m,n]} & {[m',n']} & {[x]} & {[y]} \\
    {[m\odot n]} & {[m'\odot n']} & {[\id_x]} & {[\id_y]}
    \arrow["{[\alpha,\beta]}", from=1-1, to=1-2]
    \arrow["{c_{m,n}}"', from=1-1, to=2-1]
    \arrow["{c_{m',n'}}", from=1-2, to=2-2]
    \arrow["{[f]}", from=1-3, to=1-4]
    \arrow["{u_x}", from=1-3, to=2-3]
    \arrow["{u_y}", from=1-4, to=2-4]
    \arrow["{\alpha\odot\beta}"', from=2-1, to=2-2]
    \arrow["{[\id_f]}"', from=2-3, to=2-4]
  \end{tikzcd}\]

  \item (Candidate pullback squares)
  If $x\xproto{m}y\xproto{n} z\xproto{p} w,$ we have the following commutative squares: 
  
  \[\begin{tikzcd}
    {[m,n]} & {[m]} & {[m,n,p]} & {[m,n]} \\
    {[n]} & {[y]} & {[n,p]} & {[n]}
    \arrow["{s_{m,n}}", from=1-1, to=1-2]
    \arrow["{t_{m,n}}"', from=1-1, to=2-1]
    \arrow["{t_m}", from=1-2, to=2-2]
    \arrow["{s_{m,n,p}}", from=1-3, to=1-4]
    \arrow["{t_{m,n,p}}"', from=1-3, to=2-3]
    \arrow["{t_{m,n}}", from=1-4, to=2-4]
    \arrow["{s_m}"', from=2-1, to=2-2]
    \arrow["{s_{n,p}}"', from=2-3, to=2-4]
  \end{tikzcd}\]

  \item (Maps into candidate pullbacks) 
  For any $m:x\proto y,$ we have: 
  \[\begin{tikzcd}
    {[m]} & {[\id_x,m]} & {[m]} & {[\id_x,m]} \\
    {[x]} & {[\id_x]} & {[x]} & {[m]}
    \arrow["{\ell_m}"', from=1-1, to=1-2]
    \arrow["{s_m}"', from=1-1, to=2-1]
    \arrow["{s_{\id_x,m}}", from=1-2, to=2-2]
    \arrow["{\ell_m}"', from=1-3, to=1-4]
    \arrow["{s_m}"', from=1-3, to=2-3]
    \arrow["{t_{\id_x,m}}", from=1-4, to=2-4]
    \arrow["{u_x}"', from=2-1, to=2-2]
    \arrow["{u_x}"', from=2-3, to=2-4]
  \end{tikzcd}\]
  and similarly for $r_m.$

  For composable triples $m,n,p,$ we have 
  \[\begin{tikzcd}
    {[m,n,p]} & {[m\odot n,p]} & {[m,n,p]} & {[m\odot n,p]} \\
    {[m,n]} & {[m\odot n]} & {[n,p]} & {[p]}
    \arrow["{\ell_{m,n,p}}"', from=1-1, to=1-2]
    \arrow["{s_{m,n,p}}", from=1-1, to=2-1]
    \arrow["{s_{m\odot n,p}}"', from=1-2, to=2-2]
    \arrow["{\ell_{m,n,p}}"', from=1-3, to=1-4]
    \arrow["{t_{m,n,p}}", from=1-3, to=2-3]
    \arrow["{t_{m\odot n,p}}"', from=1-4, to=2-4]
    \arrow["{c_{m,n}}", from=2-1, to=2-2]
    \arrow["{t_{n,p}}", from=2-3, to=2-4]
  \end{tikzcd}\]
  and similarly for $r_{m,n,p}.$

  For loosely composable pairs of cells $\alpha:m\to m',\beta:n\to n',$ we have
  \[\begin{tikzcd}
    {[m,n]} & {[m',n']} & {[m,n]} & {[m',n']} \\
    {[m]} & {[m']} & {[n]} & {[n']}
    \arrow["{[\alpha,\beta]}", from=1-1, to=1-2]
    \arrow["{s_{m,n}}"', from=1-1, to=2-1]
    \arrow["{s_{m',n'}}", from=1-2, to=2-2]
    \arrow["{[\alpha,\beta]}", from=1-3, to=1-4]
    \arrow["{t_{m,n}}"', from=1-3, to=2-3]
    \arrow["{t_{m',n'}}", from=1-4, to=2-4]
    \arrow["\alpha"', from=2-1, to=2-2]
    \arrow["\beta", from=2-3, to=2-4]
  \end{tikzcd}\]

  \item (Associativity and unitality of laxators) Given $x\xproto{m}y\xproto{n}z\xproto{p}w$ in $\dbl{D},$ we have the commutative 
  squares: 
  \[\begin{tikzcd}
    {[m,n,p]} & {[m\odot n,p]} & {[m]} & {[\id_x,m]} \\
    {[m,n\odot p]} & {[m\odot n\odot p]} & {[m,\id_y]} & {[m]}
    \arrow["{\ell_{m,n,p}}"', from=1-1, to=1-2]
    \arrow["{r_{m,n,p}}", from=1-1, to=2-1]
    \arrow["{c_{m\odot n,p}}"', from=1-2, to=2-2]
    \arrow["{\ell_m}", from=1-3, to=1-4]
    \arrow["{r_m}"', from=1-3, to=2-3]
    \arrow["{c_{\id_x,m}}", from=1-4, to=2-4]
    \arrow["{c_{m,n\odot p}}", from=2-1, to=2-2]
    \arrow["{c_{m,\id_y}}"', from=2-3, to=2-4]
  \end{tikzcd}\]

  \end{enumerate}
  \end{itemize}

  Now equip $\varphi(\dbl{D})$ with the structure of a limit sketch by marking the 
  candidate pullback squares listed in the relations above as pullbacks. 
\end{construction}

\begin{proposition}[Category of models is locally presentable]\label{prop:models_locally_presentable}
  The category of models of a simple double theory $\dbl{D}$ in
  $\Span$ (and their strict tight natural transformations)
  is locally presentable, admitting a finitely accessible, faithful, conservative 
  right adjoint functor $U$ to the power $\Set^{\Ob\dbl{D}+\Loose\dbl{D}}$ of the category 
  of sets to the set of objects and loose arrows in $\dbl{D}.$
  This functor $U$ is defined by $U(X)(x) \coloneqq X(x)$ and $U(X)(m) \coloneqq X(m)$.
  \end{proposition}
  \begin{proof}

  Consider a model $X:\varphi(\dbl{D})\to \Set$ of the sketch $\varphi(\dbl{D})$
  constructed above, that is, 
  a functor sending the candidate pullback squares to pullback squares. 
  This provides sets $X[x],$ functions $X[f],$ spans $X[m],$ and span maps $X[\alpha]$ for 
  objects $x,$ tight maps $f,$ loose maps $m,$ and cells $\alpha$ in $\dbl{D},$ functorial
  in $f$ and $\alpha$. 
  
  From the pullback preservation, we may infer $X([m,n])=Xm\times_{Xy} Xn$ 
  and $X([m,n,p])=Xm\times_{Xy} Xn\times_{Xz} Xp$ are sets of paths of heteromorphisms of 
  the appropriate lengths. With the first identification in mind, 
  we will thus interpret the maps $X(c_{m,n}):X([m,n])\to X([m\odot n])$ as laxators and the
  $X(u_x):X([x])\to X([\id_x])$ as unitors, which are natural due to the naturality relations
  in $\varphi(\dbl{D}).$
  
  All the other maps in $\varphi(\dbl{D})$ have their images
  determined by the pullback property and the maps-into-candidate-pullbacks relations. 
  Thus $X$ contains precisely the data of a lax double functor $\bar X$ into $\Span$.
  Then, the last class of relations in $\varphi(\dbl{D})$ guarantees precisely that 
  this data satisfies the axioms of a lax double functor $\bar X$, as desired. 

  A natural transformation $\tau:X\to Y$ between models $X,Y:\varphi(\dbl{D})\to \Set$ 
  effectively has components only at each $[x]$ and each $[m]$, 
  the components at $[m,n]$ and $[m,n,p]$ being determined
  by the pullback preservation condition. The transformation $\tau$ will be natural as soon as it
  is so at all morphisms of $\varphi(\dbl{D})$ whose codomain is not sketched as the apex of 
  a pullback. Thus $\tau$ must be natural at each $[f],$ 
  accounting for naturality of the candidate induced transformation $\bar\tau:\bar X\to \bar Y$ on 
  $\dbl{D}_0$; at the $s_m,t_m,$ accounting for $\tau_{[m]}$ inducing a span map $\bar\tau_m$; 
  at $c_{m,n}$, accounting for respect of $\bar\tau$ for loose composition; at 
  $[\alpha]$, accounting for naturality of $\bar\tau$ on $\dbl{D}_1;$ and finally at $u_x,$
  accounting for respect of $\bar\tau$ for loose units. Thus we see the map 
  $X\mapsto \bar X,\tau\mapsto\bar \tau$ gives an equivalence of categories between 
  models of $\varphi(\dbl{D})$ and models of $\dbl{D}.$

  Since $\mathsf{Mod}(\varphi(\dbl{D}))$ is the category of models of a finite limit sketch, 
  it is well known (\cite[Chapter 1]{locally-presentable-1994}) that it is finitely locally presentable, and indeed that the forgetful 
  functor into $\Set^{\varphi(\dbl{D})}$ is a fully faithful, finitely accessible right adjoint, 
  which is thus also conservative. Restricting along the inclusion 
  $\Ob\dbl{D}+\Loose\dbl{D}\to\varphi(\dbl{D})$ gives a further finitely accessible, faithful and continuous 
  functor $[\varphi(\dbl{D}),\Set]\to \Set^{\Ob\dbl{D}+\Loose\dbl{D}},$ which due to faithfulness and 
  continuity is conservative and due to local presentability is a right adjoint, as promised.
  \end{proof}
  
\begin{remark}[On pseudo and lax natural maps]
Note that, in various cases, we may care more about the pseudonatural or even (co)lax
natural transformations of models, as defined in \cite{cartesian-double-theories-2024}.
The 1-categories of such models will rarely be locally presentable, or even complete 
and cocomplete; one would rather study the 2-categories by analogy with the well-studied
2-categories of algebras for a 2-monad and pseudo, lax, or colax morphisms. We do 
not pursue this here.
\end{remark}

\section{Proof of elements correspondence}\label{app:proof-of-discrete-opfibrations}
In this section, we give the details of the proof of \cref{thm:instances-discrete-opfibrations}. 

We recall the statement of the theorem:

\emph{
Let $B$ be a span-valued model of a simple double theory $\dbl{D}$.
There is an equivalence $\nabla:\Dopf(B)\leftrightarrows \Inst(B):\int$
between the category of discrete opfibrations over $B$
and the category of instances of $B.$
}

\begin{proof}
We repeat the definition of $\nabla$ and $\int$ on objects, interpolating 
some of the checks that were initially skipped.
\begin{enumerate}
\item (The definition of $\nabla$ on objects)
Consider a discrete opfibration $p:E\to B.$ We shall define an instance $\nabla p:1\proto B$
giving the indexed view of $p.$ The material of $\nabla p$ is as follows:
\begin{itemize}
\item On objects, define $\nabla p_x \coloneqq (1\leftarrow Ex\xto{p_x} Bx).$
\item Given the tight morphism $f:x\to y$ in $\dbl{D}$, we define
\begin{equation*}
  \begin{tikzcd}
	{\nabla p_x} \\
	{\nabla p_y}
	\arrow["{\nabla p_f}"', from=1-1, to=2-1]
  \end{tikzcd}
  \quad\coloneqq\quad
  \begin{tikzcd}
	1 & Ex & Bx \\
	1 & Ey & By
	\arrow[equals, from=1-1, to=2-1]
	\arrow[from=1-2, to=1-1]
	\arrow["{p_x}", from=1-2, to=1-3]
	\arrow["Ef"', from=1-2, to=2-2]
	\arrow["Bf", from=1-3, to=2-3]
	\arrow[from=2-2, to=2-1]
	\arrow["{p_y}"', from=2-2, to=2-3]
  \end{tikzcd}.
\end{equation*}
\item Given the loose morphism $m:x\proto y$ in $\dbl{D}$, we define the cell $\nabla p_m$ to be:
\[\begin{tikzcd}
	1 & Ex & Bx & Bm & By \\
	&& Em \\
	1 && Ey && By
	\arrow[equals, from=1-1, to=3-1]
	\arrow[from=1-2, to=1-1]
	\arrow["p_x", from=1-2, to=1-3]
	\arrow[from=1-4, to=1-3]
	\arrow[from=1-4, to=1-5]
	\arrow[equals, from=1-5, to=3-5]
	\arrow[from=2-3, to=1-2]
	\arrow["\lrcorner"{anchor=center, pos=0.125, rotate=135}, draw=none, from=2-3, to=1-3]
	\arrow["{p_m}"{description}, from=2-3, to=1-4]
	\arrow[from=2-3, to=3-3]
	\arrow[from=3-3, to=3-1]
	\arrow[from=3-3, to=3-5]
\end{tikzcd}\]
Here, the marked square is cartesian because
$p$ is a discrete opfibration and the right-hand pentagon commutes because it's a distorted
version of one of the squares witnessing that $pm$ is a span map.
\end{itemize}

Now we check the instance axioms for $\nabla p.$

\begin{itemize}
\item (Functoriality on arrows) Follows from tight functoriality of $p.$
\item (Naturality of actions) Given a cell $\xinlinecell{x}{y}{z}{w}{f}{g}{m}{n}{\alpha}$ in $\dbl{D},$ 
we must show the following two cells are equal, where we have decorated the cells with elements of the most important
sets involved:
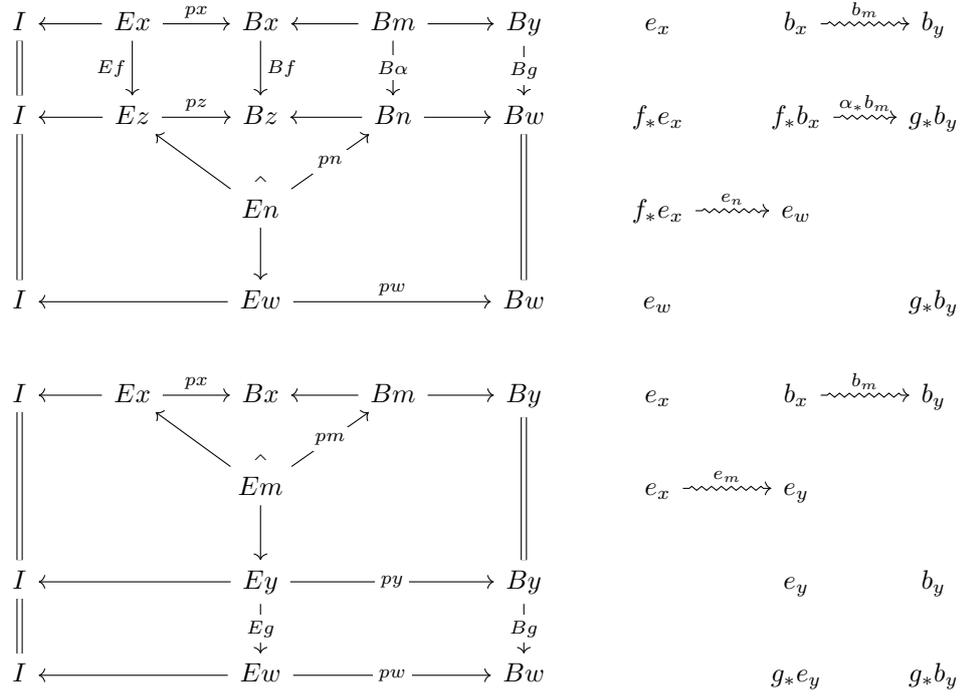
\[\begin{tikzcd}
	I & Ex & Bx & Bm & By & {e_x} & {b_x} & {b_y} \\
	I & Ez & Bz & Bn & Bw & {f_*e_x} & {f_*b_x} & {g_*b_y} \\
	&& En &&& {f_*e_x} & {e_w} \\
	I && Ew && Bw & {e_w} && {g_*b_y} \\
	I & Ex & Bx & Bm & By & {e_x} & {b_x} & {b_y} \\
	&& Em &&& {e_x} & {e_y} \\
	I && Ey && By && {e_y} & {b_y} \\
	I && Ew && Bw && {g_*e_y} & {g_*b_y}
	\arrow[equals, from=1-1, to=2-1]
	\arrow[from=1-2, to=1-1]
	\arrow["px", from=1-2, to=1-3]
	\arrow["Ef"{description}, from=1-2, to=2-2]
	\arrow["Bf"{description}, from=1-3, to=2-3]
	\arrow[from=1-4, to=1-3]
	\arrow[from=1-4, to=1-5]
	\arrow["{B\alpha}"{description}, from=1-4, to=2-4]
	\arrow["Bg"{description}, from=1-5, to=2-5]
	\arrow["{b_m}", squiggly, from=1-7, to=1-8]
	\arrow[equals, from=2-1, to=4-1]
	\arrow[from=2-2, to=2-1]
	\arrow["pz", from=2-2, to=2-3]
	\arrow[from=2-4, to=2-3]
	\arrow[from=2-4, to=2-5]
	\arrow[equals, from=2-5, to=4-5]
	\arrow["{\alpha_*b_m}", squiggly, from=2-7, to=2-8]
	\arrow[from=3-3, to=2-2]
	\arrow["\lrcorner"{anchor=center, pos=0.125, rotate=135}, draw=none, from=3-3, to=2-3]
	\arrow["pn"{description}, from=3-3, to=2-4]
	\arrow[from=3-3, to=4-3]
	\arrow["{e_m}", squiggly, from=3-6, to=3-7]
	\arrow[from=4-3, to=4-1]
	\arrow["pw"{description}, from=4-3, to=4-5]
	\arrow[equals, from=5-1, to=7-1]
	\arrow[from=5-2, to=5-1]
	\arrow["px", from=5-2, to=5-3]
	\arrow[from=5-4, to=5-3]
	\arrow[from=5-4, to=5-5]
	\arrow[equals, from=5-5, to=7-5]
	\arrow["{b_m}", squiggly, from=5-7, to=5-8]
	\arrow[from=6-3, to=5-2]
	\arrow["\lrcorner"{anchor=center, pos=0.125, rotate=135}, draw=none, from=6-3, to=5-3]
	\arrow["pm"{description}, from=6-3, to=5-4]
	\arrow[from=6-3, to=7-3]
	\arrow["{e_m}", squiggly, from=6-6, to=6-7]
	\arrow[equals, from=7-1, to=8-1]
	\arrow[from=7-3, to=7-1]
	\arrow["py"{description}, from=7-3, to=7-5]
	\arrow["Eg"{description}, from=7-3, to=8-3]
	\arrow["Bg"{description}, from=7-5, to=8-5]
	\arrow[from=8-3, to=8-1]
	\arrow["pw"{description}, from=8-3, to=8-5]
\end{tikzcd}\]
Here we understand, for instance, $p_x(e_x)=b_x$ and $Ef(e_x)=f_*e_x,$ etc.

Now, by naturality of $p$ at $\alpha,$ we have $e_n=\alpha_*e_m,$ and in particular 
$g_*e_y=e_w,$ which proves the commutativity we desire.

\item (Associativity of actions) Given $x\xproto{m}y\xproto{n}z$ in $\dbl{D}$, we are to show the following
two cells are equal:
\[\begin{tikzcd}[column sep=tiny]
	1 & Ex & Bx & Bm & By & Bn & Bz & {e_x} & {b_x} & {b_y} & {b_z} \\
	&& Em &&&&& {e_x} & {e_y} \\
	1 && Ey && By & Bn & Bz && {e_y} & {b_y} & {b_z} \\
	&&&& En &&&& {e_y} & {e_{z_1}} \\
	1 &&& Ez &&& Xz &&& {e_{z_1}} & {b_z} \\
	1 & Ex & Bx & Bm & By & Bn & Bz & {e_x} & {b_x} & {b_y} & {b_z} \\
	&&&& \bullet \\
	1 & Ex & Bx && {B(m\odot n)} && Bz & {e_x} & {b_x} && {b_z} \\
	&&& {E(m\odot n)} &&&& {e_x} & {e_{z_2}} \\
	1 &&& Ez &&& Bz && {e_{z_2}}
	\arrow[equals, from=1-1, to=3-1]
	\arrow[from=1-2, to=1-1]
	\arrow["px", from=1-2, to=1-3]
	\arrow[from=1-4, to=1-3]
	\arrow[from=1-4, to=1-5]
	\arrow[equals, from=1-5, to=3-5]
	\arrow[from=1-6, to=1-5]
	\arrow[from=1-6, to=1-7]
	\arrow[equals, from=1-6, to=3-6]
	\arrow[equals, from=1-7, to=3-7]
	\arrow["{b_m}", squiggly, from=1-9, to=1-10]
	\arrow["{b_n}", squiggly, from=1-10, to=1-11]
	\arrow[from=2-3, to=1-2]
	\arrow["\lrcorner"{anchor=center, pos=0.125, rotate=135}, draw=none, from=2-3, to=1-3]
	\arrow["pm"{description}, from=2-3, to=1-4]
	\arrow[from=2-3, to=3-3]
	\arrow["{e_m}", squiggly, from=2-8, to=2-9]
	\arrow[equals, from=3-1, to=5-1]
	\arrow[from=3-3, to=3-1]
	\arrow["py"{description}, from=3-3, to=3-5]
	\arrow[from=3-6, to=3-5]
	\arrow[from=3-6, to=3-7]
	\arrow[equals, from=3-7, to=5-7]
	\arrow["{b_n}", squiggly, from=3-10, to=3-11]
	\arrow[from=4-5, to=3-3]
	\arrow["\lrcorner"{anchor=center, pos=0.125, rotate=135}, draw=none, from=4-5, to=3-5]
	\arrow["pn"{description}, from=4-5, to=3-6]
	\arrow[from=4-5, to=5-4]
	\arrow["{e_n}", squiggly, from=4-9, to=4-10]
	\arrow[from=5-4, to=5-1]
	\arrow[from=5-4, to=5-7]
	\arrow[equals, from=6-1, to=8-1]
	\arrow[from=6-2, to=6-1]
	\arrow["px", from=6-2, to=6-3]
	\arrow[equals, from=6-2, to=8-2]
	\arrow[equals, from=6-3, to=8-3]
	\arrow[from=6-4, to=6-3]
	\arrow[from=6-4, to=6-5]
	\arrow[from=6-6, to=6-5]
	\arrow[from=6-6, to=6-7]
	\arrow[equals, from=6-7, to=8-7]
	\arrow["{b_m}", squiggly, from=6-9, to=6-10]
	\arrow["{b_n}", squiggly, from=6-10, to=6-11]
	\arrow[from=7-5, to=6-4]
	\arrow["\lrcorner"{anchor=center, pos=0.125, rotate=135}, draw=none, from=7-5, to=6-5]
	\arrow[from=7-5, to=6-6]
	\arrow[from=7-5, to=8-5]
	\arrow[equals, from=8-1, to=10-1]
	\arrow[from=8-2, to=8-1]
	\arrow["px", from=8-2, to=8-3]
	\arrow[""{name=0, anchor=center, inner sep=0}, from=8-5, to=8-3]
	\arrow[from=8-5, to=8-7]
	\arrow[equals, from=8-7, to=10-7]
	\arrow["{b_{m\odot n}}", squiggly, from=8-9, to=8-11]
	\arrow[from=9-4, to=8-2]
	\arrow["{p(m\odot n)}"{description}, from=9-4, to=8-5]
	\arrow[from=9-4, to=10-4]
	\arrow["{e_{m\odot n}}", squiggly, from=9-8, to=9-9]
	\arrow[from=10-4, to=10-1]
	\arrow["pz", from=10-4, to=10-7]
	\arrow["\lrcorner"{anchor=center, pos=0.125, rotate=135}, draw=none, from=9-4, to=0]
\end{tikzcd}\]

Now by loose functoriality and the discrete opfibration propertyof $p,$ we can conclude that 
$e_{m\odot n}=E_{m,n}(e_m,e_n),$ so in particular $e_{z_2}=e_{z_1},$ which is what we needed to show.

\item (Unitality of actions) Finally, for $x\in \dbl{D}$ we must show that 
the square below is the identity: 
\[\begin{tikzcd}
	1 & Ex & Bx & Bx & Bx \\
	\\
	1 & Ex & Bx & {B\id_x} & Bx \\
	&& {E\id_x} \\
	1 && Ex && Bx
	\arrow[equals, from=1-1, to=3-1]
	\arrow[from=1-2, to=1-1]
	\arrow["px", from=1-2, to=1-3]
	\arrow[equals, from=1-2, to=3-2]
	\arrow[equals, from=1-3, to=3-3]
	\arrow[from=1-4, to=1-3]
	\arrow[from=1-4, to=1-5]
	\arrow["{B_x}"{description}, from=1-4, to=3-4]
	\arrow[equals, from=1-5, to=3-5]
	\arrow[equals, from=3-1, to=5-1]
	\arrow[from=3-2, to=3-1]
	\arrow["px", from=3-2, to=3-3]
	\arrow[from=3-4, to=3-3]
	\arrow[from=3-4, to=3-5]
	\arrow[equals, from=3-5, to=5-5]
	\arrow[from=4-3, to=3-2]
	\arrow["{p\id_x}"{description}, from=4-3, to=3-4]
	\arrow[from=4-3, to=5-3]
	\arrow[from=5-3, to=5-1]
	\arrow["px", from=5-3, to=5-5]
\end{tikzcd}\]
This is to prove that, given a choice of $b_x\in Bx$ and a lift $e_x\in Ex$ of $b_x$ along 
$p,$ then the unique lift of $B_{\id_x}(b_x):b_x\proto b_x$ along $p$ out of $e_x$ has codomain
$e_x.$ Since $E_{\id_x}(e_x)$ certainly has codomain $e_x,$ this follows from loose 
functoriality of $p.$
\end{itemize}

\item (The definition of $\int$ on objects)
Passing in the other direction, consider an instance $P:1\proto B.$ We shall define
the fibered view $\pi:\int P\to B.$ The material of $\int P$ and $\pi$ is as follows:
\begin{itemize}
\item (On objects) Given $x\in\dbl{D},$ write $Px$ as $1\leftarrow \int P x\xto{\pi_x} Bx$ 
to define $\int P$ and $\pi$ on $x.$
\item (On tight morphisms) For $f:x\to y$ in $\dbl{D},$ we set $\int P f \coloneqq Pf$ as below,
noting that the commutative square on the right will become a strict naturality
square for $\pi$:
\[\begin{tikzcd}
	1 & {\int Px} & Bx \\
	1 & {\int Py} & By
	\arrow[equals, from=1-1, to=2-1]
	\arrow[from=1-2, to=1-1]
	\arrow["{\pi_x}", from=1-2, to=1-3]
	\arrow["{\int Pf}"', from=1-2, to=2-2]
	\arrow["Bf", from=1-3, to=2-3]
	\arrow[from=2-2, to=2-1]
	\arrow["{\pi_y}"', from=2-2, to=2-3]
\end{tikzcd}\]
\item (On loose morphisms) For $m:x\proto y$ in $\dbl{D},$ we are given the cell $Pm$ below,
where again we define $\int Pm$ and $\pi_m$ by naming appropriate components of $P.$
Note that we are ensuring that $\pi$ will be a discrete opfibration by our choice of $\int Pm$ and
this diagram already contains the fact that $\int Pm$ is a span map.
\[\begin{tikzcd}
	1 & {\int P x} & Bx & Bm & By \\
	&& {\int Pm} \\
	1 && {\int P_y} && By
	\arrow[equals, from=1-1, to=3-1]
	\arrow[from=1-2, to=1-1]
	\arrow["{\pi_x}", from=1-2, to=1-3]
	\arrow[from=1-4, to=1-3]
	\arrow[from=1-4, to=1-5]
	\arrow[equals, from=1-5, to=3-5]
	\arrow[from=2-3, to=1-2]
	\arrow["\lrcorner"{anchor=center, pos=0.125, rotate=135}, draw=none, from=2-3, to=1-3]
	\arrow["{\pi_m}"{description}, from=2-3, to=1-4]
	\arrow[from=2-3, to=3-3]
	\arrow[from=3-3, to=3-1]
	\arrow["{\pi_y}", from=3-3, to=3-5]
\end{tikzcd}\]

\item (On cells) Given a cell $\xinlinecell{x}{y}{z}{w}{f}{g}{m}{n}{\alpha}$ in $\dbl{D},$
we have available so far the data below, where all visible squares are already known to commute.
\[\begin{tikzcd}
	{\int P x} & {\int P m} & {\int P y} \\
	{\int P z} & {\int P n} & {\int P w} \\
	Bx & Bm & By \\
	Bz & Bn & Bw
	\arrow["{\int P f}", from=1-1, to=2-1]
	\arrow["{\pi_x}"{description}, curve={height=24pt}, dotted, from=1-1, to=3-1]
	\arrow["{s_m}"', from=1-2, to=1-1]
	\arrow["{t_m}", from=1-2, to=1-3]
	\arrow["{\int P\alpha}", dashed, from=1-2, to=2-2]
	\arrow["{\pi_m}"{description, pos=0.2}, curve={height=12pt}, dotted, from=1-2, to=3-2]
	\arrow["{\int P g}", from=1-3, to=2-3]
	\arrow["{\pi_y}"{description}, curve={height=-24pt}, dotted, from=1-3, to=3-3]
	\arrow["{\pi_z}"{description}, curve={height=24pt}, dotted, from=2-1, to=4-1]
	\arrow["{s_n}", from=2-2, to=2-1]
	\arrow["{t_n}"', from=2-2, to=2-3]
	\arrow["{\pi_n}"{description, pos=0.3}, curve={height=-24pt}, dotted, from=2-2, to=4-2]
	\arrow["{\pi_w}"{description}, curve={height=-24pt}, dotted, from=2-3, to=4-3]
	\arrow["Bf", from=3-1, to=4-1]
	\arrow[from=3-2, to=3-1]
	\arrow[from=3-2, to=3-3]
	\arrow["{B\alpha}", from=3-2, to=4-2]
	\arrow["Bg", from=3-3, to=4-3]
	\arrow[from=4-2, to=4-1]
	\arrow[from=4-2, to=4-3]
\end{tikzcd}\]

Now recall that we defined $\int P n \coloneqq \int Pz \times_{Bz} Bn$, so that
by the universal property of the pullback we may define 
$\int P\alpha$ to be the unique map such that $(\int P\alpha)\cdot \pi_n=\pi_m\cdot B\alpha$ 
and $(\int P\alpha)\cdot s_n=s_m\cdot (\int P f).$ These two conditions ensure 
respectively that $\pi$ is natural at $\alpha$ and that the left-hand square witnessing that
$\int P \alpha$ is a span map commutes. 

As for the right-hand square, 
recall that $t_m: \int P x\times_{Bx} Bm\to \int P y$ is precisely the action map 
from $P$; thus the desired equation $t_m\cdot \int P g=\int P\alpha\cdot t_n$ says 
that $g_*(p_x\cdot b_m)=(f_*p_x)\cdot \alpha_*b_m,$ for any 
$p_x\in Px, b_m:b_x\proto b_y,$ which is exactly the naturality of action axiom for the instance $P.$

\item (Laxators)
Given composable loose morphisms $x\xproto{m} y\xproto{n} z$ in $\dbl{D},$ we have already defined
everything except $q$ in the diagram below. This map $q$ is defined 
by the universal property of the pullback at $\bullet_B$ in such a 
way as to produce a span map between the spans with apexes $\bullet_E$ and $\bullet_B.$

\[\begin{tikzcd}
	{\int P x} & {\int P m} & {\int P y} & {\int P n} & {\int P z} \\
	&& {\bullet_E} \\
	{\int P x} && {\int P (m\odot n)} && {\int P z} \\
	Bx & Bm & By & Bn & Bz \\
	&& {\bullet_B} \\
	Bx && {B(m\odot n)} && Bz
	\arrow[equals, from=1-1, to=3-1]
	\arrow[from=1-2, to=1-1]
	\arrow[from=1-2, to=1-3]
	\arrow["\pi"{description}, dotted, from=1-2, to=4-2]
	\arrow[from=1-4, to=1-3]
	\arrow[from=1-4, to=1-5]
	\arrow["\pi"{description}, dotted, from=1-4, to=4-4]
	\arrow[equals, from=1-5, to=3-5]
	\arrow[from=2-3, to=1-2]
	\arrow["\lrcorner"{anchor=center, pos=0.125, rotate=135}, draw=none, from=2-3, to=1-3]
	\arrow[from=2-3, to=1-4]
	\arrow["{\int P_{m,n}}"', dashed, from=2-3, to=3-3]
	\arrow["q"{description}, curve={height=-24pt}, dotted, from=2-3, to=5-3]
	\arrow["{s_{m\odot n}}"{description, pos=0.7}, from=3-3, to=3-1]
	\arrow[from=3-3, to=3-5]
	\arrow["\pi"{description}, curve={height=12pt}, dotted, from=3-3, to=6-3]
	\arrow[equals, from=4-1, to=6-1]
	\arrow[from=4-2, to=4-1]
	\arrow[from=4-2, to=4-3]
	\arrow[from=4-4, to=4-3]
	\arrow[from=4-4, to=4-5]
	\arrow[equals, from=4-5, to=6-5]
	\arrow[from=5-3, to=4-2]
	\arrow["\lrcorner"{anchor=center, pos=0.125, rotate=135}, draw=none, from=5-3, to=4-3]
	\arrow[from=5-3, to=4-4]
	\arrow["{B_{m,n}}", from=5-3, to=6-3]
	\arrow[from=6-3, to=6-1]
	\arrow[from=6-3, to=6-5]
\end{tikzcd}\]

This permits us to define the laxator
\begin{equation*} \textstyle
  \int P_{m,n}:\bullet_E \to \int P(m\odot n)=B(m\odot n)\times_{Bx} \int P x
\end{equation*}
using the universal property of the pullback. Specifically, we require that
$\int P_{m,n}\cdot \pi_{m\odot n}=q\cdot B_{m,n},$ which ensures loose functoriality of $\pi$ at $m,n$,
and that $\int P_{m,n}\cdot s_{m\odot n}$ makes the left-hand pentagon commute.
To make $\int P_{m,n}$ a span map we must check that the right-hand pentagon commutes as well;
this is precisely the associativity of actions for $P$ at $m,n.$

Note that the domain 
$\int Px\times_{Bx} Bm\times_{\int Py} \int Py\times_{By} B_n$ of the laxator is canonically
isomorphic to $\int Px\times_{Bx} Bm\times_{By} Bn$ while the codomain is 
$\int Px\times_{Bx} B(m\odot n).$ In these terms the laxator 
simply maps $(p_x,b_m,b_n)$ to $(p_x,b_m\odot b_n),$ where we generalize 
$\odot$ to denote $B_{m,n}(b_m,b_n).$

\item (Unitors) 
Similarly, we are given the data to uniquely define a unitor
from the universal property of the pullback $\int P \id_x=\int P x \times_{Bx} B\id_x$
by filling the upper boundary below:
\[\begin{tikzcd}
	{\int Px} & {\int Px} & {\int Px} \\
	{\int Px} & {\int P \id_x} & {\int Px} \\
	Bx & Bx & Bx \\
	Bx & {B\id_x} & Bx
	\arrow[equals, from=1-1, to=2-1]
	\arrow[curve={height=18pt}, dotted, from=1-1, to=3-1]
	\arrow[from=1-2, to=1-1]
	\arrow[from=1-2, to=1-3]
	\arrow["{P_x}"', dashed, from=1-2, to=2-2]
	\arrow[curve={height=-18pt}, dotted, from=1-2, to=3-2]
	\arrow[equals, from=1-3, to=2-3]
	\arrow[curve={height=18pt}, dotted, from=2-1, to=4-1]
	\arrow[from=2-2, to=2-1]
	\arrow[from=2-2, to=2-3]
	\arrow[curve={height=-18pt}, dotted, from=2-2, to=4-2]
	\arrow[equals, from=3-1, to=4-1]
	\arrow[from=3-2, to=3-1]
	\arrow[from=3-2, to=3-3]
	\arrow["{B_x}"{description}, from=3-2, to=4-2]
	\arrow[equals, from=3-3, to=4-3]
	\arrow[from=4-2, to=4-1]
	\arrow[from=4-2, to=4-3]
\end{tikzcd}\]

Thus $\int P_x\cdot \pi_{\id_x}=\pi_x\cdot B_x,$ which accounts for loose functoriality of $\pi$ at $\id_x,$
while $\int P_x \cdot s_x=\id_{\int Px},$ which accounts for half of the span map property. The fact 
that this definition makes $\int P_x\cdot t_x=\id_{\int Px}$ as well follows from the fact that
$t_x$ is the action map for $P$, and unitality of this action says that acting by something in the
image of $B_x$ is trivial.
\end{itemize}

We must check that $\int P$ is a lax double functor.

\begin{itemize}
\item (Internal functoriality)
Functoriality of $\int P$ on arrows follows from that of $P.$ Functoriality 
on cells follows from that of $B,$ functoriality of $\int P$ on arrows, 
and the pullback property of $\int P$ on loose arrows.
\item (Naturality of laxators)
For this axiom, given cells 
$\xinlinecell{x}{y}{x'}{y'}{f}{g}{m}{m'}{\alpha},
\xinlinecell{y}{z}{y'}{z'}{g}{h}{n}{n'}{\beta}$ we have to show the two 
canonical maps of spans shown below coincide:
\[\begin{tikzcd}[column sep=tiny]
	{\int Px} & {\int Pm} & {\int Py} & {\int Pn} & {\int P z} & {p_x} & {b_x} & {b_y} & {b_z} \\
	{\int Px'} & {\int Pm'} & {\int Py'} & {\int Pn'} & {\int P z'} & {f_*p_x} & {f_*b_x} & {g_*b_y} & {h_*b_z} \\
	&& \bullet \\
	{\int P x'} && {\int P (m'\odot n')} && {\int P z'} & {f_*p_x} & {f_*b_x} && {h_*b_z} \\
	{\int Px} & {\int Pm} & {\int Py} & {\int Pn} & {\int P z} & {p_x} & {b_x} & {b_y} & {b_z} \\
	&& \bullet \\
	{\int P x} && {\int P (m\odot n)} && {\int P z} & {p_x} & {b_x} && {b_z} \\
	{\int Px'} && {\int P(m'\odot n')} && {\int Pz'} & {f_*p_x} & {f_*b_x} && {h_*b_z}
	\arrow["{\int P f}"', from=1-1, to=2-1]
	\arrow[from=1-2, to=1-1]
	\arrow[from=1-2, to=1-3]
	\arrow["{\int P \alpha}"', from=1-2, to=2-2]
	\arrow["{\int P g}"', from=1-3, to=2-3]
	\arrow[from=1-4, to=1-3]
	\arrow[from=1-4, to=1-5]
	\arrow["{\int P \beta}"', from=1-4, to=2-4]
	\arrow["{\int P h}"', from=1-5, to=2-5]
	\arrow["{b_m}", squiggly, from=1-7, to=1-8]
	\arrow["{b_n}", squiggly, from=1-8, to=1-9]
	\arrow[equals, from=2-1, to=4-1]
	\arrow[from=2-2, to=2-1]
	\arrow[from=2-2, to=2-3]
	\arrow[from=2-4, to=2-3]
	\arrow[from=2-4, to=2-5]
	\arrow[equals, from=2-5, to=4-5]
	\arrow["{\alpha_*b_m}", squiggly, from=2-7, to=2-8]
	\arrow["{\beta_*b_n}", squiggly, from=2-8, to=2-9]
	\arrow[from=3-3, to=2-2]
	\arrow["\lrcorner"{anchor=center, pos=0.125, rotate=135}, draw=none, from=3-3, to=2-3]
	\arrow[from=3-3, to=2-4]
	\arrow["{\int P_{m',n'}}"', from=3-3, to=4-3]
	\arrow[from=4-3, to=4-1]
	\arrow[from=4-3, to=4-5]
	\arrow["{\alpha_* b_m\odot \beta_* b_n}", squiggly, from=4-7, to=4-9]
	\arrow[equals, from=5-1, to=7-1]
	\arrow[from=5-2, to=5-1]
	\arrow[from=5-2, to=5-3]
	\arrow[from=5-4, to=5-3]
	\arrow[from=5-4, to=5-5]
	\arrow[equals, from=5-5, to=7-5]
	\arrow["{b_m}", squiggly, from=5-7, to=5-8]
	\arrow["{b_n}", squiggly, from=5-8, to=5-9]
	\arrow[from=6-3, to=5-2]
	\arrow["\lrcorner"{anchor=center, pos=0.125, rotate=135}, draw=none, from=6-3, to=5-3]
	\arrow[from=6-3, to=5-4]
	\arrow["{\int P_{m,n}}"', from=6-3, to=7-3]
	\arrow["{\int P f}"', from=7-1, to=8-1]
	\arrow[from=7-3, to=7-1]
	\arrow[from=7-3, to=7-5]
	\arrow["{\int P (\alpha\odot \beta)}"', from=7-3, to=8-3]
	\arrow["{\int P h}"', from=7-5, to=8-5]
	\arrow["{b_m\odot b_n}", squiggly, from=7-7, to=7-9]
	\arrow[from=8-3, to=8-1]
	\arrow[from=8-3, to=8-5]
	\arrow["{(\alpha\odot \beta)_*(b_m\odot b_n)}", squiggly, from=8-7, to=8-9]
\end{tikzcd}\]

In fact, thinking elementwise, it suffices to note that 
$(\alpha\odot \beta)_*(b_m\odot b_n)=\alpha_* b_m\odot \beta_* b_n,$ 
which is an instance of naturality of laxators for $B.$ That is, 
naturality of laxators for $\int P$ does not actually depend on the 
instance axioms for $P.$

\item (Naturality of unitors)
This follows easily from naturality of $\pi$ together with naturality
of unitors for $B$:
\[\begin{tikzcd}
	{\int Px} & {\int Px} & {\int Px} & {p_x} \\
	{\int Px} & {\int P _{\id x}} & {\int P x} & {p_x} & {b_x} & {b_x} \\
	{\int P y} & {\int P_{\id y}} & {\int P y} & {f_*p_x} & {f_*b_x} & {f_*b_x} \\
	{\int Px} & {\int Px} & {\int Px} & {p_x} \\
	{\int Py} & {\int P y} & {\int P y} & {f_*p_x} \\
	{\int P y} & {\int P_{\id y}} & {\int P y} & {f_*p_x} & {\pi_y(f_*p_x)} & {\pi_y(f_*p_x)}
	\arrow[equals, from=1-1, to=2-1]
	\arrow[equals, from=1-2, to=1-1]
	\arrow[equals, from=1-2, to=1-3]
	\arrow["{\int P_x}"', from=1-2, to=2-2]
	\arrow[equals, from=1-3, to=2-3]
	\arrow["{\int Pf}"', from=2-1, to=3-1]
	\arrow[from=2-2, to=2-1]
	\arrow[from=2-2, to=2-3]
	\arrow["{\int P\id f}"', from=2-2, to=3-2]
	\arrow["{\int P f}"', from=2-3, to=3-3]
	\arrow["\shortmid"{marking}, equals, from=2-5, to=2-6]
	\arrow[from=3-2, to=3-1]
	\arrow[from=3-2, to=3-3]
	\arrow["{f_*\id_{b_x}}"{inner sep=.8ex}, "\shortmid"{marking}, from=3-5, to=3-6]
	\arrow["{\int P f}"', from=4-1, to=5-1]
	\arrow[equals, from=4-2, to=4-1]
	\arrow[equals, from=4-2, to=4-3]
	\arrow["{\int P f}"', from=4-2, to=5-2]
	\arrow["{\int P f}"', from=4-3, to=5-3]
	\arrow[equals, from=5-1, to=6-1]
	\arrow[equals, from=5-2, to=5-1]
	\arrow[equals, from=5-2, to=5-3]
	\arrow["{\int P_y}"', from=5-2, to=6-2]
	\arrow[equals, from=5-3, to=6-3]
	\arrow[from=6-2, to=6-1]
	\arrow[from=6-2, to=6-3]
	\arrow["\shortmid"{marking}, equals, from=6-5, to=6-6]
\end{tikzcd}\]

\item (Associativity)
Given a chain $w\xproto{m} x\xproto{n} y\xproto{p} z$ 
in $\dbl{D}$ and some element 
$(p_w,b_w\xrightsquigarrow{b_m} b_x\xrightsquigarrow{b_n} b_y\xrightsquigarrow{b_p} b_z)$ 
of 
$\int Pm\times_{\int Px}\int Pn\times_{\int Py}\int Pp$, which is canonically isomorphic 
to $\int Pw\times_{Bw} Bm\times_{Bx}Bn\times_{By}Bp,$
we are asked to show that $b_w\cdot ((b_x\odot b_y)\odot b_z)=b_w\cdot (b_x\odot (b_y\odot b_z)).$
This again follows simply from associativity for $B.$ 

\item (Unitality)
Once more, here, given $m:x\proto y$ in $\dbl{D}$ and 
$p_m:p_x\proto p_y$ in $\int P m,$ we are asked to show that 
$\id_{p_x}\odot p_m=p_m=p_m\odot \id_{p_y}.$ Since $\int Pm=\int P x\times_{Bx} Bm,$ 
this is simply to show that $(p_x, (\id_{b_x}\odot b_m))=(p_x, b_m)=(p_x,b_m\odot\id_{b_y}),$
which is just unitality for $B.$
\end{itemize}

The axioms for $\pi$ were verified in passing during the construction above, so we have 
finished the construction of our discrete opfibration over $B.$

This completes the constructions of $\nabla$ and $\int$ on objects. We now move 
on to morphisms. 

\item (The definition of $\nabla$ on morphisms)
Now take two discrete opfibrations $p,p':E,E'\to B$ between models of $\dbl{D}$
and a map $f:E\to E',$ which is just a strict natural transformation over 
$B.$ (But note that the discrete opfibration condition means that such a morphism
is determined entirely by its object components.)

We must define $\nabla f:\nabla p\to\nabla p'.$ For each $x\in\dbl{D},$
we are given the below map of spans, which we shall use as the component
$\nabla f_x$ of $\nabla f$ at $x$:
\[\begin{tikzcd}
	1 & Ex & Bx \\
	1 & {E'x} & Bx
	\arrow[equals, from=1-1, to=2-1]
	\arrow[from=1-2, to=1-1]
	\arrow["{p_x}"', from=1-2, to=1-3]
	\arrow["{f_x}"{description}, from=1-2, to=2-2]
	\arrow[equals, from=1-3, to=2-3]
	\arrow[from=2-2, to=2-1]
	\arrow["{p'_x}", from=2-2, to=2-3]
\end{tikzcd}\]
This definition will make $\nabla$ functorial, since we shall have $\nabla(f\cdot g)_x=(f\cdot g)_x=f_x\cdot g_x$
and similarly for identities.

We check that this makes $\nabla f$ into a morphism of instances:
\begin{itemize}
\item (Equivariance)
We must show the following pastings are equal:
\[\begin{tikzcd}
	1 & Ex & Bx & Bm & By \\
	1 & {E'x} & Bx & Bm & By \\
	&& {E'm} \\
	1 && {E'y} && By \\
	1 & Ex & Bx & Bm & By \\
	&& Em \\
	1 && Ey && By \\
	1 && {E'y} && By
	\arrow[equals, from=1-1, to=2-1]
	\arrow[from=1-2, to=1-1]
	\arrow["px", from=1-2, to=1-3]
	\arrow["{f_x}"', from=1-2, to=2-2]
	\arrow[equals, from=1-3, to=2-3]
	\arrow[from=1-4, to=1-3]
	\arrow[from=1-4, to=1-5]
	\arrow[equals, from=1-4, to=2-4]
	\arrow[equals, from=1-5, to=2-5]
	\arrow[equals, from=2-1, to=4-1]
	\arrow[from=2-2, to=2-1]
	\arrow["{p'x}", from=2-2, to=2-3]
	\arrow[from=2-4, to=2-3]
	\arrow[from=2-4, to=2-5]
	\arrow[equals, from=2-5, to=4-5]
	\arrow[from=3-3, to=2-2]
	\arrow["\lrcorner"{anchor=center, pos=0.125, rotate=135}, draw=none, from=3-3, to=2-3]
	\arrow["{p'm}"{description}, from=3-3, to=2-4]
	\arrow[from=3-3, to=4-3]
	\arrow[from=4-3, to=4-1]
	\arrow[from=4-3, to=4-5]
	\arrow[equals, from=5-1, to=7-1]
	\arrow[from=5-2, to=5-1]
	\arrow["{p_x}", from=5-2, to=5-3]
	\arrow[from=5-4, to=5-3]
	\arrow[from=5-4, to=5-5]
	\arrow[equals, from=5-5, to=7-5]
	\arrow[from=6-3, to=5-2]
	\arrow["\lrcorner"{anchor=center, pos=0.125, rotate=135}, draw=none, from=6-3, to=5-3]
	\arrow["{p_m}", from=6-3, to=5-4]
	\arrow[from=6-3, to=7-3]
	\arrow[equals, from=7-1, to=8-1]
	\arrow[from=7-3, to=7-1]
	\arrow["{p_y}", from=7-3, to=7-5]
	\arrow["{f_y}"{description}, from=7-3, to=8-3]
	\arrow[equals, from=7-5, to=8-5]
	\arrow[from=8-3, to=8-1]
	\arrow["{p'_y}", from=8-3, to=8-5]
\end{tikzcd}\]
To prove this, consider $p(e_x)=b_x$ and $b_m\in Bm(b_x,b_y).$
We can use the discrete opfibration property to get $e_m:e_x\proto e_y$
over $b_m,$ or to $f_x(e_x)$ to get $e'_m:f_x(e_x)\proto e'_y,$ and
the pasting says that $f_y(e_y)=e'_y.$ This is true because both elements of
$E'y$ are the codomain of lifts of $b_m$ with domain $f_x(e_x).$ 
\item (Naturality)
Given $t:x\to y$ in $\dbl{D},$ we must show the following cells coincide:
\[\begin{tikzcd}
	1 & Ex & Bx && 1 & Ex & Bx \\
	1 & Ey & By & {=} & 1 & {E'x} & Bx \\
	1 & {E'y} & By && 1 & {E'y} & By
	\arrow[equals, from=1-1, to=2-1]
	\arrow[from=1-2, to=1-1]
	\arrow["{p_x}", from=1-2, to=1-3]
	\arrow["Et"{description}, from=1-2, to=2-2]
	\arrow["Bt"', from=1-3, to=2-3]
	\arrow[equals, from=1-5, to=2-5]
	\arrow[from=1-6, to=1-5]
	\arrow["{p_x}", from=1-6, to=1-7]
	\arrow["{f_x}"{description}, from=1-6, to=2-6]
	\arrow[equals, from=1-7, to=2-7]
	\arrow[equals, from=2-1, to=3-1]
	\arrow[from=2-2, to=2-1]
	\arrow["{p_y}", from=2-2, to=2-3]
	\arrow["{f_y}"{description}, from=2-2, to=3-2]
	\arrow[equals, from=2-3, to=3-3]
	\arrow[equals, from=2-5, to=3-5]
	\arrow[from=2-6, to=2-5]
	\arrow["{p'_x}", from=2-6, to=2-7]
	\arrow["{E't}"{description}, from=2-6, to=3-6]
	\arrow["Bt"{description}, from=2-7, to=3-7]
	\arrow[from=3-2, to=3-1]
	\arrow["{p'_y}", from=3-2, to=3-3]
	\arrow[from=3-6, to=3-5]
	\arrow["{p'_y}", from=3-6, to=3-7]
\end{tikzcd}\]
But this is just naturality of $f.$
\end{itemize}

\item (The definition of $\int$ on morphisms)
Now take two instances $P,P':1\proto B$ of $B$ and a morphism 
$\mu:P\to P'.$ We must define $\int \mu:(\pi:\int P\to B)\to (\pi':\int P'\to B).$
For each $x\in \dbl{D},$ we are given the below map of spans, which provides 
the component $\int \mu_x:$
\[\begin{tikzcd}
	1 & {\int Px} & Bx \\
	1 & {\int P'x} & Bx
	\arrow[equals, from=1-1, to=2-1]
	\arrow[from=1-2, to=1-1]
	\arrow["{\pi_x}", from=1-2, to=1-3]
	\arrow["{\mu_x}", from=1-2, to=2-2]
	\arrow[equals, from=1-3, to=2-3]
	\arrow[from=2-2, to=2-1]
	\arrow["{\pi'_x}", from=2-2, to=2-3]
\end{tikzcd}\]
The cell component of $\int \mu$ at $m:x\proto y$ in $\dbl{D}$ 
is uniquely determined by the universal property of the pullback 
$\int P'm$ as below:
\[\begin{tikzcd}
	{\int Px} & {\int Pm} & {\int Py} \\
	{\int P'x} & {\int P'm} & {\int P'y} \\
	Bx & Bm
	\arrow[from=1-1, to=2-1]
	\arrow[from=1-2, to=1-1]
	\arrow[from=1-2, to=1-3]
	\arrow[dashed, from=1-2, to=2-2]
	\arrow[from=1-3, to=2-3]
	\arrow[from=2-1, to=3-1]
	\arrow[from=2-2, to=2-1]
	\arrow[from=2-2, to=2-3]
	\arrow["\lrcorner"{anchor=center, pos=0.125, rotate=-90}, draw=none, from=2-2, to=3-1]
	\arrow[from=2-2, to=3-2]
	\arrow[from=3-2, to=3-1]
\end{tikzcd}\]
In particular, as we mentioned above, morphisms of discrete opfibrations are
uniquely determined by their object components; since $\int$ and $\nabla$ 
are by definition mutually inverse on object components, we see they will 
be fully faithful. It remains to check that $\int \mu$ is a natural 
transformation.
\begin{itemize}
\item (Naturality with respect to cells) Given a cell
$\xinlinecell{x}{y}{z}{w}{f}{g}{m}{n}{\alpha}$ in $\dbl{D},$ we must 
show the following equation:
\[\begin{tikzcd}
	{\int P x} & {\int Pm} & {\int P y} && {\int P x} & {\int Pm} & {\int Py} \\
	{\int P z} & {\int P n} & {\int P w} & {=} & {\int P' x} & {\int P'm} & {\int P'y} \\
	{\int P ' z} & {\int P' n} & {\int P' w} && {\int P' z} & {\int P'n} & {\int P'w}
	\arrow["{\int Pf}"', from=1-1, to=2-1]
	\arrow[from=1-2, to=1-1]
	\arrow[from=1-2, to=1-3]
	\arrow["{\int P \alpha}"', from=1-2, to=2-2]
	\arrow["{\int P g}"', from=1-3, to=2-3]
	\arrow["{\int \mu_x}"', from=1-5, to=2-5]
	\arrow[from=1-6, to=1-5]
	\arrow[from=1-6, to=1-7]
	\arrow["{\int \mu_m}"', from=1-6, to=2-6]
	\arrow["{\int \mu_y}"', from=1-7, to=2-7]
	\arrow["{\int \mu_z}"', from=2-1, to=3-1]
	\arrow[from=2-2, to=2-1]
	\arrow[from=2-2, to=2-3]
	\arrow["{\int\mu_n}"', from=2-2, to=3-2]
	\arrow["{\int \mu_w}"', from=2-3, to=3-3]
	\arrow["{\int P'f}"', from=2-5, to=3-5]
	\arrow[from=2-6, to=2-5]
	\arrow[from=2-6, to=2-7]
	\arrow["{\int P'\alpha}"', from=2-6, to=3-6]
	\arrow["{\int P'g}"', from=2-7, to=3-7]
	\arrow[from=3-2, to=3-1]
	\arrow[from=3-2, to=3-3]
	\arrow[from=3-6, to=3-5]
	\arrow[from=3-6, to=3-7]
\end{tikzcd}\]
Now, $\int P\alpha:\int Px\times_{Bx} Bm\to \int Pz\times_{Bz} Bn$ is defined
as $Pf\times_{Bf} B\alpha$ and 
$\int\mu_n:\int Pz\times_{Bz} Bn\to \int P'n\times_{Bz} Bn$ is defined as
$\mu_z\times_{Bz} Bn,$ so that the left-hand pasting is 
$(Pf\cdot \mu_z)\times_{Bf} B\alpha$ while, similarly, the right-hand
pasting is $(\mu_x\cdot P'f)\times_{Bf}B\alpha.$ These coincide 
by naturality of $\mu.$
\item (External functoriality)
For unitality, we must check 
\[\begin{tikzcd}[column sep=tiny]
	{\int Px} & {\int Px} & {\int Px} && {\int Px} & {\int Px} & {\int Px} \\
	{\int Px} & {\int P(\id_x)} & {\int Px} & {=} & {\int P' x} & {\int P' x} & {\int P' x} \\
	{\int P'x} & {\int P'(\id_x)} & {\int P'x} && {\int P'x} & {\int P'\id_x} & {\int P'x}
	\arrow[equals, from=1-1, to=2-1]
	\arrow[from=1-2, to=1-1]
	\arrow[from=1-2, to=1-3]
	\arrow["{\int P_x}"', from=1-2, to=2-2]
	\arrow[equals, from=1-3, to=2-3]
	\arrow["{\mu_x}"', from=1-5, to=2-5]
	\arrow[from=1-6, to=1-5]
	\arrow[from=1-6, to=1-7]
	\arrow["{\id_{\mu_x}}"', from=1-6, to=2-6]
	\arrow["{\mu_x}"', from=1-7, to=2-7]
	\arrow["{\int\mu_x}"', from=2-1, to=3-1]
	\arrow[from=2-2, to=2-1]
	\arrow[from=2-2, to=2-3]
	\arrow["{\mu_{\id_x}}"', from=2-2, to=3-2]
	\arrow["{\mu_x}"', from=2-3, to=3-3]
	\arrow[equals, from=2-5, to=3-5]
	\arrow[from=2-6, to=2-5]
	\arrow[from=2-6, to=2-7]
	\arrow["{\int P'_x}"', from=2-6, to=3-6]
	\arrow[equals, from=2-7, to=3-7]
	\arrow[from=3-2, to=3-1]
	\arrow[from=3-2, to=3-3]
	\arrow[from=3-6, to=3-5]
	\arrow[from=3-6, to=3-7]
\end{tikzcd}\]
The left-hand side is the composite 
\[\int Px \xto{(\id,\pi_x\cdot B_x)} \int P x\times_{Bx} B\id_x
\xto{\int\mu_x\times_{Bx}B\id_x}\int P'x\times_{Bx} B\id_x\]
which composes to $(\int \mu_x,\pi_x\cdot B_x).$

The right-hand side is the composite 
\[\int Px \xto{\int \mu_x} \int P'x \xto{(\id_{P'x},\pi'_x\cdot B_x)} \int P'x\times_{Bx}B\id_x\]
which composes to $(\int \mu_x,\int\mu_x\pi'_x\cdot B_x).$ 
So this follows from the fact that $\int\mu_x\cdot \pi'_x=\pi_x.$

Finally, we check respect for composition. Given 
$x\xproto{m}y\xproto{n}z,$ we must show that the following 
cells are equal:
\[\begin{tikzcd}
	{\int P x} & {\int Pm} & {\int Py} & {\int Pn} & {\int Pz} \\
	&& \bullet \\
	{\int Px} && {\int P(m\odot n)} && {\int Pz} \\
	{\int P'x} && {\int P'(m\odot n)} && {\int P' z} \\
	{\int Px} & {\int Pm} & {\int Py} & {\int Pn} & {\int Pz} \\
	{\int P' x} & {\int P' m} & {\int P' y} & {\int P' n} & {\int P' z} \\
	&& \bullet \\
	{\int P' x} && {\int P' (m\odot n)} && {\int P' z}
	\arrow[equals, from=1-1, to=3-1]
	\arrow[from=1-2, to=1-1]
	\arrow[from=1-2, to=1-3]
	\arrow[from=1-4, to=1-3]
	\arrow[from=1-4, to=1-5]
	\arrow[equals, from=1-5, to=3-5]
	\arrow[from=2-3, to=1-2]
	\arrow["\lrcorner"{anchor=center, pos=0.125, rotate=135}, draw=none, from=2-3, to=1-3]
	\arrow[from=2-3, to=1-4]
	\arrow["{\int P_{m,n}}"', from=2-3, to=3-3]
	\arrow["{\int\mu_x}"', from=3-1, to=4-1]
	\arrow[from=3-3, to=3-1]
	\arrow[from=3-3, to=3-5]
	\arrow["{\int\mu_{m\odot n}}"', from=3-3, to=4-3]
	\arrow["{\int\mu_z}"', from=3-5, to=4-5]
	\arrow[from=4-3, to=4-1]
	\arrow[from=4-3, to=4-5]
	\arrow["{\int\mu_x}"', from=5-1, to=6-1]
	\arrow[from=5-2, to=5-1]
	\arrow[from=5-2, to=5-3]
	\arrow["{\int\mu_m}"', from=5-2, to=6-2]
	\arrow["{\int\mu_y}"', from=5-3, to=6-3]
	\arrow[from=5-4, to=5-3]
	\arrow[from=5-4, to=5-5]
	\arrow["{\int\mu_n}"', from=5-4, to=6-4]
	\arrow["{\int\mu_z}"', from=5-5, to=6-5]
	\arrow[equals, from=6-1, to=8-1]
	\arrow[from=6-2, to=6-1]
	\arrow[from=6-2, to=6-3]
	\arrow[from=6-4, to=6-3]
	\arrow[from=6-4, to=6-5]
	\arrow[equals, from=6-5, to=8-5]
	\arrow[from=7-3, to=6-2]
	\arrow["\lrcorner"{anchor=center, pos=0.125, rotate=135}, draw=none, from=7-3, to=6-3]
	\arrow[from=7-3, to=6-4]
	\arrow["{\int P' _{m,n}}"', from=7-3, to=8-3]
	\arrow[from=8-3, to=8-1]
	\arrow[from=8-3, to=8-5]
\end{tikzcd}\]
The domains are 
$\int Px\times_{Bx} Bm\times_{By} Bn$
Recall that $\int P_{m,n},$ with codomain 
$\int Px\times_{Bx} B(m\odot n),$ is defined as
$\int Px\times_{Bx} B_{m,n},$ while 
$\mu_{m\odot n}:\int Px\times_{Bx} B(m\odot n)
\to \int P'x\times_{Bx} B(m\odot n)$ is defined as
$\mu_x\times_{Bx} B(m\odot n).$ Therefore the upper
composite is given by $\mu_x\times_{Bx} B_{m,n}.$

On the other hand, the lower composite is given by 
$\mu_x\times_{Bx}Bm\times_{By} Bn$ composed
with $\int P'x\times_{Bx} B_{m,n},$
which is also $\mu_x\times_{Bx} B_{m,n}.$
\end{itemize}

\item (Essential surjectivity) 
Given a choice of pullback functor in $\Set,$ 
$\nabla$ and $\int$ are actually bijective on objects.
This concludes the proof.
\end{enumerate}
\end{proof}

\end{appendices}

\end{document}